
\ifx\shlhetal\undefinedcontrolsequence\let\shlhetal\relax\fi

\input amstex
\expandafter\ifx\csname mathdefs.tex\endcsname\relax
  \expandafter\gdef\csname mathdefs.tex\endcsname{}
\else \message{Hey!  Apparently you were trying to
  \string\input{mathdefs.tex} twice.   This does not make sense.} 
\errmessage{Please edit your file (probably \jobname.tex) and remove
any duplicate ``\string\input'' lines}\endinput\fi




\catcode`\X=12\catcode`\@=11

\def\n@wcount{\alloc@0\count\countdef\insc@unt}
\def\n@wwrite{\alloc@7\write\chardef\sixt@@n}
\def\n@wread{\alloc@6\read\chardef\sixt@@n}
\def\r@s@t{\relax}\def\v@idline{\par}\def\@mputate#1/{#1}
\def\l@c@l#1X{\firstpart.#1}\def\gl@b@l#1X{#1}\def\t@d@l#1X{{}}

\def\crossrefs#1{\ifx\all#1\let\tr@ce=\all\else\def\tr@ce{#1,}\fi
   \n@wwrite\cit@tionsout\openout\cit@tionsout=\jobname.cit 
   \write\cit@tionsout{\tr@ce}\expandafter\setfl@gs\tr@ce,}
\def\setfl@gs#1,{\def\@{#1}\ifx\@\empty\let\next=\relax
   \else\let\next=\setfl@gs\expandafter\xdef
   \csname#1tr@cetrue\endcsname{}\fi\next}
\def\m@ketag#1#2{\expandafter\n@wcount\csname#2tagno\endcsname
     \csname#2tagno\endcsname=0\let\tail=\all\xdef\all{\tail#2,}
   \ifx#1\l@c@l\let\tail=\r@s@t\xdef\r@s@t{\csname#2tagno\endcsname=0\tail}\fi
   \expandafter\gdef\csname#2cite\endcsname##1{\expandafter
     \ifx\csname#2tag##1\endcsname\relax?\else\csname#2tag##1\endcsname\fi
     \expandafter\ifx\csname#2tr@cetrue\endcsname\relax\else
     \write\cit@tionsout{#2tag ##1 cited on page \folio.}\fi}
   \expandafter\gdef\csname#2page\endcsname##1{\expandafter
     \ifx\csname#2page##1\endcsname\relax?\else\csname#2page##1\endcsname\fi
     \expandafter\ifx\csname#2tr@cetrue\endcsname\relax\else
     \write\cit@tionsout{#2tag ##1 cited on page \folio.}\fi}
   \expandafter\gdef\csname#2tag\endcsname##1{\expandafter
      \ifx\csname#2check##1\endcsname\relax
      \expandafter\xdef\csname#2check##1\endcsname{}%
      \else\immediate\write16{Warning: #2tag ##1 used more than once.}\fi
      \multit@g{#1}{#2}##1/X%
      \write\t@gsout{#2tag ##1 assigned number \csname#2tag##1\endcsname\space
      on page \number\count0.}%
   \csname#2tag##1\endcsname}}

\def\multit@g#1#2#3/#4X{\def\t@mp{#4}\ifx\t@mp\empty%
      \global\advance\csname#2tagno\endcsname by 1 
      \expandafter\xdef\csname#2tag#3\endcsname
      {#1\number\csname#2tagno\endcsnameX}%
   \else\expandafter\ifx\csname#2last#3\endcsname\relax
      \expandafter\n@wcount\csname#2last#3\endcsname
      \global\advance\csname#2tagno\endcsname by 1 
      \expandafter\xdef\csname#2tag#3\endcsname
      {#1\number\csname#2tagno\endcsnameX}
      \write\t@gsout{#2tag #3 assigned number \csname#2tag#3\endcsname\space
      on page \number\count0.}\fi
   \global\advance\csname#2last#3\endcsname by 1
   \def\t@mp{\expandafter\xdef\csname#2tag#3/}%
   \expandafter\t@mp\@mputate#4\endcsname
   {\csname#2tag#3\endcsname\lastpart{\csname#2last#3\endcsname}}\fi}
\def\t@gs#1{\def\all{}\m@ketag#1e\m@ketag#1s\m@ketag\t@d@l p
\let\realscite\scite
\let\realstag\stag
   \m@ketag\gl@b@l r \n@wread\t@gsin
   \openin\t@gsin=\jobname.tgs \re@der \closein\t@gsin
   \n@wwrite\t@gsout\openout\t@gsout=\jobname.tgs }
\outer\def\localtags{\t@gs\l@c@l}
\outer\def\globaltags{\t@gs\gl@b@l}
\outer\def\newlocaltag#1{\m@ketag\l@c@l{#1}}
\outer\def\newglobaltag#1{\m@ketag\gl@b@l{#1}}

\newif\ifpr@ 
\def\m@kecs #1tag #2 assigned number #3 on page #4.%
   {\expandafter\gdef\csname#1tag#2\endcsname{#3}
   \expandafter\gdef\csname#1page#2\endcsname{#4}
   \ifpr@\expandafter\xdef\csname#1check#2\endcsname{}\fi}
\def\re@der{\ifeof\t@gsin\let\next=\relax\else
   \read\t@gsin to\t@gline\ifx\t@gline\v@idline\else
   \expandafter\m@kecs \t@gline\fi\let \next=\re@der\fi\next}
\def\pretags#1{\pr@true\pret@gs#1,,}
\def\pret@gs#1,{\def\@{#1}\ifx\@\empty\let\n@xtfile=\relax
   \else\let\n@xtfile=\pret@gs \openin\t@gsin=#1.tgs \message{#1} \re@der 
   \closein\t@gsin\fi \n@xtfile}

\newcount\sectno\sectno=0\newcount\subsectno\subsectno=0
\newif\ifultr@local \def\ultralocal{\ultr@localtrue}
\def\firstpart{\number\sectno}
\def\lastpart#1{\ifcase#1 \or a\or b\or c\or d\or e\or f\or g\or h\or 
   i\or k\or l\or m\or n\or o\or p\or q\or r\or s\or t\or u\or v\or w\or 
   x\or y\or z \fi}

\def\resetall{\global\advance\sectno by 1\subsectno=0
   \gdef\firstpart{\number\sectno}\r@s@t}
\def\resetsub{\global\advance\subsectno by 1
   \gdef\firstpart{\number\sectno.\number\subsectno}\r@s@t}
\def\newsection#1\par{\resetall\vskip0pt plus.3\vsize\penalty-250
   \vskip0pt plus-.3\vsize\bigskip\bigskip
   \message{#1}\leftline{\bf#1}\nobreak\bigskip}
\def\subsection#1\par{\ifultr@local\resetsub\fi
   \vskip0pt plus.2\vsize\penalty-250\vskip0pt plus-.2\vsize
   \bigskip\smallskip\message{#1}\leftline{\bf#1}\nobreak\medskip}


\newdimen\marginshift

\newdimen\margindelta
\newdimen\marginmax
\newdimen\marginmin

\def\margininit{       
\marginmax=3 true cm                  
				      
\margindelta=0.1 true cm              
\marginmin=0.1true cm                 
\marginshift=\marginmin
}    

\def\t@gsjj#1,{\def\@{#1}\ifx\@\empty\let\next=\relax\else\let\next=\t@gsjj
   \def\@@{p}\ifx\@\@@\else
   \expandafter\gdef\csname#1cite\endcsname##1{\citejj{##1}}
   \expandafter\gdef\csname#1page\endcsname##1{?}
   \expandafter\gdef\csname#1tag\endcsname##1{\tagjj{##1}}\fi\fi\next}
\newif\ifshowstuffinmargin
\showstuffinmarginfalse
\def\jjtags{\ifx\shlhetal\relax 
  \else
\ifx\shlhetal\undefinedcontrolseq
\else
\showstuffinmargintrue
\ifx\all\relax\else\expandafter\t@gsjj\all,\fi\fi \fi
}

\def\tagjj#1{\realstag{#1}\oldmginpar{\zeigen{#1}}}
\def\citejj#1{\rechnen{#1}\mginpar{\zeigen{#1}}}     

\def\rechnen#1{\expandafter\ifx\csname stag#1\endcsname\relax ??\else
                           \csname stag#1\endcsname\fi}

\newdimen\theight

\def\marginfont{\sevenrm}

\def\trymarginbox#1{\setbox0=\hbox{\marginfont\hskip\marginshift #1}%
		\global\marginshift\wd0 
		\global\advance\marginshift\margindelta}

\def \oldmginpar#1{%
\ifvmode\setbox0\hbox to \hsize{\hfill\rlap{\marginfont\quad#1}}%
\ht0 0cm
\dp0 0cm
\box0\vskip-\baselineskip
\else 
             \vadjust{\trymarginbox{#1}%
		\ifdim\marginshift>\marginmax \global\marginshift\marginmin
			\trymarginbox{#1}%
                \fi
             \theight=\ht0
             \advance\theight by \dp0    \advance\theight by \lineskip
             \kern -\theight \vbox to \theight{\rightline{\rlap{\box0}}%
\vss}}\fi}

\newdimen\upordown
\global\upordown=8pt
\font\tinyfont=cmtt8 
\def\mginpar#1{\smash{\hbox to 0cm{\kern-10pt\raise7pt\hbox{\tinyfont #1}\hss}}}
\def\mginpar#1{{\hbox to 0cm{\kern-10pt\raise\upordown\hbox{\tinyfont #1}\hss}}\global\upordown-\upordown}


\def\t@gsoff#1,{\def\@{#1}\ifx\@\empty\let\next=\relax\else\let\next=\t@gsoff
   \def\@@{p}\ifx\@\@@\else
   \expandafter\gdef\csname#1cite\endcsname##1{\zeigen{##1}}
   \expandafter\gdef\csname#1page\endcsname##1{?}
   \expandafter\gdef\csname#1tag\endcsname##1{\zeigen{##1}}\fi\fi\next}
\def\verbatimtags{\showstuffinmarginfalse
\ifx\all\relax\else\expandafter\t@gsoff\all,\fi}
\def\zeigen#1{\hbox{$\scriptstyle\langle$}#1\hbox{$\scriptstyle\rangle$}}


\def\margintag#1{\ifshowstuffinmargin\oldmginpar{\zeigen{#1}}\fi}

\def\marginplain#1{\ifshowstuffinmargin\mginpar{{#1}}\fi}
\def\marginbf#1{\marginplain{{\bf \ \ #1}}}

\def\(#1){\edef\dot@g{\ifmmode\ifinner(\hbox{\noexpand\etag{#1}})
   \else\noexpand\eqno(\hbox{\noexpand\etag{#1}})\fi
   \else(\noexpand\ecite{#1})\fi}\dot@g}

\newif\ifbr@ck
\def\eat#1{}
\def\[#1]{\br@cktrue[\br@cket#1'X]}
\def\br@cket#1'#2X{\def\temp{#2}\ifx\temp\empty\let\next\eat
   \else\let\next\br@cket\fi
   \ifbr@ck\br@ckfalse\br@ck@t#1,X\else\br@cktrue#1\fi\next#2X}
\def\br@ck@t#1,#2X{\def\temp{#2}\ifx\temp\empty\let\neext\eat
   \else\let\neext\br@ck@t\def\temp{,}\fi
   \def\teemp{#1}\ifx\teemp\empty\else\rcite{#1}\fi\temp\neext#2X}
\def\resetbr@cket{\gdef\[##1]{[\rtag{##1}]}}
\def\references{\resetbr@cket\newsection References\par}

\newtoks\symb@ls\newtoks\s@mb@ls\newtoks\p@gelist\n@wcount\ftn@mber
    \ftn@mber=1\newif\ifftn@mbers\ftn@mbersfalse\newif\ifbyp@ge\byp@gefalse
\def\defm@rk{\ifftn@mbers\n@mberm@rk\else\symb@lm@rk\fi}
\def\n@mberm@rk{\xdef\m@rk{{\the\ftn@mber}}%
    \global\advance\ftn@mber by 1 }
\def\rot@te#1{\let\temp=#1\global#1=\expandafter\r@t@te\the\temp,X}
\def\r@t@te#1,#2X{{#2#1}\xdef\m@rk{{#1}}}
\def\b@@st#1{{$^{#1}$}}\def\str@p#1{#1}
\def\symb@lm@rk{\ifbyp@ge\rot@te\p@gelist\ifnum\expandafter\str@p\m@rk=1 
    \s@mb@ls=\symb@ls\fi\write\f@nsout{\number\count0}\fi \rot@te\s@mb@ls}
\def\byp@ge{\byp@getrue\n@wwrite\f@nsin\openin\f@nsin=\jobname.fns 
    \n@wcount\currentp@ge\currentp@ge=0\p@gelist={0}
    \re@dfns\closein\f@nsin\rot@te\p@gelist
    \n@wread\f@nsout\openout\f@nsout=\jobname.fns }
\def\m@kelist#1X#2{{#1,#2}}
\def\re@dfns{\ifeof\f@nsin\let\next=\relax\else\read\f@nsin to \f@nline
    \ifx\f@nline\v@idline\else\let\t@mplist=\p@gelist
    \ifnum\currentp@ge=\f@nline
    \global\p@gelist=\expandafter\m@kelist\the\t@mplistX0
    \else\currentp@ge=\f@nline
    \global\p@gelist=\expandafter\m@kelist\the\t@mplistX1\fi\fi
    \let\next=\re@dfns\fi\next}
\def\symbols#1{\symb@ls={#1}\s@mb@ls=\symb@ls} 
\def\bigsymbol{\textstyle}
\symbols{\bigsymbol\ast,\dagger,\ddagger,\sharp,\flat,\natural,\star}
\def\ftnumbers{\ftn@mberstrue} \def\ftsymbols{\ftn@mbersfalse}
\def\paginal{\byp@ge} \def\resetftnumbers{\ftn@mber=1}
\def\ftnote#1{\defm@rk\expandafter\expandafter\expandafter\footnote
    \expandafter\b@@st\m@rk{#1}}

\long\def\jump#1\endjump{}
\def\ssum{\mathop{\lower .1em\hbox{$\textstyle\Sigma$}}\nolimits}

\def\qed{\nobreak\kern 1em \vrule height .5em width .5em depth 0em}
\def\newneq{\hbox{\rlap{\hbox to 1\wd9{\hss$=$\hss}}\raise .1em 
   \hbox to 1\wd9{\hss$\scriptscriptstyle/$\hss}}}
\def\subsetne{\setbox9 = \hbox{$\subset$}\mathrel{\hbox{\rlap
   {\lower .4em \newneq}\raise .13em \hbox{$\subset$}}}}
\def\supsetne{\setbox9 = \hbox{$\subset$}\mathrel{\hbox{\rlap
   {\lower .4em \newneq}\raise .13em \hbox{$\supset$}}}}

\def\vbar{\mathchoice{\vrule height6.3ptdepth-.5ptwidth.8pt\kern-.8pt}
   {\vrule height6.3ptdepth-.5ptwidth.8pt\kern-.8pt}
   {\vrule height4.1ptdepth-.35ptwidth.6pt\kern-.6pt}
   {\vrule height3.1ptdepth-.25ptwidth.5pt\kern-.5pt}}
\def\f@dge{\mathchoice{}{}{\mkern.5mu}{\mkern.8mu}}
\def\b@c#1#2{{\rm \mkern#2mu\vbar\mkern-#2mu#1}}
\def\b@b#1{{\rm I\mkern-3.5mu #1}}
\def\b@a#1#2{{\rm #1\mkern-#2mu\f@dge #1}}
\def\bb#1{{\count4=`#1 \advance\count4by-64 \ifcase\count4\or\b@a A{11.5}\or
   \b@b B\or\b@c C{5}\or\b@b D\or\b@b E\or\b@b F \or\b@c G{5}\or\b@b H\or
   \b@b I\or\b@c J{3}\or\b@b K\or\b@b L \or\b@b M\or\b@b N\or\b@c O{5} \or
   \b@b P\or\b@c Q{5}\or\b@b R\or\b@a S{8}\or\b@a T{10.5}\or\b@c U{5}\or
   \b@a V{12}\or\b@a W{16.5}\or\b@a X{11}\or\b@a Y{11.7}\or\b@a Z{7.5}\fi}}

\catcode`\X=11 \catcode`\@=12




\let\thischap\jobname

\def\partof#1{\csname returnthe#1part\endcsname}
\def\CHAPOF#1{\csname returnthe#1chap\endcsname}

\def\chapof#1{\CHAPOF{#1}}

\def\setchapter#1,#2,#3;{%
  \expandafter\def\csname returnthe#1part\endcsname{#2}%
  \expandafter\def\csname returnthe#1chap\endcsname{#3}%
}

\def\setprevious#1 #2 {%
  \expandafter\def\csname set#1page\endcsname{\input page-#2}
}


 \setchapter  E53,B,N;       \setprevious E53 null
 \setchapter  300z,B,A;       \setprevious 300z E53
 \setchapter  88r,B,I;       \setprevious 88r 300z
 \setchapter  600,B,II;       \setprevious  600 88r
 \setchapter  705,B,III;       \setprevious   705 600
 \setchapter  734,B,IV;        \setprevious   734 705
 \setchapter  300x,B,;      \setprevious   300x 734

 \setchapter 300a,A,V.A;      \setprevious 300a 88r
 \setchapter 300b,A,V.B;       \setprevious 300b 300a
 \setchapter 300c,A,V.C;       \setprevious 300c 300b
 \setchapter 300d,A,V.D;       \setprevious 300d 300c
 \setchapter 300e,A,V.E;       \setprevious 300e 300d
 \setchapter 300f,A,V.F;       \setprevious 300f 300e
 \setchapter 300g,A,V.G;       \setprevious 300g 300f

  \setchapter  E46,B,VI;      \setprevious    E46 734
  \setchapter  838,B,VII;      \setprevious   838 E46

\newwrite\pageout
\def\rememberpagenumber{\let\setpage\relax
\openout\pageout page-\jobname  \relax \write\pageout{\setpage\the\pageno.}}

\def\recallpagenumber{\csname set\jobname page\endcsname
\def\headmark##1{\rightheadtext{\chapof{\jobname}.##1}}\WRITETOC}
\def\setupchapter#1{\leftheadtext{\chapof{\jobname}. #1}}

\def\setpage#1.{\pageno=#1\relax\advance\pageno1\relax}

\def\cprefix#1{
\edef\theotherpart{\partof{#1}}\edef\theotherchap{\chapof{#1}}%
\ifx\theotherpart\thispart
   \ifx\theotherchap\thischap 
    \else 
     \theotherchap%
    \fi
   \else 
     \theotherchap\fi}

\def\sectioncite[#1]#2{%
     \cprefix{#2}#1}

\def\chaptercite#1{Chapter \cprefix{#1}}

\edef\thispart{\partof{\thischap}}
\edef\thischap{\chapof{\thischap}}

\def\lastpage of '#1' is #2.{\expandafter\def\csname lastpage#1\endcsname{#2}}

\def\yCITE[#1]#2{\cprefix{#2}.\scite{#2-#1}}

\newwrite\writetoc
\immediate\openout\writetoc \jobname.toc
\def\addcontents#1{\def\WRITETOC{\immediate\write\writetoc{\noexpand\tocentry{\chapof{\jobname}}{#1}{\number\pageno}}}}



\def\spuriousreset{}


\expandafter\ifx\csname citeadd.tex\endcsname\relax
\expandafter\gdef\csname citeadd.tex\endcsname{}
\else \message{Hey!  Apparently you were trying to
\string\input{citeadd.tex} twice.   This does not make sense.} 
\errmessage{Please edit your file (probably \jobname.tex) and remove
any duplicate ``\string\input'' lines}\endinput\fi

\sectno=-1   
\localtags
\jjtags
\NoBlackBoxes
\define\mr{\medskip\roster}
\define\sn{\smallskip\noindent}
\define\mn{\medskip\noindent}
\define\bn{\bigskip\noindent}
\define\ub{\underbar}
\define\wilog{\text{without loss of generality}}
\define\ermn{\endroster\medskip\noindent}

\define\ortp{\text{\bf tp}}
\define\sftp{\text{\rm tp}}
\define\aec{\text{abstract elementary class \,}}
\define \nl{\newline}
\newbox\noforkbox \newdimen\forklinewidth
\forklinewidth=0.3pt   
\setbox0\hbox{$\textstyle\bigcup$}
\setbox1\hbox to \wd0{\hfil\vrule width \forklinewidth depth \dp0
                        height \ht0 \hfil}
\wd1=0 cm
\setbox\noforkbox\hbox{\box1\box0\relax}
\def\unionstick{\mathop{\copy\noforkbox}\limits}
 \def\nonfork#1#2_#3{#1\unionstick_{\textstyle #3}#2}
 \def\nonforkin#1#2_#3^#4{#1\unionstick_{\textstyle #3}^{\textstyle #4}#2}     
%
\setbox0\hbox{$\textstyle\bigcup$}
\setbox1\hbox to \wd0{\hfil{\sl /\/}\hfil}
\setbox2\hbox to \wd0{\hfil\vrule height \ht0 depth \dp0 width
                                \forklinewidth\hfil}

\wd1=0cm
\wd2=0cm
\newbox\doesforkbox
\setbox\doesforkbox\hbox{\box1\box0\relax}
\def\nunionstick{\mathop{\copy\doesforkbox}\limits}

\def\fork#1#2_#3{#1\nunionstick_{\textstyle #3}#2}
\def\forkin#1#2_#3^#4{#1\nunionstick_{\textstyle #3}^{\textstyle #4}#2}     

\magnification=\magstep 1
\documentstyle{amsppt}

{    
\catcode`@11

\ifx\alicetwothousandloaded@\relax
  \endinput\else\global\let\alicetwothousandloaded@\relax\fi

\gdef\subjclass{\let\savedef@\subjclass
 \def\subjclass##1\endsubjclass{\let\subjclass\savedef@
   \toks@{\def\usualspace{{\rm\enspace}}\eightpoint}%
   \toks@@{##1\unskip.}%
   \edef\thesubjclass@{\the\toks@
     \frills@{{\noexpand\rm2000 {\noexpand\it Mathematics Subject
       Classification}.\noexpand\enspace}}%
     \the\toks@@}}%
  \nofrillscheck\subjclass}
} 


\expandafter\ifx\csname alice2jlem.tex\endcsname\relax
  \expandafter\xdef\csname alice2jlem.tex\endcsname{\the\catcode`@}
\else \message{Hey!  Apparently you were trying to
\string\input{alice2jlem.tex}  twice.   This does not make sense.}
\errmessage{Please edit your file (probably \jobname.tex) and remove
any duplicate ``\string\input'' lines}\endinput\fi

\expandafter\ifx\csname bib4plain.tex\endcsname\relax
  \expandafter\gdef\csname bib4plain.tex\endcsname{}
\else \message{Hey!  Apparently you were trying to \string\input
  bib4plain.tex twice.   This does not make sense.}
\errmessage{Please edit your file (probably \jobname.tex) and remove
any duplicate ``\string\input'' lines}\endinput\fi

\def\renewcommand{\newcommand}	       
\edef\cite{\the\catcode`@}%
\catcode`@ = 11
\let\@oldatcatcode = \cite
\chardef\@letter = 11
\chardef\@other = 12
%
%
%
%
\def\@innerdef#1#2{\edef#1{\expandafter\noexpand\csname #2\endcsname}}%
%
%
\@innerdef\@innernewcount{newcount}%
\@innerdef\@innernewdimen{newdimen}%
\@innerdef\@innernewif{newif}%
\@innerdef\@innernewwrite{newwrite}%
%
%
%
\def\@gobble#1{}%
%
%
%
\ifx\inputlineno\@undefined
   \let\@linenumber = \empty 
\else
   \def\@linenumber{\the\inputlineno:\space}%
\fi
%
%
%
\def\@futurenonspacelet#1{\def\cs{#1}%
   \afterassignment\@stepone\let\@nexttoken=
}%
\begingroup 
\def\\{\global\let\@stoken= }%
\\ 
\endgroup
\def\@stepone{\expandafter\futurelet\cs\@steptwo}%
\def\@steptwo{\expandafter\ifx\cs\@stoken\let\@@next=\@stepthree
   \else\let\@@next=\@nexttoken\fi \@@next}%
\def\@stepthree{\afterassignment\@stepone\let\@@next= }%
%
%
%
\def\@getoptionalarg#1{%
   \let\@optionaltemp = #1%
   \let\@optionalnext = \relax
   \@futurenonspacelet\@optionalnext\@bracketcheck
}%
%
%
\def\@bracketcheck{%
   \ifx [\@optionalnext
      \expandafter\@@getoptionalarg
   \else
      \let\@optionalarg = \empty
      \expandafter\@optionaltemp
   \fi
}%
\def\@@getoptionalarg[#1]{%
   \def\@optionalarg{#1}%
   \@optionaltemp
}%
%
%
%
\def\@nnil{\@nil}%
\def\@fornoop#1\@@#2#3{}%
\def\@for#1:=#2\do#3{%
   \edef\@fortmp{#2}%
   \ifx\@fortmp\empty \else
      \expandafter\@forloop#2,\@nil,\@nil\@@#1{#3}%
   \fi
}%
\def\@forloop#1,#2,#3\@@#4#5{\def#4{#1}\ifx #4\@nnil \else
       #5\def#4{#2}\ifx #4\@nnil \else#5\@iforloop #3\@@#4{#5}\fi\fi
}%
\def\@iforloop#1,#2\@@#3#4{\def#3{#1}\ifx #3\@nnil
       \let\@nextwhile=\@fornoop \else
      #4\relax\let\@nextwhile=\@iforloop\fi\@nextwhile#2\@@#3{#4}%
}%
%
%
%
\@innernewif\if@fileexists
\def\@testfileexistence{\@getoptionalarg\@finishtestfileexistence}%
\def\@finishtestfileexistence#1{%
   \begingroup
      \def\extension{#1}%
      \immediate\openin0 =
         \ifx\@optionalarg\empty\jobname\else\@optionalarg\fi
         \ifx\extension\empty \else .#1\fi
         \space
      \ifeof 0
         \global\@fileexistsfalse
      \else
         \global\@fileexiststrue
      \fi
      \immediate\closein0
   \endgroup
}%
%
%
%
%
\def\bibliographystyle#1{%
   \@readauxfile
   \@writeaux{\string\bibstyle{#1}}%
}%
\let\bibstyle = \@gobble
%
%
\let\bblfilebasename = \jobname
\def\bibliography#1{%
   \@readauxfile
   \@writeaux{\string\bibdata{#1}}%
   \@testfileexistence[\bblfilebasename]{bbl}%
   \if@fileexists
      \nobreak
      \@readbblfile
   \fi
}%
\let\bibdata = \@gobble
%
%
\def\nocite#1{%
   \@readauxfile
   \@writeaux{\string\citation{#1}}%
}%
\@innernewif\if@notfirstcitation
%
%
\def\cite{\@getoptionalarg\@cite}%
%
%
\def\@cite#1{%
   \let\@citenotetext = \@optionalarg
   \printcitestart
   \nocite{#1}%
   \@notfirstcitationfalse
   \@for \@citation :=#1\do
   {%
      \expandafter\@onecitation\@citation\@@
   }%
   \ifx\empty\@citenotetext\else
      \printcitenote{\@citenotetext}%
   \fi
   \printcitefinish
}%
\newif\ifweareinprivate
\weareinprivatetrue
\ifx\shlhetal\undefinedcontrolseq\weareinprivatefalse\fi
\ifx\shlhetal\relax\weareinprivatefalse\fi
\def\@onecitation#1\@@{%
   \if@notfirstcitation
      \printbetweencitations
   \fi
   \expandafter \ifx \csname\@citelabel{#1}\endcsname \relax
      \if@citewarning
         \message{\@linenumber Undefined citation `#1'.}%
      \fi
     \ifweareinprivate
      \expandafter\gdef\csname\@citelabel{#1}\endcsname{%
\strut 
\vadjust{\vskip-\dp\strutbox
\vbox to 0pt{\vss\parindent0cm \leftskip=\hsize 
\advance\leftskip3mm
\advance\hsize 4cm\strut\openup-4pt 
\rightskip 0cm plus 1cm minus 0.5cm ?  #1 ?\strut}}
         {\tt
            \escapechar = -1
            \nobreak\hskip0pt\pfeilsw
            \expandafter\string\csname#1\endcsname
             \pfeilso
            \nobreak\hskip0pt
         }%
      }%
     \else  
      \expandafter\gdef\csname\@citelabel{#1}\endcsname{%
            {\tt\expandafter\string\csname#1\endcsname}
      }%
     \fi  
   \fi
   \csname\@citelabel{#1}\endcsname
   \@notfirstcitationtrue
}%
%
%
\def\@citelabel#1{b@#1}%
%
%
\def\@citedef#1#2{\expandafter\gdef\csname\@citelabel{#1}\endcsname{#2}}%
%
%
%
\def\@readbblfile{%
   \ifx\@itemnum\@undefined
      \@innernewcount\@itemnum
   \fi
   \begingroup
      \def\begin##1##2{%
         \setbox0 = \hbox{\biblabelcontents{##2}}%
         \biblabelwidth = \wd0
      }%
      \def\end##1{}
      %
      %
      \@itemnum = 0
      \def\bibitem{\@getoptionalarg\@bibitem}%
      \def\@bibitem{%
         \ifx\@optionalarg\empty
            \expandafter\@numberedbibitem
         \else
            \expandafter\@alphabibitem
         \fi
      }%
      \def\@alphabibitem##1{%
         \expandafter \xdef\csname\@citelabel{##1}\endcsname {\@optionalarg}%
         \ifx\biblabelprecontents\@undefined
            \let\biblabelprecontents = \relax
         \fi
         \ifx\biblabelpostcontents\@undefined
            \let\biblabelpostcontents = \hss
         \fi
         \@finishbibitem{##1}%
      }%
      \def\@numberedbibitem##1{%
         \advance\@itemnum by 1
         \expandafter \xdef\csname\@citelabel{##1}\endcsname{\number\@itemnum}%
         \ifx\biblabelprecontents\@undefined
            \let\biblabelprecontents = \hss
         \fi
         \ifx\biblabelpostcontents\@undefined
            \let\biblabelpostcontents = \relax
         \fi
         \@finishbibitem{##1}%
      }%
      \def\@finishbibitem##1{%
         \biblabelprint{\csname\@citelabel{##1}\endcsname}%
         \@writeaux{\string\@citedef{##1}{\csname\@citelabel{##1}\endcsname}}%
         \ignorespaces
      }%
      %
      %
      \let\em = \bblem
      \let\newblock = \bblnewblock
      \let\sc = \bblsc
      \frenchspacing
      \clubpenalty = 4000 \widowpenalty = 4000
      \tolerance = 10000 \hfuzz = .5pt
      \everypar = {\hangindent = \biblabelwidth
                      \advance\hangindent by \biblabelextraspace}%
      \bblrm
      \parskip = 1.5ex plus .5ex minus .5ex
      \biblabelextraspace = .5em
      \bblhook
      \input \bblfilebasename.bbl
   \endgroup
}%
%
%
\@innernewdimen\biblabelwidth
\@innernewdimen\biblabelextraspace
%
%
%
\def\biblabelprint#1{%
   \noindent
   \hbox to \biblabelwidth{%
      \biblabelprecontents
      \biblabelcontents{#1}%
      \biblabelpostcontents
   }%
   \kern\biblabelextraspace
}%
%
%
%
\def\biblabelcontents#1{{\bblrm [#1]}}%
%
%
\def\bblrm{\rm}%
%
%
\def\bblem{\it}%
%
%
\def\bblsc{\ifx\@scfont\@undefined
              \font\@scfont = cmcsc10
           \fi
           \@scfont
}%
%
%
\def\bblnewblock{\hskip .11em plus .33em minus .07em }%
%
%
\let\bblhook = \empty
%
%
%
\def\printcitestart{[}
\def\printcitefinish{]}
\def\printbetweencitations{, }
\def\printcitenote#1{, #1}
%
%
%
\let\citation = \@gobble
%
%
%
\@innernewcount\@numparams
%
%
\def\newcommand#1{%
   \def\@commandname{#1}%
   \@getoptionalarg\@continuenewcommand
}%
%
%
\def\@continuenewcommand{%
   \@numparams = \ifx\@optionalarg\empty 0\else\@optionalarg \fi \relax
   \@newcommand
}%
%
%
\def\@newcommand#1{%
   \def\@startdef{\expandafter\edef\@commandname}%
   \ifnum\@numparams=0
      \let\@paramdef = \empty
   \else
      \ifnum\@numparams>9
         \errmessage{\the\@numparams\space is too many parameters}%
      \else
         \ifnum\@numparams<0
            \errmessage{\the\@numparams\space is too few parameters}%
         \else
            \edef\@paramdef{%
               \ifcase\@numparams
                  \empty  No arguments.
               \or ####1%
               \or ####1####2%
               \or ####1####2####3%
               \or ####1####2####3####4%
               \or ####1####2####3####4####5%
               \or ####1####2####3####4####5####6%
               \or ####1####2####3####4####5####6####7%
               \or ####1####2####3####4####5####6####7####8%
               \or ####1####2####3####4####5####6####7####8####9%
               \fi
            }%
         \fi
      \fi
   \fi
   \expandafter\@startdef\@paramdef{#1}%
}%
%
%
%
%
\def\@readauxfile{%
   \if@auxfiledone \else 
      \global\@auxfiledonetrue
      \@testfileexistence{aux}%
      \if@fileexists
         \begingroup
            \endlinechar = -1
            \catcode`@ = 11
            \input \jobname.aux
         \endgroup
      \else
         \message{\@undefinedmessage}%
         \global\@citewarningfalse
      \fi
      \immediate\openout\@auxfile = \jobname.aux
   \fi
}%
%
%
\newif\if@auxfiledone
\ifx\noauxfile\@undefined \else \@auxfiledonetrue\fi
%
%
%
%
\@innernewwrite\@auxfile
\def\@writeaux#1{\ifx\noauxfile\@undefined \write\@auxfile{#1}\fi}%
%
%
%
\ifx\@undefinedmessage\@undefined
   \def\@undefinedmessage{No .aux file; I won't give you warnings about
                          undefined citations.}%
\fi
%
%
\@innernewif\if@citewarning
\ifx\noauxfile\@undefined \@citewarningtrue\fi
%
%
%
\catcode`@ = \@oldatcatcode

\def\pfeilso{\leavevmode
            \vrule width 1pt height9pt depth 0pt\relax
           \vrule width 1pt height8.7pt depth 0pt\relax
           \vrule width 1pt height8.3pt depth 0pt\relax
           \vrule width 1pt height8.0pt depth 0pt\relax
           \vrule width 1pt height7.7pt depth 0pt\relax
            \vrule width 1pt height7.3pt depth 0pt\relax
            \vrule width 1pt height7.0pt depth 0pt\relax
            \vrule width 1pt height6.7pt depth 0pt\relax
            \vrule width 1pt height6.3pt depth 0pt\relax
            \vrule width 1pt height6.0pt depth 0pt\relax
            \vrule width 1pt height5.7pt depth 0pt\relax
            \vrule width 1pt height5.3pt depth 0pt\relax
            \vrule width 1pt height5.0pt depth 0pt\relax
            \vrule width 1pt height4.7pt depth 0pt\relax
            \vrule width 1pt height4.3pt depth 0pt\relax
            \vrule width 1pt height4.0pt depth 0pt\relax
            \vrule width 1pt height3.7pt depth 0pt\relax
            \vrule width 1pt height3.3pt depth 0pt\relax
            \vrule width 1pt height3.0pt depth 0pt\relax
            \vrule width 1pt height2.7pt depth 0pt\relax
            \vrule width 1pt height2.3pt depth 0pt\relax
            \vrule width 1pt height2.0pt depth 0pt\relax
            \vrule width 1pt height1.7pt depth 0pt\relax
            \vrule width 1pt height1.3pt depth 0pt\relax
            \vrule width 1pt height1.0pt depth 0pt\relax
            \vrule width 1pt height0.7pt depth 0pt\relax
            \vrule width 1pt height0.3pt depth 0pt\relax}

\def\pfeilsw{ \leavevmode 
            \vrule width 1pt height0.3pt depth 0pt\relax
            \vrule width 1pt height0.7pt depth 0pt\relax
            \vrule width 1pt height1.0pt depth 0pt\relax
            \vrule width 1pt height1.3pt depth 0pt\relax
            \vrule width 1pt height1.7pt depth 0pt\relax
            \vrule width 1pt height2.0pt depth 0pt\relax
            \vrule width 1pt height2.3pt depth 0pt\relax
            \vrule width 1pt height2.7pt depth 0pt\relax
            \vrule width 1pt height3.0pt depth 0pt\relax
            \vrule width 1pt height3.3pt depth 0pt\relax
            \vrule width 1pt height3.7pt depth 0pt\relax
            \vrule width 1pt height4.0pt depth 0pt\relax
            \vrule width 1pt height4.3pt depth 0pt\relax
            \vrule width 1pt height4.7pt depth 0pt\relax
            \vrule width 1pt height5.0pt depth 0pt\relax
            \vrule width 1pt height5.3pt depth 0pt\relax
            \vrule width 1pt height5.7pt depth 0pt\relax
            \vrule width 1pt height6.0pt depth 0pt\relax
            \vrule width 1pt height6.3pt depth 0pt\relax
            \vrule width 1pt height6.7pt depth 0pt\relax
            \vrule width 1pt height7.0pt depth 0pt\relax
            \vrule width 1pt height7.3pt depth 0pt\relax
            \vrule width 1pt height7.7pt depth 0pt\relax
            \vrule width 1pt height8.0pt depth 0pt\relax
            \vrule width 1pt height8.3pt depth 0pt\relax
            \vrule width 1pt height8.7pt depth 0pt\relax
            \vrule width 1pt height9pt depth 0pt\relax
      }


\def\widestnumber#1#2{}

\def\citewarning#1{\ifx\shlhetal\relax 
    \else
    \par{#1}\par
    \fi
}

\def\rm{\fam0 \tenrm}

\def\fakesubhead#1\endsubhead{\bigskip\noindent{\bf#1}\par}



%
%
%

%

\font\textrsfs=rsfs10
\font\scriptrsfs=rsfs7
\font\scriptscriptrsfs=rsfs5

\newfam\rsfsfam
\textfont\rsfsfam=\textrsfs
\scriptfont\rsfsfam=\scriptrsfs
\scriptscriptfont\rsfsfam=\scriptscriptrsfs

\edef\oldcatcodeofat{\the\catcode`\@}
\catcode`\@11

\def\Cal@@#1{\noaccents@ \fam \rsfsfam #1}

\catcode`\@\oldcatcodeofat


\expandafter\ifx \csname margininit\endcsname \relax\else\margininit\fi

\long\def\red#1\endred{}
\long\def\green#1\endgreen{}
\long\def\blue#1\endblue{}
\long\def\private#1\endprivate{}

\def\endred{ \unmatched endred! }
\def\endgreen{ \unmatched endgreen! }
\def\endblue{ \unmatched endblue! }
\def\endprivate{ \unmatched endprivate! }

\ifx\latexcolors\undefinedcs\def\latexcolors{}\fi

\def\emptycs{}
\def\evaluatelatexcolors{%
        \ifx\latexcolors\emptycs\else
        \expandafter\xxevaluate\latexcolors\xxfertig\evaluatelatexcolors\fi}
\def\xxevaluate#1,#2\xxfertig{\setupthiscolor{#1}%
        \def\latexcolors{#2}}


\font\smallfont=cmsl7
\def\rutgerscolor{\ifmmode\else\endgraf\fi\smallfont
\advance\leftskip0.5cm\relax}
\def\setupthiscolor#1{\edef\tmptmpcs{\noexpand\bgroup\noexpand\rutgerscolor
\noexpand\def\noexpand\currentcolor{#1}%
\noexpand}%
\expandafter\let\csname#1\endcsname\tmptmpcs
\def\tmptmpcs{\checkColorUnmatched{#1}\popthecolor}
\expandafter\let\csname end#1\endcsname\tmptmpcs}

\def\checkColorUnmatched#1{\def\expectcolor{#1}%
    \ifx\expectcolor\currentcolor   
    \else \edef\failhere{\noexpand\tryingToClose '\currentcolor' with end\expectcolor}\failhere\fi}

\def\currentcolor{???}

\def\popthecolor{\ifmmode\else\endgraf\fi\egroup}

\expandafter\def\csname#1\endcsname{}

\evaluatelatexcolors

 \let\outerhead\head
 \def\head{\innerhead}
 \let\innerhead\outerhead

 \let\outersubhead\subhead
 \def\subhead{\innersubhead}
 \let\innersubhead\outersubhead

 \let\outersubsubhead\subsubhead
 \def\subsubhead{\innersubsubhead}
 \let\innersubsubhead\outersubsubhead

 \let\outerproclaim\proclaim
 \def\proclaim{\innerproclaim}
 \let\innerproclaim\outerproclaim

 %
 %
 %
 %

\def\demo#1{\medskip\noindent{\it #1.\/}}
\def\enddemo{\smallskip}

\def\beginaside{\endgraf\leftskip2cm \vrule width 0pt\relax}
\def\endaside{\endgraf\leftskip0cm \vrule width 0pt\relax}

\pageheight{8.5truein}
\topmatter
\title{Introduction to: classification theory for abstract elementary class} 
\endtitle
\rightheadtext{Intro to stability theory for a.e.c.} 
\author {Saharon Shelah \thanks {\null\newline 
I would like to thank
Alice Leonhardt for the beautiful typing. \null\newline
Publication E53} \endthanks} \endauthor 

\affil{The Hebrew University of Jerusalem \\
Einstein Institute of Mathematics \\
Edmond J. Safra Campus, Givat Ram \\
Jerusalem 91904, Israel
 \medskip
 Department of Mathematics \\
 Hill Center-Busch Campus \\
  Rutgers, The State University of New Jersey \\
 110 Frelinghuysen Road \\
 Piscataway, NJ 08854-8019 USA} \endaffil

\endtopmatter
\document  
 
\pretags{300a,300b,300c,300d,300e,300f,300g,300x,300y,300z,88r,600,705,734}
\newpage

\head {Content} \endhead
 \spuriousreset
\bn
Abstract, pg.3
\bn
\S0 Introduction, pg.4
\bn
\S1 Introduction for model theorists, pg.4-15
\sn
\hskip10pt (A) Why to be interested in dividing lines, pg.5
\sn
\hskip10pt (B) Historical comments on non-elementary classes, pg.10
\bn
\S2 Introduction for the logically challenged, pg.16-39
\sn
\hskip10pt (A) What are we after?, pg.16
\sn
\hskip10pt (B) The structure/non-structure dichotomy, pg.22
\sn
\hskip10pt (C) Abstract elementary classes, pg.30
\sn
\hskip10pt (D) Toward good $\lambda$-frames, pg.34
\bn
\S3 Good $\lambda$-frames, pg.40-49
\sn
\hskip10pt (A) getting a good $\lambda$-frame, pg.40
\sn
\hskip10pt (B) the successor of a good $\lambda$-frame, pg.42
\sn
\hskip10pt (C) the beauty of $\omega$ successive good $\lambda$-frames, pg.44
\bn
\S4 Appetite comes with eating, pg.50-57
\sn
\hskip10pt (A) The empty half of the glass, pg.50
\sn
\hskip10pt (B) The full half and half baked, pg.52
\sn
\hskip10pt (C) The white part of the map, pg.54
\bn
\S5 Basic knowledge, pg.58-60
\sn
\hskip10pt (A) knowledge needed and dependency of chapters, p.58
\sn
\hskip10pt (B) Some basic definitions and notation, pg.58
\bn
\S6 Symbols, pg.61-66
\newpage

\head {Abstract} \endhead
 \spuriousreset
\bigskip

Classification theory of elementary classes  deals with first
order (elementary) classes of structures (i.e. fixing a set $T$ of
first order sentences, we investigate the class of models of  
$T$ with the elementary submodel notion). It tries to find
dividing lines, prove their consequences, prove ``structure
theorems, positive theorems" on those in the ``low side" (in
particular stable and superstable theories), and prove
``non-structure, complexity theorems" on the ``high side".  It has
started with categoricity and number of non-isomorphic models.  It is
probably recognized as the central part of model theory, however
it will be even better to have such (non-trivial) theory for non-elementary
classes.  Note also that many classes of
structures considered in algebra are not first order; some
families of such classes are close to first order (say have kind
of compactness).  But here we shall deal with a classification theory for
the more general case without assuming knowledge of the first order
case (and in most parts not assuming knowledge of model theory at all).
\newpage

\head {\S0 Introduction and notation} \endhead  \resetall \sectno=0
 \spuriousreset
\bigskip

In \S2 we shall try to explain the purpose of the book to mathematicians with
little relevant background.  \S1
describes dividing lines and gives historical
background.  In \S5 we point out the (reasonably limited) background
needed for reading various parts and some basic definitions and in \S6
we list the use of symbols.
The content of the book is mostly described in \S2-\S3-\S4 but \S4 mainly deals
with further problems and \S6 with the symbols used. 

Is this a book?  I.e. is it a book or a collection of articles?  Well,
in content it is a book but the chapters have been written as
articles, (in particular has independent introductions and there are
some repetitions) and it was not clear
that they will appear together, see \S5(A) for more on how to read them.
\newpage

\head {\S1 Introduction for model theorists} \endhead  \resetall \sectno=1
 \spuriousreset
\bn
 
\bn
(A) \ub{Why to be interested in dividing lines}?

Classification theory for first order (= elementary) classes
is so established
now that up to the last few years most people tended to forget
that there are non-first order possibilities.  There are several good
reasons to consider these other possibilities; first, it is
better to understand a more general
context, we would like to prove stronger theorems by having wider context,
classify a larger family of classes.  
Second, understanding more general contexts may shed light
on the first order one.  In particular, larger families may have
stronger closure properties (see later). 
Third, many classes 
arising in "nature" are not first order (``in nature" here means other
parts of mathematics). 

Of course, we may suspect that applying to a wider context may
leave us with little content, i.e., the proofs may essentially be
just rewording of the old proofs (with cumbersome extra conditions);
 maybe there is no nice
theory, not enough interesting things to be discovered in this
context; it seems to me that  experience has already refuted the
first suspicion. Concerning the other suspicion, we shall try to
give a positive answer to it, i.e. develop a theory; on both see the
rest of the introduction.

In any case, ``not first order" does not define our family of
classes of models as discussed below. This is both witnessed from
the history (on which this section concentrates) and suggested by
reflection;  clearly we cannot prove much on arbitrary classes, so we need
some restriction to reasonable classes.  Now there may be
incomparable cases of reasonableness and a priori it is natural to 
expect to be able to say considerably 
more on the ``more reasonable" cases.  E.g. we
expect that much more can be said on first order classes than on
the class of models of a sentence from $\Bbb L_{\omega_1,\omega}$.

We are mainly interested here in generalizing the theorems on categoricity,
superstability and stability to such contexts, in particular we
consider the parallel of  {\L}o\'s Conjecture and the (very probably much
harder) main gap conjecture as test problems.

This choice of test problem is connected to the belief in (a),(b),(c)
discussed below (that motivates \cite{Sh:c}).
\mr
\item "{$(a)$}"  It is very interesting to find dividing lines and it
is a fruitful approach in investigating quite general classes of models.
\ermn
That is, we start with a large family of (in our case) classes (e.g., the
family of elementary (= first order) classes or the family of
universal classes or the family of locally finite algebras satisfying
some equations) and we would like to find
natural dividing lines.  A dividing line is not just a good property,
it is one for which we have some things to say on both sides: the
classes having the property and the ones failing it.  In our context
normally all the classes on one side, the ``high" one, will be provably 
``chaotic" by the non-structure side of our theory, and all the
classes on the other side, the ``low" one will have a positive 
theory.  The class of models of true arithmetic is a prototypical
example for a class in the ``high" side and the class of
algebraically closed field the prototypical non-trivial example in
the ``low" side.  

Of course, not all important and
interesting properties are like that.  If $F$ is a binary function on
a set $A$, not much is known to follow from $(A,F)$ not being a
group.  In model theory introducing o-minimal theories was 
motivated by looking for parallel
to minimal theories and attempts to investigate theories close to the real
field (e.g., adding the function $x \mapsto e^x$).  Their investigation
has been very important and successful, including parallels of stability
theory for strongly minimal sets, but it does not follow our 
paradigm.  A success of the guideline of 
looking for dividing lines had been the discovery of 
being stable (elementary classes, i.e. $(\text{Mod}_T,\prec)$, \cite{Sh:1}).  
\relax From this point of view to discover a dividing line means to prove
the existence of complementary properties from each side:
\mr
\item "{$(i)$}"   $T$ is unstable
iff it has the order property (recall that $T$ has the order property
means that: some first order formula $\varphi(\bar x,\bar y)$ 
linearly orders in $M$ some infinite $\bold I \subseteq {}^{\ell
g(\bar x)}M$ in a model $M$ of $T$)
\sn
\item "{$(ii)$}"  $T$ is stable iff $A \subseteq M 
\models T$ implies $|{\bold S}(A,M)|$, the set of $1$-types on $A$ for $M$ is
not too large $(\le |A|^{|T|})$.
\ermn
A case illustrating the point of
dividing line is a precursor of the order property,
property $E$ of Ehrenfeucht \cite{Eh57}, it says that some first order formula 
$\varphi(x_1,\dotsc,x_n)$ is asymmetric on some infinite 
$A \subseteq M,M$ a model of $T$; it is stronger than the order
property (= negation of stability).
A posteriori, order on the set of $n$-tuples is simpler; this 
is not a failure, what Ehrenfeucht did was fine for his aims,
but looking for dividing lines forces you to get the ``true" notion.

Even better than stable was superstable because it seems to me to 
maximize the ``area" which we view as being how
many elementary classes it covers times how much we can say about
them.  On the other hand, it has always seemed to me more interesting than
$\aleph_0$-stable as the failure of $\aleph_0$-stability is weak,
i.e. it has a few consequences.
There is a first order superstable not $\aleph_0$-stable 
class $K$ such that a model $M \in K$ is
determined up to isomorphism by a dimension (a cardinal) and a 
set of reals.  This exemplifies that an elementary class can fail to
be $\aleph_0$-stable but still is ``low": we largely can completely
list its models.  Such a class is the class of vector spaces over
$\Bbb Z/2 \Bbb Z$ expanded by predicates $P_n$ for independent
sub-spaces of co-dimension 2.  A model $M$ in this class is determined
up to isomorphism by one cardinal (the dimension of the sub-space $V_M
= \cap\{P^M_n:n \in \Bbb N\}$) and the quotient $M/V_M$ which has size
at most continuum (alternatively the set $\{\eta_a:a \in M\},\eta_a(n) \in
\{0,1\}$ and where $\eta_a = \langle
\eta_a(0),\eta_a(1),\ldots\rangle$ and 
$\eta_a(n)=0 \Leftrightarrow a \in P^M_n$).

Of course, the guidelines of looking for dividing lines if taken
religiously can lead you astray.  
It does not seem to recommend investigation of FMR (Finite Morley
Rank) elementary classes which has covered important ground 
(see e.g. Borovik-Nessin \cite{BoNe94}).
This guideline has helped, e.g. to discover dependent and
strongly dependent elementary classes, but so far our approach 
has seemingly not succeeded too much in advancing the investigation.

See more on this in end of \S2(B), in particular Question \scite{E53-nl.2.F}.
\mr
\item "{$(b)$}"  It is desirable to have an exterior a priori existing
goal as a test problem. 
\ermn
Such a problem in model theory was \L os conjecture which says: 
if a first order class of
countable vocabulary (= language) is categorical in one 
$\lambda > \aleph_0$ (= has one and only one model of cardinality 
$\lambda$ up to isomorphisms) then it is categorical in
every $\lambda > \aleph_0$.  At least for me so was 
Morley conjecture \cite{Mo65}
which says that for first order class with countable vocabulary, 
the number of its models of cardinality $\lambda > \aleph_0$ up to 
isomorphism is non-decreasing with $\lambda$.
This motivated my research in the early seventies which eventually
appeared as \cite{Sh:a} (with several late additions like local weight in
\cite[Ch.V,\S4]{Sh:a}).  
Now having introduced ``$\aleph_\varepsilon$-saturated models", it
seems unconvincing to understand $\dot I(\lambda,K)$, the
number of models in $K$ of cardinality $\lambda$ up to isomorphism,
for $K$ the class of $\aleph_\varepsilon$-saturated models of a first
order class, hence though essentially done then, was not written till
much later.
Eventually $``\dot I(\lambda,T)$ non-decreasing"
was done for the family of classes of models of a
countable first order theory (which was the original center of
interest; see \cite{Sh:c}).  

By this solution, there are 
very few ``reasons" for such $K = \text{ Mod}_T$ to 
have many models: being unstable,
unsuperstable, DOP (dimensional order property), OTOP (omitting type
order property) and deepness
(for fuller explanation see after \scite{E53-nl.1.3};  see more,
characterizing the family of functions $\dot I(\lambda,T)$ for countable $T$ 
in Hart-Hrushovski-Laskowski \cite{HHL00}).
So the direct aim was to solve the test question (e.g., the main gap
\footnote{which says that \ub{either} $\dot I(\lambda,T) = 2^\lambda$
for every ($> |T|$, or large enough) $\lambda$ \ub{or}
$\dot I(\aleph_\alpha,T) \le 
\beth_{\gamma(T)}(|\alpha|)$ for every $\alpha$ (for
some ordinal $\gamma(T)$); see more in \scite{E53-nl.1.2}.}),
but the motivation has always been the belief that
solving it will be rewarded with discovering worthwhile dividing lines
and developing a theory for both sides of each.

The point is that looking at the number of non-isomorphic models and
in particular the main gap we hope to develop a theory.  Other
exterior problems will hopefully give rise to other interesting
theories, which may be related to stability theory or may not; 
this was the point of \cite{Sh:10}, in particular the long list of
exterior results in the end of its introduction,
 and the words ``classification
theory" in the name of \cite{Sh:a}.  But, the above point 
seemingly was slow in being noticed.

Of course, if we consider the family of classes which are ``high" by
one criterion/dividing line, we expect that with respect to other
questions/dividing lines the ``previously high ones" 
will be divided and on a significant
portion of them we have another positive theory, quite reasonably
generalizing the older ones (but maybe we shall be led to very
different theories).  
E.g. for unstable first order classes \cite{Sh:93} succeeded 
in this respect: ``low ones" are the simple theories
and the ``high ones" are theories with the tree property (on
exciting later developments, see \cite{KiPi98} or \cite{GIL02}).
\mr
\item "{$(c)$}"  successful dividing lines will throw light on 
problems not considered when suggesting them.
\ermn
The point is that the theory should be worthwhile even if you discard
the original test problems.  Stability theory is just as interesting
for some other problems as for counting number of non-isomorphic
models.  E.g.
\mr
\item "{$(*)_1$}"  the maximal number of models no one embeddable into
another.
\ermn
This sounds very close to counting, so we expect this is to have a
closely related answer.

In fact for elementary classes (with countable vocabulary) which have a
structure theorem (see \scite{E53-nl.1.2} below), this number is $<
\beth_{\omega_1}$, for the others it is very much higher (see more on
the trichotomy after \scite{E53-nl.1.3}); so the answer to $(*)_1$ turns
out to be nicer than the one concerning the number, $\lambda \mapsto
\dot I(\lambda,T)$.
\mr
\item "{$(*)_2$}"  in ${\frak K}$ there are models very similar yet
non-isomorphic.
\ermn
This admits several interpretations which in general have complete and
partial solutions quite tied up with stability theory.  One is finding
$\Bbb L_{\infty,\lambda}$-equivalent not isomorphic models of cardinality
$\lambda$.  Stronger along this line are
EF$_\lambda$-equivalent not isomorphic.  Another is that there are
non-isomorphic models of $T$ such that a forcing neither collapsing
cardinals nor adding too short sequences makes them isomorphic.  For
non-logicians we should explain that this says in a very strong sense 
that there are no reasonable invariants, see \cite{Sh:225},
\cite{Sh:225a}, Baldwin-Shelah \cite{BLSh:464}, Laskowski-Shelah 
\cite{LwSh:489}, Hyttinen-Tuuri \cite{HyTu91}, Hyttinen-Shelah-Tuuri 
\cite{HShT:428}, Hyttinen-Shelah \cite{HySh:474}, 
\cite{HySh:529}, \cite{HySh:602}.
\mr
\item "{$(*)_3$}"  For which classes $K$ do we have: its models are no more
complicated than trees (in the graph theoretic sense say 
rooted graphs with no cycle)?
\ermn
This question was specified to having a tree of submodels which 
is ``free" (= ``non-forking") and it is a decomposition, i.e., the 
whole model is prime over the tree.  
This is answered by stability theory (for Mod$_T,T$ countable)
\mr
\item "{$(*)_4$}"  similarly replacing graphs with no cycles 
by another simple class, e.g., linear orders.
\ermn
This is very interesting, but too hard at present (see more in
Cohen-Shelah \cite{CoSh:919})
\mr
\item "{$(*)_5$}"   decidable theories, e.g. we 
may note that there was much done on decidability and
understanding of the monadic theory of some structures (in particular
Rabin's celebrated theorem).  Those works concentrated on linear orders
and on trees.  Was this because of our shortcoming or for inherent
reasons? 
\ermn
We may interpret this as a call to classify classes, in particular, first order
ones by their complexity as measured by monadic logic.  
This was carried to large
extent in Baldwin-Shelah \cite{BlSh:156} for first order classes.
Now this seems a priori orthogonal to classification taking number of
models as the test question; note that the class of linear orders is
unstable but 
reasonably low for \cite{BlSh:156}, whereas any class is maximally
complicated if it has a pairing function (e.g. a one-to-one function
$F^M$ from $P^M_1 \times P^M_2$ into $P^M_3$ while $P^M_1,P^M_2$ are
infinite) and there are such classes which are categorical in every
$\lambda \ge \aleph_0$.   In spite of all this
\cite{BlSh:156} relies heavily on stability theory; see \cite{Bl85},
\cite{Sh:197}, \cite{Sh:205}, \cite{Sh:284c}
\mr
\item "{$(*)_6$}"   the ordinal $\kappa$-depth of a model (Karp
complexity).
\ermn
For a model $M$ and a partial automorphism $f$ of $M$, Dom$(f)$ of
cardinality $< \kappa$, we can define its $\kappa$-depth in $M$, an ordinal (or
$\infty$) by Dp$_\kappa(f,M) \ge \alpha$ iff for every $\beta < \alpha$ and
subsets $A_1,A_2$ of cardinality $< \kappa$, there is a partial
automorphism $f'$ of $M$ extending $f$ of $\kappa$-depth $\ge \beta$ such that
$|\text{Dom}(f')| < \kappa,A_1 \subseteq \text{\rm Dom}(f'),A_2
\subseteq \text{\rm Rang}(f')$.

Let

$$
\align
\text{Dp}_\kappa(M) = \cup\{\text{Dp}_\kappa(f,M)+1:&f \text{ a partial
automorphism of } M \text{ of cardinality} \\
  & < \kappa \text{ and Dp}_M(f) < \infty\}.
\endalign
$$
\mn
This measures the complexity of the models and Dp$_\kappa(T) =
\cup\{\text{Dp}(M)+1:M$ a model of $T\}$ is a reasonable measure of
the complexity of $T$.  With considerable efforts, reasonable
knowledge concerning this measure was gained by Laskowski-Shelah
\cite{LwSh:560}, \cite{LwSh:687}, \cite{LwSh:871} confirming to some
extent the thesis above.
\mr
\item "{$(*)_7$}"  categoricity and number of models in
$\aleph_\alpha$, in ZF (i.e., with no choice).
\ermn
See \cite{Sh:840}.

You may view in this context the question of having non-forking (= abstracts 
dependence relations), orthogonality, regularity but for me this is 
part of the inside theory rather than an external problem
\smallskip

$(d) \quad$  non-structure is not so negative.
\mn
Now this book predominantly deals with the positive side, structure
theory, so defending the honour of non-structure is not really
necessary (it is the subject of \cite{Sh:e} though).  Still first we
may note that finding the maximal family of classes for which we know
something is considerably better than finding a sufficient condition.
In particular finding ``the maximal family ... 
such that ..." is finding dividing lines
and this is meaningless without non-structure results.

Second, this forces you to encounter real difficulties and develop
better tools; also using the complicated properties of a class which
already satisfies some ``low side properties" may require using and/or
developing a positive theory.  

Last but not least, non-structure from a different
perspective is positive.  Applying ``non-structure theory" 
to modules this gives
representation theorems of rings as endomorphism rings (see 
G\"obel-Trlifaj \cite{GbTl06}; note that the ``black boxes" 
used there started from
\cite[VIII]{Sh:c}).  In fact, generally for unstable elementary
class $K$, we can find models which in some respect represent a
pregiven ordered group (see \cite{Sh:800}).  This has been applied to
clarify in some cases to which
generalized quantifiers give a compact logic (see \cite{Sh:e} and more
in \cite{Sh:800}).

It may clarify to consider an alternative strategy: we have a reasonable
idea of what we look for and we have a specific class or structure
which should fit the theory.  This works when the analysis we have in
mind is reflected reasonably well in the specific case.  It may be
misleading when the examples we have, do not reflect the complexity of
the situation, and it seems to be the case in the problems we have at
hand.  More specifically, though the ``example" of the theory of
superstable first order classes stand before us, we do not try to take
the way of trying to assume enough of its properties so that it works;
rather we try look for dividing lines.
\nl
See more on ``why dividing lines" in the end of (B) of \S2.
\bn
(B) \ub{Historical comments on non-elementary classes}:

Let us return to non-elementary classes.  Generally, on model theory
for non-elementary classes see Keisler \cite{Ke71} and the handbook
\cite{BaFe85}: closer to our interest is the forthcoming book of
Baldwin \cite{Bal0x} and the older Makowsky \cite{Mw85}, mainly around
$\aleph_1$.

Below we present the results according to the kind of classes dealt
with (rather than chronologically). 
\nl
The oldest choice of families of classes 
(in this context) is the family of class of $\kappa$-sequence 
homogeneous models for a fixed $D$.

Morley and Keisler \cite{KM67} proved that there are at most $2^{2^{|T|}}$
such models of $T$ in any cardinality.  Keisler \cite{Ke71} proved that
if $\psi \in \Bbb L_{\omega_1,\omega}$ is categorical in $\aleph_1$
and its model in $\aleph_1$ is sequence homogeneous then it is
categorical in every $\lambda > \aleph_1$; generalizing (his version of)
the proof of Morley's theorem.  In \cite{Sh:3} instead of having a
monster ${\frak C}$, i.e., a 
$\bar \kappa$-saturated model of a first order $T$,
we have a $\bar \kappa$-sequence homogeneous model ${\frak C}$.  Let $D
= D({\frak C}) = \{\sftp(\bar a,\emptyset,{\frak C}):\bar a \in
{\frak C}$; i.e., $\bar a$ a finite sequence from ${\frak C}\}$;
note that $D,\bar \kappa$ determines ${\frak C}$ and we look at the
class of $M  \prec {\frak C}$ (or the class of
$(D,\lambda)$-homogeneous $M \prec {\frak C}$).
There the stability spectrum was reasonably characterized, 
splitting and strong splitting were introduced 
(for first order theory this was later
refined to forking).  See somewhat more in \cite{Sh:54}.
\sn
Lately, this (looking at the $\prec$-submodels of a
 $(D,\lambda)$-homogeneous monster ${\frak C}$)  has become very popular, see 
Hyttinen \cite{Hy98}, Hyttinen and Shelah \cite{HySh:629},
\cite{HySh:632}, \cite{HySh:629} (the main gap for
$(D,\aleph_\varepsilon)$-homogeneous models for a good diagram $D$),  
Grossberg-Lessman \cite{GrLe02}, \cite{GrLe0x} (the main gap for 
good $\aleph_0$-stable (= totally transcendental)), \cite{GrLe00a}, 
Lessman \cite{Le0x},
\cite{Le0y} (all on generalizing geometric stability).  

We may look at contexts which are closer to first order, i.e., having
some version of compactness.  Chang-Keisler \cite{ChKe62},
\cite{ChKe66} has looked at models with truth values in a topological
space such that ultraproducts can be naturally defined.
Robinson had looked at model theory of the
classes of existentially closed models of first order universal or
just inductive theories.  Henson \cite{He74} and Stern \cite{Str76} 
have looked at
Banach spaces (we can take an ultraproduct of the spaces, throw away
the elements with infinite norm and divide by those with infinitesimal
norm).  Basically the logic is ``negation deficient", see Henson-Iovino
\cite{HeIo02}.

The aim of \cite{Sh:54} was to show that the most basic stability theory
was doable for Robinson style model theory.  In
particular it deals with case II (the models of a universal first
order theory which has the amalgamation property) and
case III (the existentially closed models of
a first order inductive $(= \Pi^1_2)$ theory); those are particular cases
of $(D,\lambda)$-homogeneous models.  Case II is a special case of
III where $T$ has amalgamation.   Lately, 
Hrushovski dealt with Robinson classes (= case II above).  A  
Ph.D. student of mine in the seventies was supposed to deal
with Banach spaces but this has not materialized.  Henson and Iovino
continued to develop model theory of Banach spaces.  Lately, interest
in the classification theory in such contexts has awakened and dealing
with cases II and III and complete metric spaces and
Banach spaces and relatives, now called continuous model theory, 
see Ben-Yaacov \cite{BY0y}, Ben-Yaacov
Usvyatsov \cite{BeUs0x}, Pillay \cite{Pi0x}, Shelah-Usvyatson \cite{ShUs:837}.

The most natural stronger (than first order)
logic to try to look at, in this context, has been 
$\Bbb L_{\omega_1,\omega}$ and even $\Bbb L_{\lambda^+,\omega}$.  By 1970
much was known on $\Bbb L_{\omega_1,\omega}$ (see Keisler's book
\cite{Ke71}); however, 
if you do not like non-first order logics, look at the class of atomic
models of a countable first order $T$.  The general question looks
hard.  At the early seventies I have clarified some things on $\psi
\in \Bbb L_{\omega_1,\omega}$ categorical in $\aleph_1$, but it was
not clear whether this leads to anything interesting.  Then the
following question of Baldwin catches my eye (question 21 of the
Friedman list \cite{Fr75})
\mr
\item "{$(*)_1$}"   can $\psi \in \Bbb L(\bold Q)$ have exactly
one uncountable model up to isomorphism? \nl
$\bold Q$ stands for the quantifier ``there are uncountably many"
\ermn
This is an excellent question, a partial answer was (\cite{Sh:48})
\mr
\item "{$(*)_2$}"  if $\diamondsuit_{\aleph_1}$ and 
$\psi \in \Bbb L_{\omega_1,\omega}(\bold Q)$ has at least one but
$< 2^{\aleph_1}$ models in $\aleph_1$ up to isomorphism \ub{then} it
has a model in $\aleph_2$ (hence has at least 2 non-isomorphic models)
\ermn
Only later the original problem (even for $\psi \in \Bbb
L_{\omega_1,\omega}(\bold Q))$ was solved in ZFC, see below.
It seems natural to ask in this case how many models $\psi$ has in
$\aleph_2$, and then successively in $\aleph_n$ (raised in
\cite{Sh:48}), but as it was hard enough, 
the work concentrates on the case of $\psi \in \Bbb L_{\omega_1,\omega}$, 
so (\cite{Sh:87a}, \cite{Sh:87b} and generalizing it to cardinals
$\lambda,\lambda^+,...$ is a major aim of this book):
\mr
\item "{$(*)_3$}"  $(a) \quad$ if $n < \omega,2^{\aleph_0} <
2^{\aleph_1} < \ldots < 2^{\aleph_n},\psi \in \Bbb L_{\omega_1,\omega},
\dot I(\aleph_\ell,\psi) < \mu_{\text{wd}}(\aleph_\ell)$, for
\footnote{$\mu_{\text{wd}}(\aleph_\ell)$ is ``almost" 
equal to $2^{\aleph_\ell}$}
\nl

\hskip25pt $\ell \le n$ and $\dot I(\aleph_1,\psi) \ge 1$ \ub{then} 
$\psi$ has a model in $\aleph_{n+1}$ and 
\nl

\hskip25pt without loss of generality $\psi$ is categorical in $\aleph_0$
\sn
\item "{${{}}$}"  $(b) \quad$ if the assumption of (a) holds for every $n
< \omega$ and $\psi$ is for simplicity 
\nl

\hskip25pt categorical in $\aleph_0$
\ub{then} the class Mod$_\psi$ is so-called excellent (see (c))
\sn
\item "{${{}}$}"  $(c) \quad$ if $\psi \in \Bbb L_{\omega_1,\omega}$ 
is excellent and is categorical in one $\lambda > \aleph_0$ then it is
\nl

\hskip25pt categorical in every $\lambda > \aleph_0$.
\ermn
Essentially, it was proved that excellent $\psi \in \Bbb
L_{\omega_1,\omega}$ are very similar to $\aleph_0$-stable (= totally
transcendental) first order countable theories (after some ``doctoring").  The
set of types over a model $M,{\Cal S}(M)$ is restricted (to not violate the
omission of the 
types which every model of $\psi$ omit).  The types themselves are as
in the first order case, set of formulas
but we should not look at complete types over any
$A \subseteq M \models \psi$, only at the cases $A = N \prec M$ or $A
= M_1 \cup M_2$ where $M_1,M_2$ are stably amalgamated over $M_0$ and more
generally at $\cup\{M_u:u \in {\Cal P}^-(n)\}$, where $\langle M_u:u
\in {\Cal P}^-(n)\rangle$ is a ``stable system".

This work was continued in Grossberg and Hart \cite{GrHa89}, (main
gap), Mekler and Shelah \cite{MkSh:366} (dealing with free algebras),
Hart and Shelah \cite{HaSh:323} (categoricity may hold for
$\aleph_0,\aleph_1,\aleph_2,\dotsc,\aleph_n$ 
but fail for large enough $\lambda$) and lately Zilber \cite{Zi0xa},
\cite{Zi0xb} (connected to his programs).
Further works on more general but not fully general are \cite{Sh:300},
Chapter II (universal classes), Shelah and Villaveces
\cite{ShVi:635}, van Dieren \cite{Va02} (abstract elementary class 
with no maximal models).  See also the closely
related Grossberg and Shelah \cite{GrSh:222}, \cite{GrSh:238},
\cite{GrSh:259}, \cite{Sh:394}, (abstract elementary class with amalgamation),
Grossberg \cite{Gr91} and Baldwin and Shelah \cite{BlSh:330},
\cite{BlSh:360}, \cite{BlSh:393}.
Lately, Grossberg and VanDieren \cite{GrVa0xa}, \cite{GrVa0xb}
Baldwin-Kueker-VanDieren \cite{BKV0x} investigate the
related tame abstract elementary class 
including upward categoricity.  They prove independently of
\marginbf{!!}{\cprefix{734}.\scite{734-am3.6}} that tame a.e.c. with amalgamation has nice
categoricity spectrum; i.e. prove categoricity in cardinals $> \mu$ in
the relevant cases; in the notation here ``tame''
means locality of orbital types over saturated model; on
\marginbf{!!}{\cprefix{734}.\scite{734-am3.6}}, see \S4(B) after $(**)_\lambda$.
Concerning $\Bbb L_{\kappa,\omega}$, see
Makkai-Shelah \cite{MaSh:285} (on cateogoricity of $T \subseteq
\Bbb L_{\kappa,\omega},\kappa$ compact starting with $\lambda$
successor), Kolman-Shelah \cite{KlSh:362} ($T \subseteq \Bbb
L_{\kappa,\omega},\kappa$ measurable, amalgamation derived from
categoricity), \cite{Sh:472} ($T \subseteq \Bbb
L_{\kappa,\omega},\kappa$ measurable, only down from successor).
See more in the book \cite{Bal0x} of Baldwin on the subject.

Going back, $(*)_3$ deals with $\psi \in \Bbb L_{\omega_1,\omega}$, it
generalizes the case $n=1$ which, however, deals with $\psi \in \Bbb
L_{\omega_1,\omega}(\bold Q)$.  On the other hand, $\psi \in \Bbb
L_{\omega_1,\omega}(\bold Q)$ is not a persuasive end of the story as
there are similar stronger logics.  Also the proof deals with $\Bbb
L_{\omega_1,\omega}(\bold Q)$ in an indirect way, we look at a related
class $K$ which has also countable models but some first order
definable set should not change when extending.  So it seems that the
basic notion is the right version of elementary extensions.
This leads to analysis which suggests the notion of
abstract elementary class, 
${\frak K}$ with LST$({\frak K}) \le \aleph_0$ which, moreover,
is PC$_{\aleph_0}$ (in \cite{Sh:88}, represented here in \chaptercite{88r}).

Now much earlier Jonsson \cite{Jn56}, \cite{Jn60} had considered
axiomatizing classes of models.  Compared with the abstract elementary
classes used (much later) in \cite{Sh:88}=\chaptercite{88r}, the main
\footnote{Jonsson axioms were, in our notations, (for a fix vocabulary $\tau$,
finite in \cite{Jn56}, countable in \cite{Jn60}), $K$ is a class of
$\tau$-models satisfying
\mr
\widestnumber\item{$(IV)'$}
\item "{$(I)$}"  there are non-isomorphic $M,N \in K$ in \cite{Jn56}
\sn
\item "{$(I)'$}"  $K$ has members of arbitrarily large cardinality in
\cite{Jn60}
\sn
\item "{$(II)$}"  $K$ is closed under isomorphisms
\sn
\item "{$(III)$}"  the joint embedding property
\sn
\item "{$(IV)$}"  disjoint amalgamation in \cite{Jn56}
\sn
\item "{$(IV)'$}"  amalgamation in \cite{Jn60}
\sn
\item "{$(V)$}"  $\cup\{M_\alpha:\alpha < \delta\} \in K$ if $M_\alpha
\in K$ is $\subseteq$-increasing
\sn
\item "{$(VI)$}"  if $N \in K$ and $M \subseteq N$ (so $|M| \ne
\emptyset$ but not necessarily $M \in K$) and $\alpha > 0,\|M\| <
\aleph_\alpha$ then there is $M' \in K$ such that $M \subseteq M'
\subseteq N$ and $\|M'\| < \aleph_\alpha$ (this is a strong form of
the LST property).
\ermn
Note that for an abstract elementary class $(K,\le_{\frak K})$, if
$\le_{\frak K} = \subseteq \restriction K$, then AxIV (smoothness)
and AxV (if $M_1 \subseteq M_2$ are $\le_{\frak K}$-submodels of $N$
then $M_1 \le_{\frak K} M_2$) of \marginbf{!!}{\cprefix{88r}.\scite{88r-1.2}} or \marginbf{!!}{\cprefix{600}.\scite{600-0.2}}
and part of AxI become trivial (hence are missing from Jonsson
axioms), the others give II,
and a weaker form of VI (specifically, for one $\aleph_\alpha$,
i.e. $\aleph_\alpha = \text{ LST}({\frak K})^+$, the other cases are proved).} 
differences are that he uses the order $\subseteq$ (being a submodel)
on $K$ (rather than an abstract order $\le_{\frak K}$) and assume the
amalgamation (and JEP joint embedding property).
His aim was to construct and axiomatize the construction of
universal and then universal homogeneous models so including
amalgamation was natural; Morley-Vaught \cite{MoVa62}
use this for elementary class.  In fact if we add amalgamation (and JEP) to
abstract elementary classes we get such theorems (see
\sectioncite[\S2]{88r}, in fact we also get uniqueness in a case of somewhat
different character, \marginbf{!!}{\cprefix{88r}.\scite{88r-2.11}}).  From our perspective
amalgamation (also $\le_{\frak K} = \subseteq$) is a heavy
assumption (but an important property, see later).  Now, model
theorists have preferred saturated to universal homogeneous and prefer
first order classes (Morley-Vaught \cite{MoVa62}, Keisler replete) 
with very good reasons, as it is better (more
transparent and give more) to deal with one element than a model.
That is, assume our aim is to show that $N$ from our class $K$ is
universal, i.e., we are given $M \in K$ of cardinality not larger than
that of $N$ and we have to construct an (appropriate) embedding of $M$
into $N$.  Naturally, we do it by approximations of cardinality
smaller than $\|M\|$, the number of elements of $M$.  Jonsson uses as
approximations isomorphisms $f$ from a submodel $M'$ of $M$ of
cardinality $<\|M\|$.  Morley and Vaught use functions from a subset
$A$ of $M$ into $N$ such that: if $n < \omega,a_0,\dotsc,a_{n-1} \in A$
satisfy a first order formula in $M$ then their image satisfies it
in $N$.  So they have to add one element at each step which is better
than dealing with a structure. 
In fact, also in this book, for a different notion of type, the types of
elements continue to play a major role (\ub{but} we use types which
are not sets of formulas over models).  So we try to have ``the best
of both approaches" - all is done over models from $K$, but we ask
existence, etc., only of singletons, for this reason in the proof of
the uniqueness of ``saturated" models we have to go ``outside" the two
models, build a third (see \marginbf{!!}{\cprefix{300b}.\scite{300b-3.10}} or \marginbf{!!}{\cprefix{600}.\scite{600-0.19}}).

Here we have chosen abstract elementary class 
as the main direction.  This includes classes
defined by $\psi \in \Bbb L_{\omega_1,\omega}$ and we can analyze
models of $\psi \in \Bbb L_{\omega_1,\omega}(\bold Q)$ in such context
by a reduction.  In \cite{Sh:88} = \chaptercite{88r} 
Baldwin's question was solved in ZFC.
Also superlimit models were introduced and amalgamation in $\lambda$
was proved assuming categoricity in $\lambda$ and $1 \le \dot I
(\lambda^+,{\frak K}) < 2^{\lambda^+}$ when $2^\lambda <
2^{\lambda^+}$.   The intention of the work was to prepare
the ground for generalizing \cite{Sh:87b}.  Note that sections \S4,\S5
from \chaptercite{88r} are harder than the parallel in \cite{Sh:87a} because we
deal with abstract elementary class 
(not just $\psi \in \Bbb L_{\omega_1,\omega}(\bold Q)$).

Now \cite{Sh:300} deals with universal classes.  This family is
incomparable with first order and \cite{Sh:155} gives hope it will be
easier.  Note that in excellent classes the types are set of formulas
and this is true even for \chaptercite{88r} though the so-called 
materializing replaces
realizing a type.  In \cite{Sh:300} (orbital)-type 
is defined by $\le_{\frak K}$-mapping.
Surprisingly we can still show ``$\lambda$-universal homogeneous" is
equivalent to $\lambda$-saturated under the reasonable interpretations
(so have to find an element rather than
a copy of a model) what was a strong argument for sequence homogenous
models (rather than model homogeneous).

In \cite{Sh:576}, which is a prequel of the work here, (redone in
\cite{Sh:E46}) we generalize
\cite{Sh:88} to any abstract elementary class ${\frak K}$ having no remnant of
compactness, see on it below.  On \chaptercite{600}, \chaptercite{705} see later.

I thank the institutions in which various parts of this book 
were presented and the
student and non-students who heard and commented.
Earlier versions of
\cite{Sh:300a}, \cite{Sh:300b}, \cite[III]{Sh:e}, \cite{Sh:300c},
\cite{Sh:300d}, \cite{Sh:300e} were presented in Rutgers in 1986; some
other parts were represented some other time.  In
Helsinki 1990 a lecture was on the indiscernibility from \cite{Sh:300f},
\cite{Sh:300g}.  First version of 
\cite{Sh:576} was presented in seminars in the
Hebrew University, Fall '94.
The G\"odel lecture in Madison Spring 1996 was on \cite{Sh:576} and 
\chaptercite{600}.  The author's lecture in the logic
methodology and history of science, Kracow '99, was on \chaptercite{600} and
\chaptercite{705}.
In seminars at the Hebrew University, \chaptercite{88r} was presented in
Spring 2002, \cite{Sh:576} was presented in 98/99, \chaptercite{600} +
\chaptercite{734} were presented in 99/00, \chaptercite{600} + \chaptercite{705} were
presented in 01/02 and my lecture in the 
Helsinki 2003 ASL meeting was on good $\lambda$-frames and \chaptercite{734}.  

I thank John Baldwin, Emanuel Dror-Farajun, Wilfred Hodges, Gil Kalai,
Adi Jarden, Alon Siton, Alex Usvyatsov, Andres Villaveces for 
many helpful comments and error detecting in the introduction
(i.e. Chapter N).

Last, but not least, I thank Alice Leonhardt for beautifully typesetting
the contents of this book.
\newpage

\head {\S2 Introduction for the logically challenged} \endhead  \resetall \sectno=2
 \spuriousreset
\bn
(This is recommended reading for logicians too, but there are some repetitions
of part (A) of \S1).

This is mainly an introduction to \chaptercite{600}, \chaptercite{705}. 
\sn
\beginaside
We assume the reader knows the notion of an infinite cardinal 
\ub{but not} that he knows about first
order logic (and first order theories); for reading (most of) the
book, not much more is needed, see \S5.
\mn
Paragraphs assuming more knowledge or are not so essential will be in
indented, e.g. when a result is explained ignoring some qualifications
and we comment on them in indented text.
\endaside
\bn
(A) \ub{What are we after}?

This introduction is intended for a general mathematical audience.  We may
view our aim in this book as developing a theory dealing with abstract
classes of mathematical structures that will also be referred to as
models.  Examples of structures are the field $\Bbb R$, any group and any
ring.  The classes of models we consider are called 
``abstract elementary classes" or briefly a.e.c.  An abstract
elementary class ${\frak K}$ is a class of structures denoted by
$K$ together with an order relation denoted by $\le_{\frak K}$ which
distinguishes for each structure $N$ a certain family $\{M \in K:M
\le_{\frak K} N\}$ of substructures (= submodels).

First, rather than giving a formal definition, we will give several examples:
\nl
\margintag{E53-nl.0.6X}\ub{\stag{E53-nl.0.6X} Examples}:
\mr
\widestnumber\item{$(iii)$}
\item "{$(i)$}"  the class of groups where the order relation is
``being a subgroup".
\ermn
In this example $\le_{\frak K}$ is simply being substructures.
(In the sequel when we do not specify the order relation is means
simply to take all substructures). 
\mr
\item "{$(ii)$}"  The class of algebraically closed fields with
characteristic zero
\sn
\item "{$(iii)$}"  the class of rings
\sn
\item "{$(iv)$}"  the class of nil rings, i.e. ring $R$ such that
for every $x \in R$, $x^n=0$ for some $n \ge 1$
\sn
\item "{$(v)$}"  the class of torsion $R$-modules for a fixed ring $R$
\sn
\item "{$(vi)$}"  the class of $R$-modules for a fixed ring $R$ but
unlike the previous cases the relation of $\le_{\frak K}$ is not just
being a submodule, but
being a ``pure submodule" \footnote{A left $R$-module $M$ is a pure
submodule of a left $R$-module $N$ when if $rx=y,x \in N$ and $y \in M$ then
$rx'=y$ for some $x' \in M$}
\sn
\item "{$(vii)$}"  the class of rings but $R_1 \le_{\frak K} R_2$
means here: $R_1$ is a subring of $R_2$ and if $R'_2$ is a finitely
generated subring of $R_2$ then $R_1 \cap R'_2$ is a finitely
generated subring of $R_1$
\sn
\item "{$(viii)$}"  the class of partial orders
\sn
\item "{$(ix)$}"  concerning Hill Lemma, Baldwin, Eklof and Trlifaj
\cite{BETp06} show it fit in a.e.c. context.
\ermn
Abstract elementary class form an extension of the notion of
elementary class which mean a class of structures which are models of
a so-called first order theory.  
The notion of abstract elementary classes, while
more general, does not rely on elementary classes and indeed, for
reading this introduction we do not assume knowledge of first order
logic.

We will be mainly interested in this book in finding parallel to the
``superstability theory" which is part of the ``classification theory"
(this is explained below; on the first order case 
see, e.g. \cite{Sh:c}, \cite{Sh:200} or other books on the subject,
e.g. Baldwin \cite{Bal88}).

Superstability theory can be described as dealing with elementary
classes of structures for which there is a good dimension theory; 
see on our broader aim below.

A structure $M$ will have a so-called vocabulary $\tau_M$ (this is its
``kind", e.g. is it a ring or a group).  
Note that for each class ${\frak K} = (K,\le_{\frak K})$ we shall
consider, all $M \in K$ has the same vocabulary 
(sometimes called language), which we denote by $\tau =
\tau_{\frak K}$, e.g., for a class of fields it is $\{+,\times,0,1\}$
where $+,\times$ are binary functions symbols interpreted in each field as
two-place functions and similarly $0,1$ are individual 
constant symbols.  We may
have also relations, (in example $(viii)$ the partial order is a
relation), note that relation symbols are usually called predicates.
The reader may restrict himself to the case of countable or even
finite vocabulary with function symbols only.  
We certainly demand each function symbol
to have finitely many places (and similarly for relation symbols).

We try now, probably prematurely, to give exact definitions of some
basic notions toward what long term goal we would like to advance,
probably it will make more sense after/if the reader continues to read
the introduction.  (But most of this will be repeated and expanded).

We think that the family of abstract elementary classes ${\frak K}$
(defined in \scite{E53-nl.0.7} below) can be divided, in some ways, so that we
can say significant things both on the ``low", simple side and on the
``high, complicated" side.  This sounds vague, can we already state a
conjecture?  It seems reasonable that a class $K$ with a unique member
(up to isomorphism, of course)
in a cardinality $\lambda$ is simple; but what can be the class of
cardinals for which this holds?  This class is called
the ``categoricity spectrum of the abstract elementary class ${\frak K}$"
(see Definitions \scite{E53-nl.0.7}, \scite{E53-nl.0.9} below), we conjecture
that is a simple set, e.g. contains every
large enough cardinal \ub{or} does not contain every large enough
cardinal.  Moreover, this also applies to the so-called superlimit spectrum of
${\frak K}$ (see Definition \scite{E53-nl.0.10}).  In the ``low, simple"
case we have, e.g. a dimension theory for ${\frak K}$, and in the 
``high case" we can prove the class is complicated and so cannot have
such a nice theory (this paragraph will be explained/expanded later).

Here we make some advances in this direction.
\nl
First, what exactly is an abstract elementary class?  It is much easier to
explain than the so-called ``elementary classes" which is defined
using (first order) logic.  A major feature are closure under
isormorphism and unions.
\definition{\stag{E53-nl.0.7} Definition}  ${\frak K} = (K,\le_{\frak K})$
is an abstract elementary class \ub{when}
\mr
\item "{$(A)$}"  $(a) \quad K$ is a 
class of structures all of the same ``kind", i.e. vocabulary;
\nl

\hskip25pt  e.g. they can be all rings or all graphs, $\tau$ denote a
vocabulary 
\sn
\item "{${{}}$}"  $(b) \quad K$ is closed under isomorphisms
\sn
\item "{${{}}$}"  $(c) \quad \le_{\frak K}$ is a partial order of $K$, also
closed under isomorphisms
\nl

\hskip25pt  and $M \le_{\frak K} N$ implies $M \subseteq N,M$ a
substructure of $N$ and, of 
\nl

\hskip25pt course, $M \in K \Rightarrow M \le_{\frak K} M$ 
\sn
\item "{${{}}$}"  $(d) \quad K$ (and $\le_{\frak K}$) are closed under direct
limits, or, what is equivalent,
\nl

\hskip25pt  by unions of $\le_{\frak K}$-increasing chains, 
i.e. if $I$ is a linear order and
\nl

\hskip25pt $M_t(t\in I)$ is
$\le_{\frak K}$-increasing with $t$ then $M = \cup\{M_t:t \in I\}$ belongs 
\nl

\hskip25pt to $K$ and; morever, $t \in I \Rightarrow M_t \le_{\frak K} M$
\sn
\item "{${{}}$}"  $(e) \quad$ similarly to clause (d)
inside $N \in K$, i.e., if $t \in I \Rightarrow M_t \le_{\frak K} N$ 
\nl

\hskip25pt then $M \le_{\frak K} N$.
\ermn
Two further demands are only slightly heavier
\mr
\item "{$(B)$}"  $(f) \quad$ if 
\footnote{this certainly holds if $\le_{\frak K}$
is defined as $\prec_{{\Cal L}(\tau({\frak K}))}$ for some logic
${\Cal L}$} $M_\ell \le_{\frak K} N$ for $\ell=1,2$ and $M_1
\subseteq M_2$ then $M_1 \le_{\frak K} M_2$
\sn
\item "{${{}}$}"  $(g) \quad (K,\le_{\frak K})$ has countable 
character, which means that every structure 
\nl

\hskip25pt can be approximated by countable ones; i.e., if
$N \in K$ then 
\nl

\hskip25pt every countable set of elements of $N$ is included in some
\nl

\hskip25pt  countable $M \le_{\frak K} N$
(in the book but not in the 
\nl

\hskip25pt  introduction we allow replacing
``countable" by ``of cardinality 
\nl

\hskip25pt $\le \text{ LST}({\frak K})"$ for some fixed cardinality
LST$({\frak K}))$.
\ermn
Not all natural classes are included, e.g. the class of Banach spaces
is not, as completeness is not preserved by unions of increasing
chains.  Still it seems very broad and the question is can we
prove something in such a general setting.
\enddefinition
\bigskip

\definition{\stag{E53-nl.0.9} Definition}  1) ${\frak K}$ (or $K$) is
categorical in $\lambda$ \ub{when} it has one and only one model of
cardinality $\lambda$ up to isomorphism.
\nl
2) The categoricity spectrum of ${\frak K}$, cat$({\frak K})$, is the
 class of cardinals $\lambda$ in which ${\frak K}$ is categorical.
\sn

A central notion in model theory is elementary classes or first order
classes which are defined using so called first order logic (which the
general reader is not required here to know, it is explained in the 
indented text below).

Each such class is the class of models of a first order theory with
the partial order $\prec$.  

Among elementary classes, a major division is between the so-called superstable
ones and the non-superstable ones, and for each superstable one there
is a dimension theory (in the sense of the dimension of a vector
space).  Our long term aim in restricted terms is to find such
good divisions for abstract elementary classes, though we do not like
to dwell on this further now, it seems user-unfriendly not to define them at
all, so for the time being noting 
that for elementary classes being superstable is equivalent to
having a superlimit model in every large enough cardinality; also
noting that superstability for abstract elementary classes suffer
from schizophrenia, i.e. there are several different definitions which
are equivalent for elementary classes, the one below is one of them.
\enddefinition
\bigskip

\definition{\stag{E53-nl.0.10} Definition}  Let ${\frak K}$ be an abstract
elementary class.
\nl
1) We say $f$ is a $\le_{\frak K}$-embedding of $M$ into $N$ when $f$
   is an isomorphism of $M$ onto some $M' \le_{\frak K} N$.
\nl
2) ${\frak K}_\lambda = (K_\lambda,\le_{{\frak K}_\lambda})$ where
$K_\lambda = \{M \in {\frak K}:\|M\| = \lambda\}$ and
$\le_{{\frak K}_\lambda} = \le_{\frak K} \restriction K_\lambda$.
\nl
3) An abstract elementary class ${\frak K}$ is superstable \ub{iff} for every
large enough $\lambda$, there is a superlimit structure $M$ for
${\frak K}$ of cardinality $\lambda$; where 
\nl
4) We say that $M$ is a superlimit (for ${\frak K}$) 
when for some (unique) $\lambda$
\mr
\item "{$(a)$}"  $M \in {\frak K}$ has cardinality $\lambda$
\sn
\item "{$(b)$}"  $M$ is $\le_{\frak K}$-universal, i.e., if $M' \in
K_\lambda$ then there is a $\le_{\frak K}$-embedding of $M'$ into $M$,
in fact with range $\ne M$
\sn
\item "{$(c)$}"  for any $\le_{\frak K}$-increasing chain of models
isomorphic to $M$ with union of cardinality $\lambda$, the union is isomorphic
to $M$.
\ermn
5) The superlimit spectrum of ${\frak K}$ is the class of $\lambda$
such that there is a superlimit model for ${\frak K}$ of
cardinality $\lambda$.  
\sn
We shall return to those notions later.

What about the examples listed above?  Concerning 
the strict definition of elementary classes as classes of the
form $(\text{Mod}_T,\prec)$ defined below, 
among the examples in \scite{E53-nl.0.6X} 
the class of algebraically closed fields 
(example $(ii)$) is an elementary class since it can be proved
that being a sub-field is equivalent to being an elementary
substructure for such fields.

In the example $(i)$, the class of models is elementary, i.e., equal to
Mod$_T$: the class of groups, but the order is not $\prec$ but
$\subseteq$.  This is true also in the examples $(iii)$, rings and
$(viii)$, partial orders.
\nl
In the example $(vi)$, the class of torsion $R$-modules is not a 
first order class as we have to say $(\forall x) \dsize \bigvee_{r
\in R \backslash \{0\}} rx = 0$ and we really need to use an infinite
disjunction.  
The situation is similar for the class of nil rings (example (iv)).  In
example $(vii)$, the class of rings with $\le_{\frak K}$ defined using
finitely generated subrings not only is the class of structures not
elementary but $\le_{\frak K}$ is neither $\prec$ nor $\subseteq$.
In the example $(vii)$, $R$-modules, $K$ is elementary but
$\le_{\frak K}$ is different.
\enddefinition
\bn
\beginaside
Recall \footnote{we urge the logically challenged: when lost, jump ahead}
the traditional frame of model theory are the so-called elementary (or
first order) classes.
That is, for some vocabulary $\tau$, and set $T$ of so-called
sentences in first order logic in this vocabulary, $K = \text{ Mod}_T
= \{M:M$ a $\tau$-structure satisfying every sentence of $T\}$ and
$\le_{\frak K}$ being $\prec$, ``elementary submodel".
Recall that $M \prec N$ if $M \subseteq N$ and for every first order
formula $\varphi(x_0,\dotsc,x_{n-1})$ in the (common) vocabulary,
i.e., from the language $\Bbb L(\tau)$ and
$a_0,\dotsc,a_{n-1} \in M,\varphi(a_0,\dotsc,a_{n-1})$ is satisfied by $M$,
(symbolically $M \models \varphi[a_0,\dotsc,a_{n-1}]$) \ub{iff} $N$
satisfies this.  

Now here an elementary class is one of the form (Mod$_T,\prec$), any
such class is an abstract elementary class (see below).  A different
abstract elementary class 
derived from $T$ is (Mod$_T,\subseteq$) but then we should restrict
ourselves to $T$ being a set of universal sentences or 
just $\Pi_2$-sentences as we like to have closure under direct
limits.  For each such $T$ another abstract elementary class which can be
derived from it is $(\{M \in \text{ Mod}_T:M$ is existentially
closed$\},\subseteq$).

We are \ub{not disputing} the choice of first order classes as \ub{central} in
model theory but there are many interesting other classes.
Most notably for algebraists are classes of locally finite structures
and for model theorists are (Mod$_\psi,\prec_{\Cal L})$ where 
$\psi$ belongs to the logic denoted by 
$\Bbb L_{\omega_1,\omega}(\tau)$ or just $\psi \in \Bbb
L_{\lambda^+,\omega}(\tau)$ for some $\lambda$ where ${\Cal L}$ is a 
fragment of this logic to which the sentence $\psi$
belongs; if $\psi \in \Bbb L_{\omega_1,\omega}(\tau)$ we may choose a
countable such ${\Cal L}$.

(This logic may seem obscure to non-logicians but it just means that we
allow to say $\dsize \bigwedge_{i \in I}
\varphi_i(x_0,\dotsc,x_{n-1})$ where $I$
has at most $\lambda$ members so enable us to say ``a ring is nill,
locally finite, etc.", \ub{but} not ``$<$ is a well ordering").

In some sense if we look at classification theory of elementary
classes as a building, we note that several ``first floors"
disappear (in the context of abstract elementary class) 
but we aim at saving considerable part of
the rest (of course not all) by developing a replacement for those
lower floors.  

We may put in the basement the
downward LS theorem (there are small $N \prec M$); it survives.  
But not so the compactness theorem even very weak forms like 
``if $\bar a = \langle a_n:n
\in \Bbb N\rangle,\bar b = \langle b_n:n \in \Bbb N\rangle$ are
sequences of members of $M$ and $f_n$ is an automorphism of $M$ mapping
$\bar a \restriction n$ to $\bar b \restriction n$ then some
$\le_{\frak K}$-extension of $M$ has an automorphism mapping $\bar a$
to $\bar b$" do not hold in arbitrary a.e.c.  
(Note that for ``$(D,\lambda)$-homogeneous models" (e.g.
\cite{Sh:3}) such forms of compactness hold.  The
point of \cite{Sh:394} is to start investigating classes for which
all is nice except that types are not determined by their small
restrictions, that is, defining $\Bbb E^\kappa_N = \{(p,q):p,q \in
{\Cal S}(N)$ and $M \in K_\kappa \Rightarrow p \restriction M = q
\restriction M\}$, this is, a priori, not the equality
(\cite[1.8,1.9,pg.4]{Sh:394}).  
We lose as well the upward LST theorem (every model has a
proper  $<_{\frak K}$-extension); (those fit the first floor).

Also in abstract elementary classes the roles of 
formulas disappear.  Hence we lose the
notion of the type of an element
$a$ over a set $A$ inside a model $M$; so the second floor
including the ``$\kappa$-saturated model" (in the traditional sense) goes
 down the drain as the types disappear.

What is saved?  (I.e. not by definitions but in the positive case of a
dividing line which has a non-structure result.)
In a suitable sense, we save: non-forking amalgamation of 
models, prime models, a decomposition of a model
over a non-forking tree of models (a relative of free amalgamation),
and for a different notion of type, being (saturated and) 
orthogonal, regular and eventually
the main gap for the parallel of $\aleph_\varepsilon$-saturated model
of a superstable $T$.
\endaside
\bn
We now try to describe our aim in broad terms; if this seems vague, in
(B) below we describe it in a restricted case more concretely.
Our aim is to consider a family of classes ${\frak K}$ (all the
``reasonable" classes) and try to \ub{classify} them in the sense of
taxonomy, we look for \ub{dividing lines} among them.  This means
dividing the family to two, one part are those which are ``high",
``complicated".  Typically we have for each ${\frak K}$ in the ``high side" a
\ub{non-structure} result, saying there are many complicated such
models $M \in K$ (in suitable sense).  
Those in the other side, the ``low" one have some
``positive" theory, we have to some extent understood those models,
e.g. they have a good dimension theory.

A reader interested to see more quickly what is done rather than
why it is done and what are our hopes should go to $(C)$ below.

A good dividing line of a family of classes is such that we really
can say something on both sides,
ideally it also should help us prove things on all $K$'s by division
to cases.  So it seems advisable to prove the equivalence of an external
property (like not having many models) and an internal property (some
understanding of models of $K$).  Now clearly such a dividing line is
interesting but, of course, there are properties which are interesting
for other reasons. (See more on this in the end of (A) of \S1). 
\bn
(B) The structure/non-structure dichotomy

More specifically we may ask: which classes have a structure theory?
By a structure theory we mean ``determined up to isomorphism by an
invariant called the
dimension or several dimensions or something like that".  A
non-structure property (or theorem) will be a strong witness that
there is no structure theory.  So the question is:
\nl
\margintag{E53-nl.1.P}\ub{\stag{E53-nl.1.P} Question}:  When does a class ${\frak K}$ of models have a
structure theory?  In particular, each model from ${\frak K}$ is
characterized up to isomorphism by a ``complete set of reasonable
invariants" like those of Steinitz (for algebraically closed fields)
and Ulm (for countable torsion abelian groups).

This is still quite vague, and it takes some explanation (and choices) 
to make it concrete.  Instead we shall be even more
specific.  We shall explain two more concrete questions:  categoricity
and the main gap and the solution in the known (first order countable
vocabulary) case.  
Counting the number of models in a class seems very
natural and to make sense we have to count them in each cardinality
separately.  If the reader is not enthusiastic about this counting,
some alternative questions lead us to the same place: e.g.: having
models which are almost isomorphic but not really isomorphic (see more
in $(*)_2$ from \S1(B)(c)).
\bigskip

\definition{\stag{E53-nl.1.0} Definition}  For a class $K$ of 
models and infinite cardinal $\lambda$ let $\dot
I(\lambda,K)$ be the number of models in $K$ of cardinality $\lambda$
up to isomorphism.  So for any $K$ it is a function from Card,
the class of cardinals to itself; we may write ${\frak K} =
(K,\le_{\frak K})$ instead of $K$.

Now a priori we may get quite arbitrary functions.
But it seems reasonable to hope that all our classes ${\frak K}$ 
will have a simple function $\lambda \mapsto \dot I(\lambda,{\frak
K})$ and classes with a ``structure theory" will have such functions with
small values.  It seems more hopeful to try to first investigate 
the most extreme cases (being one and being maximal), 
considering both our chances to solve and for
getting an interesting answer; also we expect the ``upper" one to give
the important dividing lines.  It is most natural to start asking
about the spectrum of existence, i.e., being non-zero, 
i.e., what can be $\{\lambda:{\frak K}_\lambda \ne
\emptyset\}$?  This had been answered quite satisfactorily (see
\marginbf{!!}{\cprefix{88r}.\scite{88r-1.10}},\marginbf{!!}{\cprefix{88r}.\scite{88r-1.11}} above LST$({\frak K})$, it is an initial
segment with a known bound),
and it seems easier at least from the present perspective.

Considering this, the number one 
naturally has a place of honor; this is categoricity.  Recall $K$ is said
to be categorical in $\lambda$ iff $\dot I(\lambda,K)=1$.

A natural thesis is
\nl
\margintag{E53-nl.1.B}\ub{\stag{E53-nl.1.B} Thesis}:  If we really
understand when a (reasonable) class is categorical in $\lambda$ it should have
little dependence on $\lambda$, ignoring ``few, exceptional" cardinals.  
\nl
[Why?  How can we understand why ${\frak K}$ is categorical in
$\lambda$? We should know so much on the class so that given two
models from $K$ of cardinality $\lambda$ we can construct in a
coherent way an isomorphism from one onto the other; but this should work
for any other (large enough) cardinal.  Also being categorical implies
the model is a very simple one, analyzable.  

This is, of course, not
true for every class of, e.g. if $K$ is the class of $\{(I,<)$: $<$ well
orders $I$ and if $|I|$ is a successor cardinal then every
initial segment has cardinality $< |I|\}$.
This class is categorical in $\aleph_\alpha$ iff $\aleph_\alpha$ is a
limit cardinal (we could change it to ``$\alpha$ even", etc).
However, we have to restrict ourselves to ``reasonable" classes.]

An antagonist argument against the thesis \scite{E53-nl.1.B} is that
for first order $T$, the class $\{\lambda:T$ has in $\lambda$ a rigid
model, i.e., one without (non-trivial) automorphism$\}$, e.g. can be
``any class of cardinals" in some sense, 
e.g., $\{\aleph_3,\aleph_{762},\beth_{\omega_3}$,
first inaccessibly cardinality$\}$.  This class may be, essentially,
 any $\Sigma^1_2$ class of cardinals (see \cite{Sh:56}).

We may answer that rigidity implies a complicated model so we may have
$T$ coding a definition of a complicated class, of cardinals, whereas being
categorical implies the models are simple.  The antagonist may answer
that allowing enough classes of models it would not work, the categoricity
spectrum will be weird and probably \L os (see below) has no
good enough reasons for his conjecture (of course we can argue till the
problem is resolved).  We may answer that \L os conjecture implicitly
says that first order classes (of countable vocabulary) 
are ``nice", ``analyzable".  So
\scite{E53-nl.1.B} begs the question of which classes are reasonable and this
book contend that abstract elementary classes are.

Of course, there may be reasonable classes for which ``${\frak K}$ is
categorical" depends on simple properties of the cardinal 
(e.g. being strong limit).

More specifically we may ask: is it true for every (relelvant) ${\frak
K}$, either ${\frak K}$ is categorical in almost every $\lambda$ or 
non-categorical in almost
every $\lambda$?  Indeed \L os had conjectured that if an elementary class  
${\frak K}$ with countable vocabulary is categorical in one 
$\lambda > \aleph_0$ then ${\frak K}$ is categorical in every
$\lambda > \aleph_0$, having in mind the example of algebraically
closed fields of a fixed characteristic.  
A milestone in mathematical logic 
history was Morley's proof of this conjecture.  The 
solution forces you to understand such ${\frak K}$.

We may ask: Is $\dot I(\lambda,{\frak K})$ a non-decreasing function?
Of course, this is a question on $K$ but the assumptions are on
${\frak K} = (K,\le_{\frak K})$. 
This sounds very reasonable as ``having more space we have more
possibilities".  For elementary ${\frak K}$ with countable vocabulary this was
conjectured by Morley (for $\lambda > \aleph_0$).  
It is not clear how to prove it directly so it seemed to me a
reasonable strategy is to find some relevant dividing lines: the
complicated classes will have the maximal number of models, 
the less-complicated ones
can be investigated as we understand them better.  This may lead us to
look at the dual to categoricity, the other extreme - when $\dot
I(\lambda,T)$ is maximal (or just very large).
\enddefinition
\bigskip

\definition{\stag{E53-nl.1.D} Definition}  The main gap conjecture for $K$
says that either
$\dot I(\lambda,K)$ is maximal (or at least large) for almost all
$\lambda$ \ub{or} the number is much smaller for almost all $\lambda$;
for definiteness we choose to interpret ``almost all $\lambda$" as for
every $\lambda$ large enough.

(We cheat a little: see \scite{E53-nl.1.2}).

This seems to me preferable to ``$\dot I(\lambda,K)$ is non-decreasing"
being more robust; this will be even more convincing if we 
succeed in proving the stronger statement:
\sn
\margintag{E53-nl.1.1}\ub{\stag{E53-nl.1.1} The structure/non-structure Thesis}  For every 
reasonable class 
either its models have a complete set of cardinal
invariants \ub{or} its models are too complicated to have
such invariants.  

This had been accomplished for elementary
classes (= first order theories) with countable vocabularies.  We 
suggest that the main gap problem is closely connected to \scite{E53-nl.1.1}.

So ideally, for classes ${\frak K}$ with structure for every model 
$M$ of ${\frak K}$ we should be able to find a set
of invariants which is complete, i.e., determines $M$ up to
isomorphism.  Such an invariant is the isomorphism type, so we should
restrict ourselves to more reasonable ones, and the natural
candidates are cardinal invariants or reasonable generalizations of
them. E.g. for a vector space over $\Bbb Q$ we need one cardinal (the
dimension = the cardinality of any basis).  For a 
vector space over an algebraically closed field, two cardinals; (the
dimension of the vector space and the transcedence degree (= maximal 
number of algebraically independent elements) of the field, both can be any
cardinal; of course, we have also to say what the characteristic of
the field is).  
For a divisible abelian group $G$, countably many cardinals
(the dimension of $\{x \in G:px =0\}$ for each prime $p$ and the rank
of $G/\text{Tor}(G)$ where Tor$(G)$ is the subgroup consisting of the
torsion members of $G$, i.e. $\{x \in G:nx=0$ for some $n>0\}$).  
For a structure with countably many
one-place relations $P_n$ (i.e., distinguished subsets), we need
$2^{\aleph_0}$ cardinals (the cardinality of each intersection of the form
$\cap\{P^M_n:n \in u\} \cap \{M \backslash P^M_n:n \notin u\}$) for
$u$ a set of natural numbers).

We believe the reader will agree that every structure of the form 
$(|M|,E)$, where $E$
is an equivalence relation, has a reasonably complete set of
invariants: namely, the function saying, for each cardinal $\lambda$,
how many equivalence classes of this cardinality occur.  Also, if we enrich
$M$ by additional relations which relate only $E$-equivalent members and
such that each $E$-equivalence class becomes a structure with a complete set
of invariants, then the resulting model will have a complete set of
invariants.  We know that even if we allow such generalized cardinal
invariants, we cannot have such a structure theory for every relevant
class (e.g. the class of linear orders has no such cardinal invariants).  
So if we have a real dichotomy as we hope for, we should have 
a solution of (a case of) the main gap conjecture
which says each class $K$ either has such invariant or is provably
more complicated.
\bn
\beginaside
Let us try to explicate this matter.  We define what is a $\lambda$-value of
depth $\alpha$ by induction on the ordinal $\alpha$: 
for $\alpha=0$ it is a cardinal $\le \lambda$, for
$\alpha = \beta +1$ it is a sequence of length $\le 2^{\aleph_0}$ of
functions from the set of $\lambda$-values of depth $\beta$ to the set
of cardinals $\le \lambda$ or a $\lambda$-value of depth $\beta$, and
for $\alpha$ a limit ordinal it is a $\lambda$-value of some depth $<
\alpha$.

An invariant [of depth $\alpha$] for models of $T$ is a function
giving, for every model $M$ of $T$ of cardinality $\lambda$, some
$\lambda$-value [of depth $\alpha$] which depends only on the
isomorphism type of $M$.  If we do not restrict $\alpha$, the set of
possible values of the invariants is known, in some sense, to 
be as complicated as the set of all models.  
\enddefinition
\bn
This leads to:
\nl
\margintag{E53-nl.1.2}\ub{\stag{E53-nl.1.2} Main Gap Thesis}:  1) A class $K$ has a 
structure theory if there
are an ordinal $\alpha$ and invariants (or sets of invariants) of
depth $\alpha$ which determines every structure (from $K$) up to
isomorphism.
\nl
2)  If $K$ fails to have a structure theory it should have
``many" models and we expect to have reasonably definable such invariants.
\bn
We can prove easily, by induction on the ordinal $\alpha$, that
\demo{\stag{E53-nl.1.2D} Observation}  The number of $\aleph_\gamma$-values
of depth $\alpha$ has a bound $\beth_\alpha(|\tau_K|+|\gamma|)$ where

$$
\beth_\beta(\mu) = \mu + \dsize \prod_{\varepsilon < \beta}
2^{\beth_\varepsilon(\mu)}.
$$
\enddemo
\bigskip

\proclaim{\stag{E53-nl.1.3} Corollary of the thesis}  If ${\frak K}$ has a
structure theory by the interpretation of 
\scite{E53-nl.1.2} \ub{then} there is an ordinal $\alpha$ such
that for every ordinal $\gamma,{\frak K}$ has 
$\le \beth_\alpha(|\tau_{\frak K}| + |\gamma|)$ 
non-isomorphic models of cardinality $\aleph_\gamma$.
\endproclaim
\bn
It is easy to show, assuming e.g., the G.C.H., that for every
$\alpha$ there are many $\gamma$'s such that $\beth_\alpha(|\omega +
\gamma|) < 2^{\aleph_\gamma}$ and even $< \aleph_\gamma$.  
Thus, if one is able to show that
${\frak K}$ has $2^{\aleph_\gamma}$ models of cardinality $\aleph_\gamma$, this
establishes non-structure.  
\endaside

In the case in which the main gap was proved, it turns out that there 
are only a few ``reasons" for  an elementary class ${\frak K}$ with
countable vocabulary to have the maximal number of models:
\mr
\item "{$(a)$}"  ${\frak K}$ is so called unstable, 
prototypical example are the
class of infinite linear orders and the class of random graphs
[formally: in some model from ${\frak K}$ some
first order formula $\varphi(\bar x,\bar y)$ with $\ell g(\bar x) = m
= \ell(\bar y)$ for every linear order $I$ there is $M \in {\frak K}$
and an $m$-tuple $\bar a_t$ from $M$ for each $t \in I$ such
that $\varphi[\bar a_s,\bar a_t]$ is satisfied in $M$ iff $s <_I t$]
\sn
\item "{$(b)$}"  ${\frak K}$ has the so called OTOP, it is similar to
(a), but the order is defined in a different way, not by a so-called
first order formula but by a formula of the form $(\exists \bar z) \dsize
\bigwedge_n \varphi_n(\bar x,\bar y,\bar r)$.  The prototypical example is
straightforward but somewhat cumbersome.
\sn
\item "{$(c)$}"  it has the DOP, this is harder to define. 
An easy example is: two cross-cutting equivalence relations.
It means that in some members $M$ of ${\frak K}$,
we can define large linear orders by using dimensions
\nl
\beginaside
A proto-typical example is: for some infinite $I$ and $R
\subseteq I \times I,M_{I,R}$ has universe $I \cup \{(s,t,\alpha):s
\in I,t \in I,\alpha < \omega_1$ and $(s,t) \in R \Rightarrow \alpha <
\omega\}$ and relation $P^M = \{(s,t,a):a = (s,t,\alpha)$ for some
$\alpha\}$.  So $R$ can be defined in $M_{I,R}$ (though is not a relation
of $M$) as $\{(s,t)$: the set $\{x:M_{I,R} \models P(s,t,x)\}$ is
uncountable$\}$.
But the definition is not first order, it speaks on dimension 
(actually we can also interpret any graphs).
Note that $T = \text{\rm Th}(M_{I,R})$ does not depend on $R$.  
\endaside
\sn
\item "{$(d)$}"  ${\frak K}$ is so called unsuperstable;  
proto-typical example $({}^\omega I,E_n)_{n <
\omega}$  where ${}^\omega I$ is the set of functions from $\Bbb N$
into $I$ and $E_n = 
\{(\eta,\nu):\eta,\nu \in {}^\omega I$ and $\eta \restriction n = \nu
\restriction n\}$
\sn
\item "{$(e)$}"  $T$ is deep, proto-typical example is
the class of graphs which are trees (i.e. with no cycles).
\ermn
We return to the more concrete
question: the main gap and the thesis \scite{E53-nl.1.1}.
We can hope that a non-structure theorem
should imply $\dot I(\lambda,K)$ is large, whereas a structure theorem 
should enable us
to show it is small and even allow us to show it is non-decreasing,
and to compute it.
\bigskip

\beginaside
Actually the picture of the ``non-structure'' side (in the resolved
case) is more complicated.  In some classes ``reasons" (a)-(d) fail but
``reason" (e) holds.  In this case the members of ${\frak K}$
are essentially as complicated
as graphs which are trees (i.e., no cycle); for them we get the maximal
number of non-isomorphic models, but we have a ``handle" on
understanding the models.  The following result illustrates this kind
of understanding:  possibly 
\footnote{formally: if some (mild) large cardinal exists}
for some $\lambda$ we cannot find $\lambda$ models (of any
cardinality) no one embeddable into the others.  If one of clauses
(a)-(d) holds, there are
stronger results in the inverse direction (e.g. we can  code
stationary sets modulo the club filter).  So it seemed that we end
up with a trichotomy rather than a dichotomy.  That is, for the
question of counting the number of models up to isomorphism the middle
family behaves more like the high one: has maximal number.  But for the
question mentioned above and also for
questions of the form: ``are there two very similar non-isomorphic
models in the class" the middle family behaves like the low 
(e.g. we can build reasonable invariants when not restricting the 
ordinal depth).  Still
there are clear results for each of the three families.
\endaside

It was (and is) our belief that there is such a theory even for
abstract elementary classes and that we should
look at what occurs at large enough cardinals, as in small cardinals
various ``incidental" facts interfere.  Notice that a priori there
need not be a solution to the structure/non-structure problem or to
the spectrum of categoricity 
problem: maybe $\dot I(\lambda,T)$ can be any one of a
family of complicated functions, or, worse, maybe we cannot
characterize reasonably those functions, or, maybe the
question of which functions occur is independent of the usual axioms
of set theory.
\mn

Now, of course, the aim of classification is not just those specific
questions.  We rather think and 
hope that trying to solve them will on the way give
interesting dividing lines among the classes.  
A class $K$ here may have too many models but still we can say
much on the structure of its models.

Now the thesies underlining the above is
\nl
\margintag{E53-nl.2.0A}\ub{\stag{E53-nl.2.0A} Thesis}
\mr
\item "{$(a)$}"  dividing lines are interesting, and obviously
reasonable test questions are a good way to find them
(and we try to use test questions of self-interest)
\sn
\item "{$(b)$}"  good dividing lines throw light also on questions
which seem very different from the original test questions
\sn
\item "{$(c)$}"  in particular, investigating $\dot I(\lambda,K)$ (and
more profoundly, characterizing the classes with complete set of invariants)
is a good way to find interesting dividing lines, but naturally there are
other ways to arrive at them and
\sn
\item "{$(d)$}"  there are measures of complexity of a class 
(other than $\dot I(\lambda,K)$) which lead to interesting 
dividing lines and some such work was done
on elementary classes (see \S1).
\endroster
\bn
Behind the discussion above also stands
\nl
\margintag{E53-Th.CarD}\ub{\stag{E53-Th.CarD} Thesis}:   To investigate classes $K$ it
is illuminating to look for each $\lambda$, at problems on ${\frak
K}_\lambda,{\frak K}$ which is restricted to cardinal $\lambda$ and
\mr
\item "{$(a)$}"  to try to prove that the answer does not depend on
$\lambda$ or at least depends just on a small amount of information on
$\lambda$
\sn
\item "{$(b)$}"  to discard too small cardinals (essentially to look
at asymptotic behaviour)
\ermn
This seems to be successful in discovering stability (and superstability).
\sn
\beginaside
An illustration is that Rowbottom had defined $\lambda$-stable
(i.e. $A \subseteq M \wedge |A| = \lambda \Rightarrow |{\Cal S}(A,M)|
\le \lambda$) but it seems to me only
having (\cite{Sh:1}) the characterization of $\{\lambda:T$ stable in
$\lambda\}$ and the equivalence with the order property and 
defining ``$T$ stable" started stability theory.  
(Of course, for his aims this was irrelevant).
\endaside

The rationale is that if the answer is the same for ``most $\lambda$",
this points to a profound property of the class and 
it forces you to find inherent principles which you 
may not be so directly led to otherwise.
Hence it probably
will be interesting even if you care little about these cardinals.  A
parallel may be that even low dimension algebraic topologists were 
interested in the solution of Poincare conjecture for dimension $\ge 5$.  
Also the behaviour in too small cardinals may be ``incidental".
So the class of dense linear order with neither first nor last
element and the class of atomless Boolean Algebra or the class of
random enough graphs are categorical in $\aleph_0$, but have many 
complicated models in higher ones.   (One may feel these are 
low theories.  This is true by some other criterions, other test
problems; in fact, there are dividing
lines among the elementary classes for which they are low.  Still, for
the test questions considered here, provably those classes are
complicated, e.g., in a strong sense do not have a set of cardinal invariants
characterizing the isomorphism type).

You may wonder: 
\sn
\margintag{E53-nl.2.F}\ub{\stag{E53-nl.2.F} Question}:  Do we recommend dividing lines everywhere? (in
mathematics) or is this something special for model theory?

Now dividing lines are meaningful in many circumstances.  But on the
one hand it is better to list all simple finite groups than to find a dividing
line among them.  Similarly for the elementary classes categorical in
\ub{every} $\lambda \ge \aleph_0$.  On the other hand, surely for many 
directions there are no fruitful dividing lines.  The thesis that
appeared here means that for broad front in
model theory this is fruitful.  (Not everywhere: too
strong infinitary logics are out).  It seemed that this has been
vindicated for stability (and to some extent for simplicity and
hopefully for (the family of) dependent elementary classes).  

It may be helpful to compare this to
alternative approaches in model theory.  One extreme position 
will say that there is
a central core in mathematics (built around classical analysis and
geometry; and number theory of course) and other areas have to 
justify themselves by contributing
something to this central core.  Dealing with cardinals is pointless
bad taste, and while some interaction of elementary classes with
cardinals had been helpful, its time has passed.  

It seemed to me that the criterion and its application
leave out worthwhile directions.
We all know that some neighboring subjects are just hollow noise and
sometimes we are even right.  So an excellent witness for a
mathematical theory to be worthwhile is its ability to solve problems
from others, preferably classical areas or problem from other
sciences.  
Certainly a sufficient
condition.  What is doubtful is whether it is a necessary condition;
we do not agree. 

However, even within this narrow criterion, the direct attack is not
the only way to look for applications to other areas.  Not so seldom
do we find that only after developing strong enough theory, deep
applications become possible, the history of model theory seems to
support this (in particular, lately in works of Hrushovski and
Zilber).  Looking at large enough cardinals serve as asymptotic
behaviour, in which it is more transparent what are the general
outlines of the picture.

The reader may wonder how this work is related, e.g. to category theory?
universal algebra?  soft model theory?.  For category theory this work, 
in short, is closer than classical
model theory but still not really close, similarly in 
category theory each class ${\frak K}$ is equipped with a
notion of mapping (rather than $\le_{\frak K}$ 
being defined from $K$ by some specific logic as in classical
model theory).  But here we restrict ourselves to embeddings (this is not
unavoidable but things are already hard enough without this) and the main
difference is that we do not forget the elements.

What about universal algebra?  A traditional model \footnote{but no
universal algebraist agree} theorist definition of 
model theory is combining universal algebra and logic,
so a large part of this work is, by that definition, 
in universal algebra.  I do not
see any reason to disagree but still the methods and results are 
well rooted in the model theoretic tradition. 

What about soft model theory?  Though our work itself does not need
soft model theory, it fits well there (and \chaptercite{88r}, 
\chaptercite{734} use infinitary logics hence are not discussed in this part).

\beginaside
First, for many important logics ${\Cal L}$, for theories $T \subseteq
{\Cal L}(\tau)$ the class 
$(\text{Mod}_T,\prec_{{\Cal L}(\tau)})$ or variants are
abstract elementary classes (certainly for the logic 
$\Bbb L_{\lambda^+,\omega})$ and by choosing the
$\le_{\frak K}$ appropriately also $\Bbb L(\bold Q^{\text{card}}_{\ge
\lambda})$; in fact they were the original motivation to look at
abstract elementary classes.  
So if you ask for the part of soft model theory dealing with
classification theory or at least investigate categoricity, you arrive
here.  Also not just varying the logic, but fixing a class Mod$_T$
fits it well.
\endaside

This work certainly reflects the author's preference to find something
in the white part of our map, the ``terra incognita" 
rather than understand perfectly what we
have reasonably understood to begin with (which is exemplified by
looking at abstract elementary classes 
on which our maps reflect our having little to say on
them, rather than FMR theories or o-minimal theories, cases where we
had considerable knowledge and would like to complete it).
Anyhow, by experience, there will not be many complaints on lack of
generality and broadness. 

Note that we would like to get results, not consistency results and allowing
definability of well ordering or completeness runs into set-theoretic
independence results so restricting ourselves to an abtract elementary
class, a framework
which excludes well ordering and complete spaces is reasonable.
But we shall not really object to cardinal arithmetic assumptions 
like weak forms of GCH.

\beginaside
In fact, having the \ub{non-structure} results depend on the universe
of set theories is not desirable but is reasonable, as they still
witness the impossibility of a positive theory.  It is reasonable to
adopt this as part of the rules of the games.  In some cases,
consistency results forbid us to go further (see, e.g. \cite{Sh:93}).
But still the positive side should better be in ZFC.
\endaside
\bn
(C) \ub{Abstract elementary classes}

We now return to the question: With which classes of structures we
shall deal?
Obviously, ``a class of structures" is too general.  Getting down to business
we concentrate on
\mr
\item "{$\boxtimes$}"  $(a) \quad$ abstract elementary classes 
\sn
\item "{${{}}$}"  $(b) \quad$ good $\lambda$-frames
\sn
\item "{${{}}$}"  $(c) \quad$ beautiful $\lambda$-frames.
\ermn
In short, in $\boxtimes(a)$, see below, we suggest 
abstract elementary classes (a.e.c.)  as 
our framework, i.e., the family of classes we try to classify; it 
clearly covers
much ground and seems, at least to me, very natural.  What needs
justification is whether we can say on it interesting things, have non-trivial
theorems.

Among elementary (= first order) classes we know which classes have reasonable
dimension theory, the so called superstable elementary classes; and we
like to understand the case for the family of \aec. \, In 
$\boxtimes(c)$, see below in \S3(C), we suggest 
beautiful $\lambda$-frames as our ``promised land", as a context where
we have reasonable understanding, e.g., have
dimension theory, can prove the main gap, etc. (but of course  more
wide families ``on our way" probably will be interesting per se).
Now it is very unsurprising that if we assume enough axioms, we shall
regain paradise (which means here quite full fledged analog to the so
called superstability
theory, at least for my taste).  Hence the problem in justifying the
choice of $\boxtimes(c)$ is mainly not in pointing to many good properties but
to show that there are enough such frames and it helps
prove theorems not mentioning it.  On the second (i.e. prove theorems
not mentioning them...), see
e.g. \scite{E53-nl.2.7}(2) below.  In our context ideally the first (i.e. ``there
are enough such frames...") means to show that they are the 
only ones, i.e., the
broadest family of \aec \, which has such a good dimension theory.  We are
far from this, still we would according to our ``guidelines" 
like at least to get beautiful frames by choosing to consider the
classes which fall on the ``low" side (in the
elementary classes case) by dividing lines (= dichotomies) inside a family of
classes which is large and natural, here among abstract elementary
classes.  That is, 
the program is to suggest some dividing lines, for the high side to
prove the so-called non-structure theorems and for the low side to have
some theory.  Being always in the low sides we should 
arrive to beautiful frames.

But most of our work falls under $\boxtimes(b)$, good $\lambda$-frames.  So it
needs double justification: on the one hand we have to show it arises
naturally from our program.  
\nl
[In detail, a weak case for ``arising naturally" is to start with an
abstract elementary classes 
satisfying some external condition of being ``low" like
categoricity, and prove that ``inside ${\frak K}$" we can find good
frames.  A strong case is to find a dividing line such that for each
low ${\frak K}$ we can find inside it ``enough" good frames, and
for all other ``few".  There is another meaning of ``arising naturally"
which would mean that we have looked at some natural examples and
extracted the definition from their common properties; this is not what
we mean.  We rather try to solve questions on the number of models but
of course the first order case was before our eyes as first
approximation to the paradise we would like to arrive to.]

On the other
hand for such frames, possibly with more assumptions justified
similarly we can say something significant.

\beginaside  
In fact, we see good $\lambda$-frames essentially
as the rock-bottom analogs of the family of elementary classes called
superstable mentioned above.
\endaside

We shall discuss $\boxtimes(a)$ and $(b)$ and $(c)$ in more detail.
We start with
\mr 
\item "{$\boxtimes(a)$}"  abstract elementary classes.  
\ermn
Recall the definition of abstract
elementary classes Definition \scite{E53-nl.0.7}.
\bn
\margintag{E53-nl.2.1}\ub{\stag{E53-nl.2.1} Explanation}:  An \aec \, is easy to 
explain (probably much simpler than elementary (=
first order) class).  
Such ${\frak K}$ consists of a class $K$ of structures = models, all of the
same ``kind", e.g. all rings have the same kind, but a group has a
different kind.  We express this by saying ``all members of $K$ has
the same vocabulary $\tau = \tau_K$".  E.g., $K$ consists of objects 
of the form $M = (A^M,F^M_0,F^M_1,Q^M),A^M$ its universe, a
non-empty set, $F^M_\ell$ a binary function on it, $Q^M$ a binary
relation.  ${\frak K}$ has also an order $\le_{\frak
K}$ on $K$, its notion of being a sub-structure (which refines the
standard notion).  Now $(K,\le_{\frak K})$  have to satisfy
some requirements: preservation under isomorphisms, $\le_{\frak K}$ being
an order, preserved by direct limits and also direct limits inside $N
\in K$, remembering that our mapping are embedding.  Also 
if $M_1 \subseteq M_2$
are both $\le_{\frak K}$-sub-structures of $N$ then $M_1 \le_{\frak K}
M_2$, and lastly we demand every $M \in K$ has a countable $\le_{\frak
K}$-sub-structure including any pregiven countable set of elements (or
replace countable by a fix cardinality, we ignore this point in the
introduction; see \sectioncite[\S1]{600}).
\bn
\beginaside
Concerning ``$M_\ell \le_{\frak K} N,(\ell=1,2),M_1 \subseteq M_2
\Rightarrow M_1 \le_{\frak K} M_2$" note that if we define 
$\le_{\frak K}$ as $\prec_{\Cal L}$ for any logic, this will hold.
\bn
\endaside
For elementary classes ${\frak K}$, because of the so-called compactness and
L\"owenheim-Skolem-Tarski theorems, the situation in all cardinals is to a
significant extent similar.  

In particular, if ${\frak K}$ is an elementary class (with countable
vocabulary) and $\lambda_1,\lambda_2$ are
(infinite) cardinals then there is $M \in {\frak K}$ of cardinality
$\lambda_1$ iff there is $M \in {\frak K}$ of cardinality
$\lambda_2$.  So recalling that $K_\lambda = \{M \in K:M$ has cardinality
$\lambda\}$ and ${\frak K}_\lambda = (K_\lambda,\le_{\frak K}
\restriction K_\lambda)$ we have $K_{\lambda_1} \ne \emptyset
\Leftrightarrow K_{\lambda_2} \ne \emptyset$.
Moreover, any infinite $M \in {\frak
K}$ has $\le_{\frak K}$-extension in every larger cardinality.  
But for abstract elementary classes it is not
necessarily true, and even if $(\forall \lambda) K_\lambda \ne
\emptyset$ there may be many $\le_{\frak K}$-maximal models, i.e., $M
\in K$ such that $M \le_{\frak K} N \Rightarrow M=N$.  This (and more)
makes the theory very different.

The context of \aec \, may seem so general, we 
may doubt if anything interesting can be
said about it; still note that this context does not allow the class of Banach
spaces as the union of an increasing chain is not necessarily
complete.  Certainly a loss.  Also the 
class $(W,\subseteq)$, the class of well orders, is not an \aec;
(recall $I$ is a well order if it is a linear order such that every
non-empty set has a first element).  
Similarly the class $(K^{\text{fgi}},\subseteq)$ where 
$K^{\text{fgi}} =$ the class of
rings (or even integral domains) in which every 
ideal is finitely generated, is not an \aec (where $\le_{\frak K}$ is
being a subring).  However, we get an \aec \, when we consider only
$K_{\le n} =$ the class of rings in which every ideal is
generated by $\le n$ elements.
\nl
\beginaside
\nl
We may like to replace $n$ by a countable ordinal $\alpha$, i.e.,
$K^{\text{fgi}}_{\le \alpha} = 
\{M \in K: \text{\rm dp}_M(\emptyset) \le \alpha\}$; where
for a ring $M$ we define dp:$\{u:u \subseteq M$ finite$\}
\rightarrow$ the ordinals by dp$_M(u) = \cup\{\text{\rm dp}(w)+1:u
\subseteq w$ and $w$ is not included in the ideal of $M$ which $u$
generates$\}$. But then we have problems with closure under unions; a
reasonable remedy is to have an appropriate $\le_{\frak K}$:
$M \le_{\frak K} N$ if $M,N$ are rings and for every finite $u
\subseteq M$ we have dp$_N(u) = \text{\rm dp}_M(u)$.
\nl
Why have we restricted ourselves to ``countable $\alpha$"?  Only
because in clause (g) of Definition \scite{E53-nl.0.7} we have used
``countable".
\nl
\endaside
\nl

But the family of abstract elementary classes includes all the examples 
listed in \scite{E53-nl.0.6X} in the beginning (of this section, 2).
\beginaside
Also, other abstract elementary classes are $(K,\prec)$ where $K$ is the 
class of locally finite models of a first order theory $T$.
Another example is (Mod$_\psi,\prec_{\Cal L})$ where $\psi$ is a 
sentence from logic $\Bbb L_{\lambda^+,\omega}$ with ${\Cal L}$
the set of subformulas of $\psi$.  Also $(K,\prec)$ where $P
\in \tau_{\frak K}$ is a unary predicate, 
$T$ first order and $K = \{M \in \text{\rm Mod}_T:
P^M = \Bbb N$, the natural numbers$\}$.
\endaside
\newline

A natural property to consider is amalgamation.  We say that ${\frak
K}$ has the amalgamation property when for any $M_\ell \in {\frak
K},\ell=0,1,2$ and $\le_{\frak K}$-embedding $f_1,f_2$ of $M_0$ into
$M_1,M_2$ respectively (this means that $f_\ell$ is an isomorphism
from $M_0$ onto some $M'_\ell \le_{\frak K} M_\ell$) there are $M_3
\in {\frak K}$ and $\le_{\frak K}$-embeddings $g_1,g_2$ of $M_1,M_2$
into $M_3$ respectively such that $g_1 \circ f_1 = g_2 \circ f_2$.
Should we adopt it?
Now it is a very important property, we would
like to have it, but it is a strong restriction (our prototyical
problem, models of $\psi \in \Bbb L_{\omega_1,\omega}$ fails it); 
so we do not assume it, but it will appear as a dividing line.
\nl
So the thesis is
\nl
\margintag{E53-nl.2.2}\ub{\stag{E53-nl.2.2} Thesis}:
\mr
\item "{$(a)$}"  In the context of abstract elementary classes we 
can answer some
non-trivial questions
\sn
\item "{$(b)$}"  In particular we can 
say something on the categoricity spectrum
\sn
\item "{$(c)$}"  In the long run a parallel to the main gap
will be found.
\ermn
A reasonable reader may require an example of results.  First we quote
\cite{Sh:576} represented here in \cite{Sh:E46}:
\proclaim{\stag{E53-nl.2.3} Theorem}  Assume $2^{\aleph_\alpha} 
< 2^{\aleph_{\alpha +1}} < 2^{\aleph_{\alpha+2}}$ and ${\frak K}$
is an abstract elementary class categorical in $\aleph_\alpha$, 
in $\aleph_{\alpha +1}$ and has
an ``intermediate" number of models in $\aleph_{\alpha +2}$, 
\ub{then} ${\frak K}$ has at least one model in $\aleph_{\alpha +3}$.
\endproclaim
\bn
Note that
\demo{\stag{E53-nl.2.3E} Notation}  If $\lambda = \aleph_\alpha$ we let
$\lambda^{+n} = \aleph_{\alpha +n}$, so can write this theorem in such
a notation, similarly later.
\enddemo
\bn

So it is an example for \scite{E53-nl.2.2}(a)+(b): not ``every function"
can occur as $\lambda \mapsto \dot I(\lambda,{\frak K})$.

Note that this theorem gives a weak conclusion, but with very weak
assumptions.  In fact at first glance it seems we are facing a wall:
our assumptions are so weak to exclude all possible relevant 
methods of model theory, in particular all relatives of compactness. 
\nl

\beginaside
I.e., we have no compact
(even just $\aleph_0$-compact) logic defining our class.  Of course, the
upward LST cannot be used, it does not make sense: the 
desired conclusion is a weak form of
it.  As for the downward L\"owenheim Skolem-Tarski theorem, with only three
cardinals available it seems to say very little.

We do not have formulas hence no types and no saturated models.  Here we
cannot use versions of ``well ordering is undefinable" as in previous
cases (see \chaptercite{88r}; if 
$\aleph_\alpha = \aleph_0$ and ${\frak K}$ is reasonable we have used 
``no $\psi \in \Bbb L_{\omega_1,\omega}(\bold Q)$ defines well
ordering (in a richer vocabulary)"; this does not apply in
\cite{Sh:576}, i.e. \cite{Sh:E46}, even when $\lambda = \aleph_0$, as 
we demand 
only LST$({\frak K}) \le \aleph_0$ rather than ``${\frak K}$ is a
PC$_{\aleph_0}$-class"; and we certainly like to allow any
$\aleph_\alpha$).  Also in general we cannot
find Ehrenfeuch-Mostowski models 
(another way to say well orders are not definable).
Also we do not assume the existence of relevant so called
large cardinals, e.g. ${\frak K}$ is definable in some $\Bbb
L_{\kappa,\omega},\kappa$ a compact or just a measurable cardinal.
So indeed no remnants of compactness are available here.
\endaside

The proof of \scite{E53-nl.2.3} leads us to our 
second framework, good $\lambda$-frames, which
has a crucial role in our investigations, see below.  The
main neatly stated result in \chaptercite{600} (part (1) of \scite{E53-nl.2.7}),
\chaptercite{705}(part (2) of \scite{E53-nl.2.7}) is:
\nl
\beginaside 
(omitting a
weak set theoretic assumption which will be eliminated in the full
version of \cite{Sh:838}).
\endaside
\bigskip

\proclaim{\stag{E53-nl.2.7} Theorem}  Assume ${\frak K}$ is an \aec . 
\nl
1) ${\frak K}$ has a member in $\aleph_{\alpha +n+1}$ \ub{if} ($n \in
\Bbb N$ and)
\mr
\item "{$(a)$}"  $n \ge 2$ and $2^{\aleph_\alpha} <
2^{\aleph_{\alpha+1}} < \ldots < 2^{\aleph_{\alpha +n}}$
\sn
\item "{$(b)$}"  ${\frak K}$ is categorical in $\aleph_\alpha$ and in
$\aleph_{\alpha +1}$
\sn
\item "{$(c)$}"  ${\frak K}$ has a model in $\aleph_{\alpha +2}$
\sn
\item "{$(d)$}"  $\dot I(\aleph_{\alpha +m},{\frak K})$ is not too
large for $m=2,\dotsc,n$.
\ermn
2) If (a)-(d) holds for every $n$ \ub{then} ${\frak K}$ is categorical
in every $\aleph_\beta \ge \aleph_\alpha$.
\endproclaim
\bn
\beginaside
Actually in this theorem ``${\frak K}$ having
L\"owenheim-Skolem-Tarski number $\le
\lambda$" (rather than $\aleph_0$) is enough.
\endaside
\bn
(D) \ub{Toward Good $\lambda$-frames (i.e. $\boxtimes(b)$}:
\mn
\margintag{E53-nl.3.1}\ub{\stag{E53-nl.3.1} Thesis}  Good $\lambda$-frames are a right context to
start our ``positive" structure theory.
\sn
\beginaside
They are a rock-bottom parallel of superstable elementary classes.
\endaside
\sn
Now compared to abstract elementary classes, much more has to be 
said in order to explain what they
are and how to justify them.  We describe good $\lambda$-frames ${\frak s}$ 
in several stages.  We need several choices to specify our context.
Usually in model theory we fix an elementary class ${\frak K}$ and
consider $M \in {\frak K}$.  Here we concentrate on one cardinal
$\lambda$, that is, we usually investigate ${\frak K}_\lambda =
(K_\lambda,\le_{{\frak K}_\lambda})$ where $K_\lambda = \{M \in
K:M$ has cardinality $\lambda\}$ and $\le_{{\frak K}_\lambda}$ is
defined by $M \le_{{\frak K}_\lambda} N$ iff $M \le_{\frak K} N,M \in
K_\lambda$ and $N \in K_\lambda$.  This is not a clear cut deviation,
also for elementary classes we sometimes fix $\lambda$, and here we
usually look at least at ${\frak K}_\lambda$ and ${\frak K}_{\lambda^+}$
together, still the flavour is different.  So (the notion ``choice"
may be seemingly problematic but a better alternative was not found).
\bn
\margintag{E53-3.1}\ub{\stag{E53-3.1} Choice}:  We concentrate on ${\frak K}_\lambda$, an
\aec \, restricted to one cardinal.

This seems reasonable because as noted above, 
transfer from one cardinal to another is central, but in our context
quite hard, so we may know various ``good" properties only around $\lambda$.  
Also there are ${\frak K}$ which in some cardinals are model
theoretically ``very simple" but in other (e.g. larger) 
cardinals complicated, and we may like to say what we can say 
about ${\frak K}_\lambda$ in $\lambda$ for which ${\frak K}_\lambda$
is ``simple".
\bn
\margintag{E53-3.2}\ub{\stag{E53-3.2} Choice}:  We concentrate here on ${\frak K}_\lambda$
with amalgamation and the JEP (joint embedding properties).  

But is amalgamation not a very strong/positive property?  Yes, \ub{but}
amalgamation for models of cardinality $\lambda$ only is much weaker
and its failure in some reasonable circumstances leads 
to non-structure results,
so it can serve as a dividing line.  More specifically, we 
know that if ${\frak K}$ is categorical in $\lambda \ge
\text{ LST}({\frak K})$ and ${\frak K}_\lambda$ fails amalgamation
and ${\frak K}_{\lambda^+} \ne \emptyset$ then in 
${\frak K}_{\lambda^+}$ we have many complicated models 
(provided that $2^\lambda < 2^{\lambda^+}$; see \chaptercite{88r}).
\bn
\margintag{E53-3.3}\ub{\stag{E53-3.3} Choice}:  In ${\frak K}_\lambda$ there is a
\ub{superlimit} model $M^*$ which means that: $M^* \in {\frak K}_\lambda$
is \ub{universal}, (i.e., any $M' \in {\frak K}_\lambda$ can be
 $\le_{\frak K}$-embedded into it), has a 
proper $<_{\frak K}$-extension and if $M$ is the
union of a $<_{\frak K}$-increasing chain of models isomorphic to
$M^*$ and $M$ is of cardinality $\lambda$, then $M$ is isomorphic to
$M^*$.

Can we give a natural example of a superlimit model?  For the abstract
elementary class of
linear orders, the rational order $(\Bbb Q,<)$ is superlimit (in
$\aleph_0$).  However, this is somewhat misleading as in larger
cardinals it is much ``harder", in fact, for the \aec \, of linear orders
there is no superlimit model in $\lambda > \aleph_0$.  By categoricity
the \aec \, of
algebraically closed fields of some fixed character has a superlimit model
in every $\lambda \ge \aleph_0$.  The class of
$\{(A,E):E$ an equivalence relation on $A\}$ is a bit more
informative.  Easily $(A,E)$ is
superlimit in it iff the number of $E$-equivalence classes as well as
the cardinality of each $E$-equivalence class is the number of
elements of $A$.

Of course, if ${\frak K}$ is categorical in $\lambda$ then every $M
\in {\frak K}_\lambda$ is superlimit (if it is not 
$\le_{\frak K}$-maximal in which case every $M \in {\frak K}$ has
cardinality $\le \lambda$), but having a superlimit is a much
weaker condition and it seems a right notion of generalizing
superstability (or, probably, a good first approximation).  
This may surely look tautological in view of
Definition \scite{E53-nl.0.10}(3), but that definition is misleading.  There
are several properties, which for elementary classes are equivalent to
being superstable, and we have chosen the existence of superlimit.
However, so far the existence of a superlimit model in $\lambda$
has few consequences.

Why the choice?  As this is an exterior way to say that our class is
``simple, low"; it is weaker than categoricity and we next demand much more.
\sn
\beginaside
Note that if ${\frak K}$ is an elementary class and $\lambda =
\lambda^{\aleph_0} + |\tau_{\frak K}|$ or $\lambda \ge \beth_\omega + 
|\tau_{\frak K}|$, \ub{then} $M \in K_\lambda,M$ is superlimit iff
$M$ is saturated and the theory is superstable; see \cite[3.1]{Sh:868}.
\endaside
\sn
Now we are very interested in the existence of something like
``free amalgamation", which in our context will be called non-forking
amalgamation.  That is, we are interested in saying when 
``$M_1,M_2$ are freely amalgamated over $M_0$ inside
$M_3$" (all in ${\frak K}_\lambda$).  In our main example we have to
use a more restrictive notion, having quadruples 
$(M_0,M_1,a,M_3)$ is non-forking where
$M_0 \le_{\frak K} M_1 \le_{\frak K} M_3,a \in M_3 \backslash M_1$.
This says that ``inside $M_3$ the element $a$ and the model $M_1$ 
are freely amalgamated over $M_0$".  (Mainly in \cite{Sh:576},
i.e. \cite{Sh:E46}, we use so called ``minimal types", which give rise to
such quadruples). 

This leads us to define a central notion here: $\ortp_{\frak
K}(a,M,N)$ denotes the ``orbit" of $a \in N$ over $M \le_{\frak K} N$.  
We express $(M_0,M_1,a,M_3)$ is non-forking also as
$``\nonfork{}{}_{}(M_0,M_1,a,M_3)"$ and also as
``$\ortp_{\frak s}(a,M_1,M_3)$ does not
fork over $M_0$" because it is analogous to the non-forking in first order
model theory.  But this background is not needed, as non-forking
is an abstract, axiomatic relation in our context.
\sn
\beginaside
This replaces here the notion of type in the investigation of
elementary (= first order) classes.  But there the types are defined
as $\sftp(\bar a,A,N) = \{\varphi(\bar x,\bar b):\bar b \subseteq
A,\varphi(\bar x,\bar y)$ is a first order formula and $N \models
\varphi[\bar a,\bar b]\}$.  Note: the case $A$ is the universe of $M
\le_{\frak K} N$ is not excluded but is not particularly
distinguished.  In fact, it was unnatural there to make the
restriction as there are theorems using our ability to restrict the
type to any subset of $A$ (e.g. for inductive proof) and it is
important to have results on any $A$.  
\endaside

We let ${\Cal S}_{{\frak K}_\lambda}(M) = 
\{\ortp_{{\frak K}_\lambda}(a,M,N):M 
\le_{{\frak K}_\lambda} N$ and $a \in N\}$ be called the
set of types over $M$.
The set of axioms (i.e., Definition \marginbf{!!}{\cprefix{600}.\scite{600-1.1}}) of good
$\lambda$-frames expresses the intuition of ``non-forking" as a free
amalgamation (in fact we are allowed to restrict the non-forking
relation to types 
$\ortp_{\frak s}(a,M_1,M_3)$ which are, so called basic ones, they
should mainly be ``dense" enough).  We
may consider these axioms per se, but we feel obliged to find evidence 
of their naturality of the form indicated above.   So
\bigskip

\definition{\stag{E53-3.4} Definition}  A good $\lambda$-frame 
${\frak s}$ consists of 
\mr
\item "{$(a)$}"  an \aec \, ${\frak K} = 
{\frak K}^{\frak s}$ and let\footnote{Note ${\frak K}^{\frak s}$ may
have models in many cardinals, whereas ${\frak K}_{\frak s}$ has models in
 only one cardinal} ${\frak K}_{\frak s} = {\frak K}_\lambda$ with
 LST$({\frak K}^s) \le \lambda$
\sn
\item "{$(b)$}"   for $M \in {\frak K}_\lambda$
we have ${\Cal S}^{\text{bs}}_{\frak s}(M)$, a subset of 
${\Cal S}_{{\frak K}_\lambda}(M)$ with LST$(K^{\frak s}) \le \lambda$
\sn
\item "{$(c)$}"   a notion of 
``$p \in {\Cal S}^{\text{bs}}(M_2)$ does not fork over 
$M_1 \le_{{\frak K}_\lambda} M_2$" satisfying some reasonable axioms.
\ermn
How does this help us in proving Theorem \scite{E53-nl.2.7}?  Relying on
the main results of \cite{Sh:576}, \cite{Sh:E46}, using the assumption
of \scite{E53-nl.2.7} we in 
\sectioncite[\S3]{600} prove that there is a good
$\lambda^+$-frame ${\frak s}$ with ${\frak K}_{\frak s} = 
{\frak K}_{\lambda^+}$.  Also in \sectioncite[\S3]{600} using a 
similar theorem from \chaptercite{88r} for the
case $\lambda = \aleph_0$ with a little different assumptions, we get a good
$\aleph_0$-frame ${\frak K}$.

We take a spiral approach: we look at a good
$\lambda$-frame ${\frak s}$, suggest a question, i.e., dividing lines,
if ${\frak s}$ falls under the complicated
side we prove a non-structure theorem.  If not, we know some things
about it and we can continue to investigate it, after we have enough knowledge
we ask another question.  In \sectioncite[\S5]{600} we start with a good
$\lambda$-frame, gain some knowledge and if there are not enough
essentially unique amalgamations we get many complicated models in
$\lambda^{++}$.  If ${\frak s}$ avoids this, we call it weakly
successful and  understand ${\frak K}_{\frak s}$ better.
In particular, we define the 
promised ``$M_1,M_2$ are non-forking amalgamated over $M_0$
inside $M_3$", we call this relation 
NF = NF$_\lambda$ = NF$_{\frak s}$ and prove that
it has the properties hoped for.  Listing its desired properties, it is unique.
But this has a price: we have to restrict ${\frak K}_{\frak s}$ to
isomorphic copies of the superlimit models.  After showing that if $S$
has a second non-structure property there are again many models in
$\lambda^{++}$, we are ``justified" in assuming ${\frak s}$ fails also
this non-structure a property.  We then succeed to find for
$\lambda^+$ another good frame, ${\frak s}^+$ such that
$K^{{\frak s}^+}_{\lambda^+} \subseteq K^{\frak s}_{\lambda^+}$.
Recall $K^{\frak s}$ is the a.e.c. lifting $K_{\frak s}$ and $K^{\frak
s}_{\lambda^+} = (K^{\frak s})_{\lambda^+}$.

What have we gained?  Have we not worked hard just to find ourselves in the
same place?  Well, ${\frak s}^+$ is a
good $\lambda^+$-frame and $\dot I(\mu,K^{{\frak s}^+}) \le \dot
I(\mu,K^{\frak s})$ for every $\mu \ge \lambda^+$ and
\mr
\item "{$(*)$}"   for every $\chi$ and good $\chi$-frame 
${\frak t},K^{\frak t}$ has models of cardinality $\chi^+$ and 
moreover of cardinality $\chi^{++}$.  
\ermn
So this is enough to prove the Theorem \scite{E53-nl.2.7}(1), by induction
on $n$.

Let us compare this to \cite{Sh:87a}, \cite{Sh:87b}.  There in 
stage $n$ we have
some knowledge on models in ${\frak K}_{\aleph_\ell}$ for $\ell \le n$
but our knowledge decreases with $\ell$.  Now (all in \cite{Sh:87b}) 
dealing with $n+1$ we have to consider
a question on models of cardinality $\lambda = \aleph_0$, for which our
specific tools for $\aleph_0$ (the omitting types theorem and the
assumption that ${\frak K}$ is $(\text{Mod}_\psi,\prec)$ where 
$\psi \in \Bbb L_{\omega_1,\omega})$ enable us to have proved a 
dichotomy, each side implied additional information concerning
$\aleph_\ell$ for $\ell \le n$, again decreasing with $\ell$.

\beginaside
[We elaborate: for each $\ell < n$ we can define so called full
stable $({\Cal P}^-(m),\aleph_\ell)$-systems $\langle M_u:u \in {\Cal
P}^-(m)\rangle$ for $m \le (n-\ell)$ where ${\Cal P}^-(m) = \{u:u
\subset \{0,\dotsc,m-1\}\}$.  So our knowledge ``decreases" with $\ell$:
we can handle only systems of lower ``dimension".  We ask on such
systems whether we can find suitable $M_{\{0,\dotsc,n-1\}}$.  Is it 
weakly unique (up to embedding)?  Is it unique?  Is there a prime one?
We can transfer up a
positive property from $({\Cal P}^-(m),\aleph_\ell)$ to
$({\Cal P}^-(m-1),\aleph_{\ell+1})$, and also negative ones if
$2^{\aleph_\ell} < 2^{\aleph_{\ell +1}}$.  A crucial point is the
existence of a strong dichotomy in the cardinality $\aleph_0$, 
either we have a prime 
solution or we have $2^{\aleph_0}$ pairwise incompatible ones.
\nl
Note that in \cite{Sh:87a}, \cite{Sh:87b}, we deal with types
as in elementary classes (i.e. as set of formulas) but only 
over models or $\cup\{M_u:u \in
{\Cal P}^-(n)\}$ when $\langle M_u:u \in {\Cal P}^-(n)\rangle$ is so
called stable.]
\endaside
\bn
The proofs of \chaptercite{600} seem neater than \cite{Sh:87a},
\cite{Sh:87b}: because we are ``poorer", we do
not have the special knowledge on the first $\lambda$.  So we do not have to
look back, we can forget ${\frak s}$ when advancing to ${\frak s}^+$.
This is nice for its purpose but suppose that we have a good
$\lambda^{+n}$-frame ${\frak s}^n$ for $n < \omega,{\frak s}^{n+1}$
being gotten from ${\frak s}^{+n}$ as above.  For this purpose, 
forgetting the past costs
us the future - we cannot say anything on models of cardinality $\ge
\lambda^{+ \omega}$.  This is rectified in \chaptercite{705}.
\bn
\beginaside
So in \chaptercite{705} we investigate the ${\frak K}_{{\frak s}^{+n}}$
for every $n$ large enough.  A priori
it is fine to do this for $n \ge 756$, and increasing the number as
 we continue to investigate.  But in spite of this knowledge,
considerable effort was wasted on small $n$, i.e., 
assuming little on ${\frak s}$, and in
\sectioncite[\S2-\S11]{705} we get the theory of prime, independence,
dimension, regular types and orthogonality we like (see, maybe,
\cite{Sh:F735} on what we really need to assume).  

But for going up we need to deal with ${\Cal P}^-(n)$-amalgamation - their
existence and uniqueness.  Then we can go up, see \sectioncite[\S12]{705}.
\endaside
\enddefinition
\newpage

\head {\S3 On Good $\lambda$-frames} \endhead  \resetall \sectno=3
 \spuriousreset
\bn
This continues \S2 and should be ``non-logician friendly" too, though
it may well be more helpful after some understanding/reading of the
material itself.
\bn
(A) Getting a good $\lambda$-frame
\sn
We try below to describe in more details the proof of Theorem
\scite{E53-nl.2.7}(1) + (2) proved in \chaptercite{600}, \chaptercite{705}, so we somewhat
repeat what was said before in (D) of \S2.  We have to start by 
getting good $\lambda$-frames.  We
could have concentrated on the case $\lambda = \aleph_0$ and rely on
\chaptercite{88r}, but as this does not fit the ``for non-logicians" we
instead rely on \cite{Sh:576}, \cite{Sh:603}, that is on \cite{Sh:E46}
and the non-structure from \cite{Sh:838}, at least the ``lean"
version.
\sn
\beginaside
For presentation we cheat a little in the non-structure part, saying
we prove results like $\dot I(\mu^{++},{\frak K}) = 2^{\mu^{++}}$ when ${\frak
K}$ satisfies some ``high" property and say $2^{\mu^+} <
2^{\mu^{++}}$.  One point is that this 
relies on using an extra set theoretic assumption on $\mu^+$:
the weak diamond ideal on $\mu^+$ not being
$\mu^{++}$-saturated.  This is a very weak assumption, it is not clear
whether its failure is 
consistent when $\mu \ge \aleph_1$ and in any case its failure 
has high consistency strength, that is, if the ideal is
$\mu^{++}$-saturated then there are inner models with quite large
cardinals.  We can eliminate this extra set theoretic assump[tion as
done in \cite{Sh:838} (see later part of the introduction).  The 
second point is we prove 
only that there are $\ge 
\mu_{\text{unif}}(2^{\mu^{++}},2^{\mu^+})$ many non-isomorphic models in
$\mu^{++}$.  This number is always $>2^{\mu^+}$
(recall we are assuming $2^{\mu^+} < 2^{\mu^{++}}$), and is equal to
$2^{\mu^{++}}$ when $\mu \ge \beth_\omega$ and conceivably the
statement ``$2^{\mu^+} < 2^{\mu^{++}} \Rightarrow
\mu_{\text{unif}}(2^{\mu^{++}},2^{\mu^+}) = 2^{\mu^{++}}$" is provable in ZFC.

Of course, below LST$({\frak K}) \le \lambda$ suffices instead of
LST$({\frak K}) = \aleph_0$.
\endaside
\sn
So first assume
\mr
\item "{$\boxdot_1$}"  ${\frak K}$ is an abstract elementary class, 
and for simplicity $2^\lambda < 2^{\lambda^+} < \ldots < 2^{\lambda^{+n}} <
2^{\lambda^{+n+1}} < \ldots,{\frak K}$ is categorical in
$\lambda,\lambda^+$, has a model in $\lambda^{++}$, and $\dot
I(\lambda^{+2},{\frak K}) < 2^{\lambda^{+2}}$.
\ermn
We shall now describe how to get a good $\lambda$-frame (or
$\lambda^+$-frame) from this assume, but it takes some time.
We can deduce that ${\frak K}_\lambda$ and ${\frak K}_{\lambda^+}$ have
amalgamation.
\nl
(Why?  Otherwise it has many complicated models in
$\lambda^+,\lambda^{++}$, respectively).  Now we consider the class
$K^{3,\text{na}}_\lambda$ of 
triples $(M,N,a),M \le_{{\frak K}_\lambda} N,a \in N
\backslash M$ with the (natural) order which is $(M_1,N_1,a_1) \le
(M_2,N_2,a_2)$ iff $a_1 = a_2$ (yes! equal) and $M_1 \le_{{\frak
K}_\lambda} M_2$ and $N_1 \le_{{\frak K}_\lambda} N_2$.

We may look at them as representing the ``orbit (or type of) 
$a$ over $M$ inside $N,\ortp_{\frak K}(a,M,N)$", which is not 
defined by formulas but by
mappings, (i.e. types are orbits over $M$) so if $M \le_{{\frak
K}_\lambda} N_\ell$ and $a_\ell \in N_\ell \backslash M$ then 
$\ortp_{{\frak K}_\lambda}(a_1,M,N_1) 
= \ortp_{\frak K}(a_2,M,N_2)$ iff for some
$\le_{{\frak K}_\lambda}$-extension $N_3$ of $N_2$ there is a
$\le_{{\frak K}_\lambda}$-embedding $h$ of $N_1$ into $N_3$ over $M$
which maps $a_1$ to $a_2$, recalling ${\frak K}_\lambda$ has amalgamation.

Why do we consider $K^{3,\text{na}}_\lambda := 
\{(M,N,a):M \le_{{\frak K}_\lambda}
N,a \in N \backslash M\}$ instead of 
${\Cal S}^{\text{na}}_{{\frak K}_\lambda}(M) := 
\{\ortp_{{\frak K}_\lambda}(a,M,N):(M,N,a) \in K^{3,\text{na}}_\lambda\}$?  
(The types $\ortp_{{\frak K}_\lambda}(a,M,N)$ when $a \in M$ are called 
algebraic (and na stands for non-algebraic) and are trivial, so 
${\Cal S}^{\text{na}}_{{\frak K}_\lambda}(M)$ is the rest.)  
Now ${\Cal S}_{{\frak K}_\lambda}(M)$ is very important and for 
$M_1 \le_{{\frak K}_\lambda} M_2,p \in {\Cal S}_{{\frak
K}_\lambda}(M_2)$ we can define its restriction to $M_1,p
\restriction M_1 \in {\Cal S}_{{\frak K}_\lambda}(M_1)$, with some natural
properties, and this mapping is onto (= surjective) as ${\frak K}_\lambda$ has
the amalgamation property.  But it is not clear that
an increasing sequence of types of length 
$\delta < \lambda^+$ of types has a bound (when cf$(\delta) >
\aleph_0$); see Baldwin-Shelah \cite{BlSh:862}.
For $K^{3,\text{na}}_\lambda$ this holds.  That is, if the sequence
$\langle(M_\alpha,N_\alpha,a_\alpha):\alpha < \delta\rangle$ is
increasing in $K^{3,\text{na}}_\lambda$, so 
$\alpha < \delta \Rightarrow a_\alpha
= a_0$, then it has a lub: the triple $(\cup\{M_\alpha:\alpha <
\delta\},\cup\{N_\alpha:\alpha <\delta\},a_0)$.

Some types (and triples) are in some sense better understood: here the
ones representing minimal types; where
\mr
\item "{$(*)$}"  $p \in {\Cal S}^{\text{na}}_{{\frak K}_\lambda}(M)$
is minimal if for every $\le_{{\frak K}_\lambda}$-extension $N$ of $M$
the type $p$ has at most one extension in 
${\Cal S}^{\text{na}}_{{\frak K}_\lambda}(N)$.
\ermn
Note that $p$
always has at least one  extension in ${\Cal S}_{{\frak K}_\lambda}(N)$ by 
amalgamation and we can prove that $p$ has at least one
from ${\Cal S}^{\text{na}}_{{\frak K}_\lambda}(N)$ in our context,
and recall that we have discarded the algebraic types, i.e. those of
 $a \in M$.
\bn
It is too much to expect that every $p \in 
{\Cal S}^{\text{na}}_{{\frak K}_\lambda}(M)$ is minimal, but what about
\nl
\margintag{E53-qu.7.1}\ub{\stag{E53-qu.7.1} Question}:  Is the class of minimal types dense,
i.e., for every $p_1\in {\Cal S}^{\text{na}}_{{\frak K}_\lambda}(M_1)$
are there $M_2 \in {\frak K}_\lambda$ and a minimal $p_2
\in {\Cal S}^{\text{na}}_{{\frak K}_\lambda}(M_2)$ such that
$M_1 \le_{{\frak K}_\lambda} M_2$ and $p_2$ extends $p_1$?

As we are assuming categoricity in $\lambda$ and $\lambda^+$, this is
not unreasonable and its failure implies having large 
${\Cal S}^{\text{na}}_{{\frak K}_\lambda}(M)$.  Now
\cite[\S3,\S4]{Sh:E46} relying on \cite{Sh:838} (earlier:
\cite{Sh:603} and part of \cite{Sh:576}) are 
dedicated to proving that the minimal 1-types are dense. 
(This requires looking more into the set theoretic side but also the
model theoretic one; an example of a property which we consider is: 
given  $M_0 <_{{\frak K}_\lambda} M_1$ is there 
$M_2,M_0 <_{{\frak K}_\lambda} M_2$ such that $M_1,M_2$ can be
amalgamated over $M_0$ uniquely?).
\bn
So we assume the answer to \scite{E53-qu.7.1} is yes; that is, we 
make the hypothesis:
\demo{\stag{E53-ans.7.1} Hypothesis}   The answer to question
\scite{E53-qu.7.1} is yes: the minimal types are dense.

Having arrived here, further investigation shows
\mr
\item "{$(*)$}"  ${\Cal S}^{\text{na}}_{{\frak K}_\lambda}(M)$ has
cardinality $\le \lambda$.
\ermn
Now it is natural to define $(M_0,M_1,a,M_3)$ is a non-forking
quadruple or $\nonfork{}{}_{\frak s}(M_0,M_1,a,M_3)$
iff $M_0 \le_{{\frak K}_\lambda} M_1 \le_{{\frak K}_\lambda} M_3,a \in
M_3 \backslash M_1$ and $\ortp_{{\frak K}_\lambda}(a,M_0,M_3)$ is
minimal.  Recalling Candide we note that 
having chosen the unique non-trivial extension, we
certainly have made the free choice:  we have no freedom left on what is 
$\ortp_{{\frak K}_\lambda}(a,M_1,M_3)$!  Now we find a good $\lambda$-frame
${\frak s}$, with ${\frak K}_{\frak s} = {\frak K}_\lambda$ and
${\frak K}^{\frak s} = {\frak K}[{\frak s}]$ will denote ${\frak K}_{\ge
\lambda} = {\frak K} \restriction \{M \in K:\|M\| \ge
\lambda\}$ and the set of basic types, is 
${\Cal S}^{\text{bs}}_{{\frak K}_\lambda}(M)$, 
the set of minimal $p \in {\Cal S}^{\text{na}}_{{\frak K}_\lambda}(M)$.  
Note that good $\lambda$-frame is defined in 
\sectioncite[\S2]{600}, existence in our case is proved in \sectioncite[\S3]{600}.

\beginaside
More accurately, in \sectioncite[\S3]{600} we prove in our present context the
existence of a good $\lambda^+$-frame ${\frak s}$ with ${\frak
K}_{\frak s} = {\frak K}_{\lambda^+}$, and we rely on having developed
NF$_\lambda$ in \cite[\S8]{Sh:576}.  But something parallel to
\cite[\S8]{Sh:576} is done in \sectioncite[\S6]{600} and described below.
Moreover, in \cite{Sh:E46} this is circumvented at the price of arriving
to almost good $\lambda$-frame and then by \cite{Sh:838} it is even a
good $\lambda$-frame and it converges with the description here.
\endaside

We assume here that ${\frak K}_{\frak s} (= {\frak K}^{\frak
s}_\lambda$) is categorical; in the present
context this is reasonable (e.g., as otherwise you
restrict yourself to $\{M \in {\frak K}_{\frak s}:M$ is superlimit$\}$).
\bn
(B) The successor of a good $\lambda$-frame
\sn
Now we look at our good $\lambda$-frame ${\frak s}$, and the 
${\frak s}$-basic types in this case 
are the minimal types.  But we can forget the minimality and
just use the properties required in the definition of a good
$\lambda$-frame (i.e. we are in \chaptercite{600}).  
Now as $M \in {\frak K}_{\frak s} \Rightarrow {\Cal
S}^{\text{bs}}_{\frak s}(M)$ has cardinality $\le \lambda$ (by the
definition of a good $\lambda$-frame) we can find
$\le_{\frak s}$-increasing chains $\langle M_i:i \le \lambda \times
\delta \rangle$ such that for every $i < \lambda \times \delta$ every $p
\in {\Cal S}^{\text{bs}}_{\frak s}(M_i)$ is realized in $M_{i+1}$.  
In such a case we say that $M_{\lambda \times \delta}$ is brimmed over $M_0$.
It follows that $M_{\lambda \times \delta}$ is determined uniquely up to
isomorphisms over $M_0$ (seemingly, depending on cf$(\delta) := \text{ Min}
\{\text{otp}(C):C \subseteq \delta$ unbounded$\}$).  
Eventually we succeed to prove that
the choice of the limit ordinal $\delta ( < \lambda^+)$ is immaterial,
see \sectioncite[\S1,\S4]{600}.  
\beginaside
(These are relatives of universal homogeneous, saturated models and
special models.)
\endaside
\sn
We define $K^{3,\text{bs}}_{\frak s}$ as the class of triples $(M,N,a)$
such that $M \le_{{\frak K}_{\frak s}} N$ and $\ortp_{{\frak K}_{\frak
s}}(a,M,N) \in {\Cal S}^{\text{bs}}_{\frak s}(M)$.
By the axioms of ``good $\lambda$-frames"
for $(M_1,N_1,a) \in K^{3,\text{bs}}_{\frak s}$ and
$M_2$ such that $M_1 \le_{\frak s} M_2$ we can find $M'_2 \in {\frak
K}_\lambda$ isomorphic to $M_2$ over $M_1$ and $N_2 \in {\frak
K}_\lambda$, which is $\le_{\frak K}$-above $M'_2$ and $N_1$ and 
$\ortp_{\frak s}(a,M'_2,N_2)$
does not fork over $M_1$.  In this case we say $(M_1,N_1,a) 
\le_{\frak s} (M'_2,N_2,a)$, (or use $\le_{\text{bs}} = \le^{\frak
s}_{\text{bs}}$ instead $\le_{\frak s}$).  

Having existence is nice, but having also 
uniqueness is better.  So we become interested in
$K^{3,\text{uq}}_{\frak s}$, the class of $(M,N,a) \in
K^{3,\text{bs}}_{\frak s}$ satisfying: if
$(M_*,N_*,a) \in K^{3,\text{bs}}_{\frak s}$ is $\le_{\frak s}$-above
$(M,N,a)$, then the way $M_*,N$ are amalgamated over $M$ inside 
$N_*$ is \ub{unique} (up to common embeddings).
\sn
\beginaside
For the first order case this means ``$\sftp(N,M \cup \{a\})$ is weakly 
orthogonal to $M$"; (i.e., domination).
\endaside
\enddemo
\bn
\margintag{E53-qu.7.2}\ub{\stag{E53-qu.7.2} Question}:  1) (Density) Do we have
``$K^{3,\text{uq}}_{\frak s}$ is dense in 
$K^{3,\text{bs}}_{\frak s}$ (under $\le_{\frak s}$)"? 
\nl
2) (Existence)  Assume $p \in {\Cal S}^{\text{bs}}_{\frak
s}(M)$, can we find $a,N$ such that $(M,N,a) \in
K^{3,\text{uq}}_{\frak s}$ and $\ortp_{\frak s}(M,N,a) = p$?  

As ${\frak K}_{\frak s}$ is categorical, we can prove that 
density implies existence.

``Have we not been here before?" the reader may wonder.  This is the spiral
phenomena: in \scite{E53-qu.7.1} we were interested in a different kind of 
uniqueness.  Now we prove that the non-density is a non-structure 
property and as a token of our pleasure, ${\frak s}$ with positive
answer is called weakly successful.
\bigskip

\demo{\stag{E53-ans.7.2} Hypothesis}  The answer to \scite{E53-qu.7.2} is yes, 
enough triples in $K^{3,\text{uq}}_{\frak s}$ exist.

So we have some cases of uniqueness of the non-forking amalgamation.
When we (in \sectioncite[\S6]{600}) close this family of cases of uniqueness,
under transitivity and monotonicity we get a four-place relation
NF$_\lambda = \text{ NF}_{\frak s}$ on ${\frak K}_\lambda$.  Working
enough we show that NF$_{\frak s}$ conforms reasonably 
with ``$M_1,M_2$ and are in non-forking ($\equiv$
free) amalgamation over $M_0$ inside $M_3$".  We justify the definition
showing that some natural properties it satisfies has at most one
solution (for any good $\lambda$-frame).

Now we start to look at models in $K^{\frak s}_{\lambda^+}$; in an
attempt to find a good $\lambda^+$-frame ${\frak s}^+ = s(+)$,
a successor of ${\frak s}$.  There are some models in $K^{\frak
s}_{\lambda^+}$; in fact, there is a universal homogeneous one $M^*$
and it is unique so if there is a superlimit $M \in K^{\frak
s}_{\lambda^+}$ then $M \cong M^*$.  Now if $\langle M_i:i < \lambda^+
\rangle$ is $\le_{{\frak K}^{\frak s}_{\lambda^+}}$-increasing $M_i
\cong M^*$ then $\cup\{M_i:i < \lambda^+\} \cong M^*$ but it is not
clear if, e.g., $\cup\{M_i:i < \omega\} \cong M^*$.  So we consider
another choice of being a substructure in $K^{\frak s}_{\lambda^+}:M_1
\le^*_{\lambda^+} M_2$ iff $M_1,M_2 \cong M^*$ and for some
$\le_{\frak K}$-representations (also called $\le_{\frak K}$-filtrations)
$\langle M^\ell_\alpha:\alpha < \lambda^+ \rangle$ of
$M_\ell$ for $\ell=1,2$ we have NF$_{\frak s}
(M^1_i,M^2_i,M^1_j,M^2_j)$ for every $i < j < \lambda^+$. 
\sn
\beginaside
[We say that $\langle M_\alpha:\alpha < \lambda^+\rangle$ is a
$\le_{\frak K}$-representation or
$\le_{\frak K}$-filtration of $M \in {\frak K}_{\lambda^+}$ when
$M_\alpha \in {\frak K}_\lambda$ is $\le_{{\frak
K}_\lambda}$-increasing continuous for $\alpha < \lambda^+$ and $M =
\cup\{M_\alpha:\alpha < \lambda^+\}$.]
\endaside
\sn
We would love to understand ${\frak K}_{\lambda^+}$, but this seems too
hard, so presently so we restrict ourselves to isomorphic copies of the
model we do understand, $M^*$.
\sn
\beginaside
This conforms with the strategy of first
understanding the quite saturated models.
\endaside

This helps to prove
``$M^*$ is superlimit" but with a price: we have to consider the
following question.
\enddemo
\bn
\margintag{E53-qu.7.3}\ub{\stag{E53-qu.7.3} Question}:  Assume $\langle M_i:i \le \delta \rangle$
is $\le^*_{\lambda^+}$-increasing continuous, $\delta$ a limit ordinal
$< \lambda^{++}$ and $i < \delta \Rightarrow M_i \cong M^*$ and $i <
\delta \Rightarrow M_i \le^*_{\lambda^+} N$ and $N \cong M^*$.  Does
it follow that $M_\delta \le^*_{\lambda^+} N$? 

This is an axiom of an abstract elementary class, so 
we know that it holds for $({\frak K}_{\lambda^+},
\le_{\frak K})$ but not necessarily for $\le^*_{\lambda^+}$.  
This is another
dividing line: if the answer is no, we get a non-structure theorem.
If the answer is yes, we call ${\frak s}$ successful.
\bigskip

\demo{\stag{E53-ans.7.3} Hypothesis}  ${\frak s}$ is successful.

We go on and prove that ${\frak s}^+$ is a good $\lambda^+$-frame.
Well, the reader may wonder: all this work and you just end up where
you have started, just one cardinal up?  True, but if ${\frak s}$ is a good
$\lambda$-frame then $K^{\frak s}_{\lambda^{++}} \ne \emptyset$, so for a
successful ${\frak s}$, applying this to the good $\lambda^+$-frame
${\frak s}^+$ we get ${\frak K}_{\lambda^{+3}_{\frak s}} \ne
\emptyset$.  Having ``arrived to the same place one cardinal up" is
enough to prove part (1) of Theorem \scite{E53-nl.2.7}!  

More elaborately, under the assumptions of \scite{E53-nl.2.7} there is a
good $\lambda^+$-frame ${\frak s}_1$ with ${\frak K}^{{\frak s}_1}
\subseteq {\frak K}^{\frak s}$.  Second, if we prove by induction on
$k=1,\dotsc,n-1$ that there is a good $\lambda^{+k}$-frame ${\frak
s}_k$ with $K^{{\frak s}_k} \subseteq K^{{\frak s}_{k-1}}$, the 
induction step is what
we have proved.  For $k=n-1, ``K^{{\frak s}_k}$ has a model in
$\lambda^{++}_{{\frak s}_k}"$ means that $K_{\lambda^{+n+1}} \ne
\emptyset$ as asked for in \scite{E53-nl.2.7}(1).
All this is \chaptercite{600},
so its proof proceeds by ``forgetting" the previous ${\frak s}$ when
advancing ${\frak s}^+$ and $\lambda^+_{\frak s}$.  Next assume
\mr
\item "{$\boxdot_2$}"  ${\frak s}$ is a $\lambda$-good frame, $\dot
I(\lambda^{+n},{\frak K}^{\frak s}) < 2^{\lambda^{+n}}$ and
$2^{\lambda^{+n}} < 2^{\lambda^{+n+1}}$ for $n <
\omega$.
\ermn
We now define by induction on $n$ a good $\lambda^{+n}$-frame ${\frak
s}^{+n} = {\frak s}(+n)$.  Let ${\frak s}^0 = {\frak s}$ and having
defined ${\frak s}^{+n}$, it has to be successful by the previous
argument so ${\frak s}^{+(n+1)} := ({\frak s}^{+n})^+$ is a well defined
good $\lambda^{+(n+1)}$-frame.
We can prove by induction on $n$ that $K_{{\frak s}(+n)} \subseteq 
K^{\frak s}$ and $m < n \Rightarrow K^{{\frak s}(+n)} \subseteq 
K^{{\frak s}(+m)}$.

Note that if $K^{\frak s}$ is the class of $(A,E)$ where $|A| \ge
\lambda$ and $E$ is an equivalence relation on $A$ then 
$K^{{\frak s}^{+n}}$ is the class of $(A,E) \in K^{\frak s}$ 
such that $E$ has $\ge \lambda^{+n}$ equivalence classes each of
cardinality $\ge \lambda^{+n}$. 
\bn
(C) The beauty of $\omega$ successive good $\lambda$-frames
\sn
What about part (2) of \scite{E53-nl.2.7}, i.e., models in 
cardinalities $\ge \lambda^{+
\omega}$?  The connection between ${\frak s}^{+n},{\frak s}^{+(n+1)}$
is not strict enough.  Now though we have $K^{{\frak s}^{+n+1}} \subseteq
K^{{\frak s}^{+n}}$,  we do not know whether 
$\le_{{\frak s}(+n+1)}$ is $\le_{{\frak K}[{\frak s}(+n)]}
\restriction K_{{\frak s}(+n+1)}$ and whether ${\frak K}_{{\frak
s}(+n+1)} = K^{{\frak s}(+n)}_{\lambda^{+n+1}}$.  We can overcome the first
problem.  We show that if ${\frak s}$ is so called good$^+$ then
$\le_{{\frak s}(+)} = \le_{{\frak K}[{\frak s}]} \restriction
K_{{\frak s}(+)}$ (and ${\frak s}$ is good$^+$ ``usually" holds e.g., 
if ${\frak s} = {\frak t}^+,{\frak t}$ is good$^+$ and
successful, see \sectioncite[\S1]{705}).
In this case $\langle {\frak K}^{{\frak
s}^{+n}}:n < \omega \rangle$ is decreasing and even $\langle K^{{\frak
s}^{+m}}_{\lambda^{+n}}:m \le n \rangle$ is decreasing in $m$, but the orders
agree when well defined.  The crux of the matter is in the end
(\sectioncite[\S12]{705}, relying on what we prove earlier in \chaptercite{705}), 
to show that for some ${\frak s}^{+ \omega},K_{{\frak s}^{+ \omega}} 
= \cap\{{\frak K}^{{\frak s}(+n)}_{\lambda^{+ \omega}}:n <
\omega\}$ and ${\frak s}^{+ \omega}$ is so called beautiful, so at 
last we shall arrive to ``the promised land" from 
$\boxtimes(c)$ from the beginning of \S2(C).
But this comes
only at the very end.  In particular before starting we have to know much
on the ${\frak K}_{{\frak s}(+n)}$'s.  It is enough to prove any of 
the nice things we like to know on ${\frak K}_{{\frak
s}(+n)}$ just for ``$n < \omega$ large enough".
A priori we may have from time to time to say ``if ${\frak s}$ has the
desirable properties $(A)_1,\dotsc,(A)_{\ell-1}$ then ${\frak s}^{+n}$
has $(A)_\ell$ (as we are assuming all ${\frak s}^{+n}(n < \omega)$ are
successful), and so when we prove a desirable property $X$ we prove it for
${\frak s}^{+n}$ when $n \ge n_X$".  Originally we
were using $n \ge 2$ or $n \ge 3$, but try to use 
little, say ``${\frak s}$ is weakly successful" (which means $n$ is 0
or 1) and lately try just to finish.

Note also that \wilog \, ${\frak s}$ is type-full, i.e. 
${\Cal S}^{\text{bs}}_{\frak s}(M) = {\Cal S}^{\text{na}}_{\frak s}(M)$, as we
can use our knowledge on NF$_{\frak s}$ to define when ``$p \in
{\Cal S}^{\text{na}}_{\frak s}(N)$ does not fork over $M \le_{\frak s}
N$" and prove that ${\frak t}$ is a good $\lambda$-frame when we
define ${\frak t}$ by ${\frak K}_{\frak t} = {\frak K}_{\frak s},{\Cal
S}^{\text{bs}}_{\frak t} = {\Cal S}^{\text{na}}_{\frak s}$, and non-forking
as above.  As we can replace ${\frak s}$ by ${\frak t}$ the
``w.l.o.g." above is justified.

Note that the ${\frak K}_{{\frak s}(+n)}$ are categorical, but this is
deceptive: ${\frak K}_{{\frak s}(+n)}$ is, but $K^{{\frak
s}(+n)}_{\lambda^{+n+1}}$ is not necessarily categorical.  So in order
to eventually understand the categoricity spectrum in
\sectioncite[\S2]{705} we sort
out when is ${\frak K}^{\frak s}_{\lambda^+}$ categorical (for a
successful good $\lambda$-frame ${\frak s}$). 

We define several (variants of) ${\frak s}$ is uni-dimensional, prove
the equivalence with $``K^{\frak s}$ is categorical in $\lambda^+_{\frak
s}"$ and show that (for successful ${\frak s}$) ${\frak s}$ is
uni-dimensional iff ${\frak s}^+$ is uni-dimensional (so this applies to
${\frak s}^{+n}$ and ${\frak s}^{+(n+1)}$ when well defined).  So in the
case we have chosen, 
${\frak s}^+,{\frak s}^{+2},...$ are uni-dimensional and 
$K^{{\frak s}(+n)}_{\lambda^{+n}} = K^{\frak s}_{\lambda^{+n}}$ 
so in the beautiful (see below) 
case it implies categoricity in all $\mu > \lambda$.

We now review \chaptercite{705} in more detail.  We define and 
investigate ``$\bold J$ is a set of elements in $N \backslash
M$ which is independent over $M$" in symbols 
$(M,N,\bold J) \in K^{3,\text{bs}}_{\frak s}$.  
The idea is that if $\langle M_i:i
\le \alpha \rangle$ is $\le_{\frak s}$-increasing, $a_i \in M_{i+1}
\backslash M_i$ and $\ortp_{\frak s}(a_i,M_i,M_{i+1})$ does not fork over
$M_i$ for $i < \alpha$, then $(M_0,M_\alpha,\{a_i:i < \alpha\}) \in
K^{3,\text{bs}}_{\frak s}$ and even $(M_0,M',\{a_i:i < \alpha\}) \in
K^{3,\text{bs}}_{\frak s}$ if $M \cup\{a_i:i < \alpha\} \subseteq M'
\le_{\frak s} M_\alpha$.
\ub{But} we have to prove that this notion has the expected
properties, e.g., the finite character (see \sectioncite[\S5]{705}).

We know about $(M,N,a) \in K^{3,\text{uq}}_{\frak s}$, but also
important is $(M,N,a) \in K^{3,\text{pr}}_{\frak s}$: the triple is prime,
i.e., such that if $(M,N',a') \in K^{3,\text{bs}}_{\frak s}$ and
$\ortp_{\frak s}(a,M,N) = \ortp_{\frak s}(a',M,N')$ then there is a
$\le_{\frak s}$-embedding of $N$ into $N'$ over $M$ mapping $a$ to
$a'$.  We prove existence in enough cases (mainly for ${\frak s}^+$) 
and eventually define and investigate also ``$N$ is prime over $M
\cup \bold J$" when $(M,N,\bold J) \in K^{3,\text{bs}}_{\frak s}$ and
$\bold J$ is maximal.

Next we develop orthogonality: assume $p_\ell \in {\Cal
S}^{\text{bs}}_{\frak s}(M)$ for $\ell=1,2$. Then $p_1 \perp p_2$ when:
if $(M,N,a) \in K^{3,\text{uq}}_{\frak s}$ and $p_1 = 
\ortp_{\frak s}(a,M,N)$ then $p_2$ has a
unique extension in ${\Cal S}_{\frak s}(N)$.  This means that there is
no connection, no interaction between $p_1$ and $p_2$.  It implies that
$(M,N,\{a_i:i < \alpha\}) \in K^{3,\text{bs}}_{\frak s}$, i.e.,
is independent iff for each $j < \alpha,(M,N,\{a_i:i < \alpha,p_j \perp
p_i\})$ is independent where $p_i = \ortp_{\frak s}(a_i,M,N)$.
We prove that this behaves reasonably; in particular, is preserved by
non-forking extensions.  We similarly define $p \perp M$ (when $M
\le_{\frak s} N,p \in {\Cal S}^{\text{bs}}_{\frak s}(N))$.  Because of
the categoricity (and ${\frak s} = {\frak t}^+$)
we can prove $K^{3,\text{pr}}_{\frak s} = K^{3,\text{uq}}_{\frak s}$.

In those terms we can characterize when 
$(M,N,a) \in K^{3,\text{bs}}_{\frak s}$ has uniqueness (i.e., 
$\in K^{3,\text{uq}}_{\frak s}$), under
the assumption that there are primes.  It holds 
\ub{iff} there is a decomposition $\langle(M_i,a_j):i \le \alpha,j
< \alpha \rangle$ of $(M,N)$, i.e., $M_0 = M,M_\alpha =
N,(M_i,M_{i+1},a_i) \in K^{3,\text{pr}}_{\frak s}$ such that $a_0=a$ and
$i \in (0,\alpha) \Rightarrow \ortp_{\frak s}(a_i,M_i,M_{i+1})
\perp M_0$.  We can define regular types such that: for $M \le_{\frak
s} N$ and regular $p \in {\Cal S}^{\text{bs}}_{\frak s}(M)$ the
dependence relation on $\bold I_{M,N} = \{a \in N:a$ realizes $p\}$
behaves as independence in vector spaces (for
others it behaves like sets of finite sequences from a vector space),
and regular types are dense (i.e., if $M <_{\frak s} N$ then for some $a \in N
\backslash M,\ortp_{\frak s}(a,M,N)$ is regular).
So $a \in \bold I_{M,N}$ depends on $\bold J \subseteq \bold I_{M,N}$
iff there are
$M_1 \le_{{\frak K}_{\frak s}} N_1$ such that $M \le_{{\frak K}_{\frak s}} 
M_1,N \le_{{\frak K}_{\frak s}} N_1,\bold J \subseteq M_1$, the triple
$(M,N,\bold J)$ has uniqueness and
$\ortp_{\frak s}(a,M_1,N_1)$ forks over $M$.  It has local
character (if $a \in \bold I_{M,N}$ depends on $\bold J$ then it depends on
some finite subsets of it), monotonicity, transitivity (if $a \in \bold
I_{M,N}$ depends on $\bold J' \subseteq \bold I_{M,N}$ 
and each $b \in \bold J'$ depends on 
$\bold J \subseteq \bold I_{M,N}$ then $a$ depends on $\bold J$)
and satisfies the exchange lemma.
Then we can define (and prove the relevant properties) 
when ``$\{M_i:i < \alpha\}$ is independent over $M$ inside
$N$" and we can deal similarly with 
``$\langle M_\eta:\eta \in {\Cal T}\rangle$ is independent
inside $N$" when ${\Cal T} \subseteq {}^{\omega >}(\lambda_{\frak s})$ is
closed under initial segments.

We may now consider the main gap in this context (but mostly this is
delayed).  From some
perspective this is ridiculous: ${\frak K}_{\frak s}$ is categorical
in $\lambda_{\frak s}$.  But we analyze $\{N:M_* \le_{\frak s} N\}$
for a fixed $M_*$, so all the models in this class have cardinality
$\lambda_s$.   (In this still there is some degeneration, but we can
analyze models from ${\frak K}^{\frak s}_{\lambda^+}$, in this case there
is no real difference between what we do and the actual main gap
theorem, so again all models have a fixed cardinality.  
And if ${\frak s}$ is beautiful, see below, we can do the same for
${\frak K}^{\frak s}$.  If ${\frak s}$ is ``good enough up to
$\lambda^{+n}$ we can deal similarly with $K^{\frak s}_{\le \lambda^{+n}})$.

So if $M \le_{\frak s} N$ (assuming, e.g. ${\frak s}$ is a successful
$\lambda$-frame with primes; less is needed), we can find a decomposition
$\langle N_\eta,a_\nu:\eta \in {\Cal T},\nu \in {\Cal T} \backslash
\{<>\}\rangle$ of $N$ which means
\mr
\item "{$\circledast$}"  $(a) \quad {\Cal T} \subseteq
{}^{\omega>}(\lambda_{\frak s})$ is non-empty closed under initial
segments
\sn
\item "{${{}}$}"  $(b) \quad N_\eta \le_{\frak s} N$
\sn
\item "{${{}}$}"  $(c) \quad \nu \triangleleft \eta \Rightarrow N_\nu
\le_{\frak s} N_\eta$
\sn
\item "{${{}}$}"  $(d) \quad (N_\eta,N_{\eta \char 94
<\alpha>},a_{\eta \char 94 <\alpha>}) \in K^{3,\text{pr}}_{\frak s}$
if $\eta \char 94 <\alpha> \in {\Cal T}$
\sn
\item "{${{}}$}"  $(e) \quad \{a_{\eta \char 94 <\alpha>}:\eta
\char 94 <\alpha> \in {\Cal T}\}$ is independent in $(M_\eta,N)$ and is
\nl

$\qquad$ a maximal such set (with no repetitions, of course)
\sn
\item "{${{}}$}"  $(f) \quad N_{<>} = M$,
\sn
\item "{${{}}$}"  $(g) \quad$ if $\cup\{N_\eta:\eta \in {\Cal T}\}
\subseteq N' <_{\frak s} N,p = \ortp_{\frak s}(a,N',N) \in 
{\Cal S}^{\text{bs}}_{\frak s}(N')$ \nl

$\qquad$ then $p \pm N_\eta$ for some $\eta \in {\Cal T}$
\sn
\item "{${{}}$}"  $(h) \quad$ if $\nu \triangleleft \eta \triangleleft
\eta \char 94 \langle \alpha \rangle \in {\Cal T}$ then $\ortp_{\frak
s}(a_{\eta \char 94 \langle \alpha\rangle},N_\eta,N_{\eta \char 94
\langle \alpha\rangle}) \bot N_\nu$.
\endroster
\enddemo
\bn
\margintag{E53-qu.7.4}\ub{\stag{E53-qu.7.4} Question}:  Is always $N$ prime and/or minimal over
$\cup\{N_\eta:\eta \in {\Cal T}\}$?

The answer is yes iff whenever ${\Cal T} = \{<>,<0>,<1>\}$ the answer is
yes and we then say that ${\frak s}$ have the so-called NDOP.  
Moreover, its negation DOP is a strong
non-structure property: for every $R \subseteq \lambda \times \lambda$
we can find $N_R \in K^{\frak s}_{\lambda^{++}}$ and $\bar a_\alpha,\bar
b_\alpha \in {}^{\lambda_{\frak s}}(N_R)$ for $\alpha < \lambda$ such
that some condition (preserved by isomorphism) is satisfied by $\bar
a_\alpha \char 94 \bar b_\beta$ in $N_R$ iff
$(\alpha,\beta) \in R$.  Also the NDOP holds for ${\frak s}^+$ iff
it holds for ${\frak s}$ when 
${\frak s}$ is successful from DOP.  We can get $\dot I
(\lambda^{++}_{\frak s},K^{\frak s}) = 2^{\lambda^{++}_{\frak s}}$ and
more.
\bn
\centerline {$* \qquad * \qquad *$}
\bigskip

How does all this help us to go up?  That is, we assume ${\frak
s}^{+n}$ is well defined and successful for every $n$ (equivalently
${\frak s}$ is $n$-successful for every $n$) and we would like to 
understand the models in
${\frak K}^{{\frak s}(+ \omega)}$, (so they have cardinality $\ge
\lambda^{+ \omega}$ and are close to being
$\lambda^{+\omega}$-saturated).  The going up is done in the framework
of stable ${\Cal P}^{(-)}(n)$-system of models $\langle M_u:u \in
{\Cal P}^-(n)\rangle,{\Cal P}^-(n) = \{u:u \subset\{0,\dotsc,n-1\}$;
 explained below.  This is done in \sectioncite[\S12]{705} (which should be
helpful for completing \cite{Sh:322}).

In short, to understand existence/uniqueness of models (and of
amalgamation) in $\lambda$, we consider such properties for some
$n$-dimensional systems of models in every large enough 
$\mu \le \lambda$.  So for
$n=0,1,2$ we get the original problems but understanding the $n$-th
case given in $\lambda$ is intimately connected to understand the
$(n+1)$-case for every large enough $\mu < \lambda$.  So
for $\lambda = \mu^+$ we get a positive property for $(\mu^+,n)$ from
one for $(\mu,n+1)$.

Why do we need such systems?  Consider $\lambda_* \ge \mu_* \ge
\lambda_{\frak s}$ and we try to analyze models of cardinality
$\in [\mu_*,\lambda_*]$ by pieces of cardinality 
$\mu_*$ or $\mu' \in [\mu_*,\lambda_*)$
(in the end we consider
$\mu_* = \lambda^{+ \omega}_{\frak s}$, but most of the analysis is
for the case $\lambda_*,\mu_* \in [\lambda_{\frak s},
\lambda^{+ \omega}_{\frak s}))$.
We can analyze a model $M$ from ${\frak K}$ of cardinality 
$\lambda_0 \in (\mu_*,\lambda_*]$ by a $\le_{\frak K}$-increasing 
continuous sequence $\langle M_\alpha:\alpha <
\lambda_0 \rangle,\mu_* \le \|M_\alpha\| = 
\|M_{\alpha +1}\| < \lambda_0$, with $M =
\cup\{M_\alpha:\alpha < \lambda_0\}$; so it suffices to analyze
$M_{\alpha+1}$ over $M_\alpha$ for each $\alpha$.  We can analyze $M_1$ over
$M_0$ for a pair of models $M_0 \le_{\frak K} M_1$ of the same
cardinality which we call $\lambda_1$ when $\lambda_1 > \mu_*$
by an $(\le_{\frak K})$-increasing continuous sequence of pairs
$\langle (M^0_i,M^1_i):i < \lambda_1 \rangle$ where $\|M^0_i\| =
\|M^0_{i+1}\| = \|M^1_i\| = \|M^1_{i+1}\| < \lambda_1$, and we have to
analyze $M^1_{i+1}$ over $\langle M^0_i,M^1_i,M^0_{i+1})$ for each
$i$.  In the next stage we have $8=2^3$ models and have to analyze the
largest over the rest.  Eventually we arrive to the case that all of
them have cardinality $\mu_*$.

In short, we have to consider suitable ${\Cal P}(n)$-systems $\langle M_u:u \in
{\Cal P}(n)\rangle$ where ${\Cal P}(n) = \{u:u \subseteq
\{0,\dotsc,n-1\}\},
u \subseteq v \Rightarrow M_u \le_{{\frak K}^{\frak s}} M_v$ and
$\|M_u\| = \|M_0\| \in [\mu_*,\lambda_*]$.

We would like to analyze $M_{\{0,\dotsc,n-1\}}$ over
$\cup\{M_u:u \in {\Cal P}^-(n)\}$ where ${\Cal P}^-(n) = {\Cal P}(n)
\backslash \{0,\dotsc,n-1\}$.  Such analysis of a ``big" system
of small models naturally help proving cases of uniqueness, e.g.,
uniqueness of non-forking-amalgamations suitably defined.  So if for
$\mu_*$ we have positive answers for every $n$, then this holds for
every $\lambda \in [\mu_*,\lambda_*]$.

But we are interested as well in existence proofs.  (Note that in the
proof we have to deal with uniqueness, existence (and some relatives)
simultaneously.)  
For the existence we need for a given suitable system $\langle
M_u:u \in {\Cal P}^-(n) \rangle$ to complete it by finding
$M_{\{0,\dotsc,n-1\}}$.  Well, but what are the suitable systems?
Those are defined, by several demands including 
$u \subseteq v \Rightarrow M_u \le_{{\frak K}^{\frak s}}
M_v$ (and many more restrictions which hold if the sequence of
approximations chosen above are ``fast" enough).  We called them the
stable ones.   For each $n,k$ we can ask on ${\frak
s}^{+n}$ some questions on ${\Cal P}(k)$-systems: mainly versions of
existence and
uniqueness.  A  major point is that failure of uniqueness for
$\lambda^{+n},{\Cal P}(m+1)$ implies failure for
$\lambda^{+n+1},{\Cal P}(m)$ (using $2^{\lambda^{+n}} <
2^{\lambda^{+n+1}}$). But to get strong dichotomy we have to use
systems which have the right amount of brimmness.  
At last we have a glimpse of ``paradise", we can define when ${\frak
s}$ is $n$-beautiful essentially when it satisfies all the good
properties on stable ${\Cal P}(m)$-systems for $m \le n$.  In the end we prove
that ${\frak s}^{+n}$ is $(n+2)$-beautiful, i.e. has all the desired
properties for $m \le n+2$ but for this we use ${\frak s}^{n+ \ell}$ being
successful for $\ell \le n$.

Having all this we can prove that ${\frak s}^{+ \omega}$ has all the
good properties (but we have to work on changing the brimmness
demands) so is $\omega$-beautiful.  This now can be lifted up, in
particular ${\frak K}^{{\frak s}(+ \omega)}$ has amalgamation and the
types $\ortp_{{\frak K}[{\frak s}(+ \omega)]}(a,M,N)$ are $\mu$-local for
$\mu = \lambda^{+ \omega}$ (in fact $\mu = \lambda$ is enough) where 
\mr
\item "{$(*)$}"  ${\frak K}$  an \aec \, with amalgamation, is
$\mu$-local when for $M \le_{\frak K} N$ and $a_1,a_2 \in N$ we have:
\nl
$\ortp_{\frak K}(a_1,M,N) = \ortp_{\frak K}
(a_2,M,N)$ \ub{iff} for every $M' \le_{\frak K} M$ of cardinality
$\mu,\ortp_{\frak K}(a_2,M',N) = \ortp_{\frak K}(a_2,M',N)$.
\ermn
Now for a beautiful ${\frak s}$, in particular we have
amalgamation/stable amalgamation, prime models over a triple of models
in stable amalgamation.  In particular we can prove the main gap.
However, here we just present the characterization of the categoricity
spectrum (see \scite{E53-nl.2.7}(2)) and delay the rest.
\nl
On \chaptercite{734} and \cite{Sh:838} see \S4(B).

We may wonder how excellent classes and beautiful good
$\lambda$-frames are related.  We explain this by comparing each to
the first order cases.

If $T$ is a (first order complete) theory $T$ which is
$\aleph_0$-stable then $(\text{Mod}_T,\prec)$ is an excellent class.
Pedantically assume $T$ categorical in $\aleph_0$. Better expand each
$M \in \text{ Mod}_T$ to $M^+ = (M,P^{M^+}_P)_{P \in D(T)}$ with
$P^{M^+}_{p(x)} = \{\bar a \in {}^{\ell g(\bar a)}M:\text{tp}(\bar
a,\emptyset,M) = p(x)\}$; as written in \cite{Sh:82} categoriticity in
$\aleph_0$, is assume but this is not used in any real way.

For an excellent class for simplicity is, for notational transparency,
the class of almac models of a first order $T$.  Now we continue to
use ``classical" types and (respecting atomicity) have primes even
primary models but only over sets like $\cup\{M_u:u \in p/n\}$ inside
$M$ where $\bar M = \langle M_u:u \in {\Cal P}(n)\rangle$ is a stable
system of models.  Stable non-forking is defined only for such $\bar
M$'s but still using formulas and (classical) types.  Big differences,
but usually any concept that is not obviously irrelevant can be
developed in this context.

If $T$ is a superstable first order $T$ stable in $\lambda$ \ub{then}
there is a beautiful $\lambda$-frame ${\frak s} = {\frak
s}_{T,\aleph_0},K^{\frak s}$ is a class of
$\aleph_\varepsilon$-saturated models of $T$ of cardinality $\ge
\lambda,\le_{\frak s}$ is $\prec$ on this class.  But for general
beautiful types are defined without formulas - by orbits (unlike
excellent classes).  As in excellent classes stable or non-forking
systems $\langle M_u:u \in {\Cal P} \subseteq {\Cal P}(u)\rangle$ are
central, but their definition is not direct, certainly not referring
to classical types.  Also ``$M_1$ is prime over $M^+$" is again
central, but defines by arrow (no ``primary model".  

Again usually whatever we can prove for ${\frak s}$ of the form
$\varepsilon - {\frak s}_{T,\aleph_\varepsilon}$, see above, and is
not obviously irrelevant can be proved in this content.
\newpage

\head {\S4 Appetite comes with eating} \endhead  \resetall \sectno=4
 \spuriousreset
\bigskip

Here we mainly review open questions, \chaptercite{734}, \cite{Sh:838} and
further relevant works which could have been
part of this book but were not completely ready; so decided not to wait
because my record of dragging almost finished books is 
bad enough even without this case.  Note that \chaptercite{734} use
infinitary logics and most of \cite{Sh:838} has largely set theoretic
character hence does not fit with \S2,\S3.

But we begin by looking at what has been described so far 
has not accomplished.  (By this division we end up dealing 
with some issues more than once.)
\bn
(A) \ub{The empty half of the glass}:
\sn
$(a)$  \ub{Categoricity in one large enough $\lambda$}:

We have here concentrated on going up in cardinality, (assuming that 
 in $\omega$ successive cardinals there are not too many models
 without even assuming the
existence of models of cardinality $\ge \lambda^{+3}!$).  We use 
weak instances of GCH \, $(2^\lambda < 2^{\lambda^+})$ and prove a
generalization of \cite{Sh:87a}, \cite{Sh:87b}.  But originally, and
it still seems a priori more reasonable,  probably even more
 central case should be to start assuming categoricity
in some high enough cardinal.  There are several 
approximations in Makkai-Shelah \cite{MaSh:285}, Kolman-Shelah
\cite{KlSh:362}, \cite{Sh:472} using so called ``large cardinals".
\sn
\beginaside
(Compact cardinals in the first, measurable cardinal in the second and
third).
\endaside
\mn
$(b)$ \ub{Main Gap}:

If we assume that for some ``large enough" $\lambda$, we do not
have ``many very complicated models", we expect to be able to show the
class is ``managable", hence has a structure theory.  But the proofs
described above, do not do that job.  Not only do we usually start with
categoricity assumptions, in our main line here we learn whatever we learn
only on the $\lambda^{+ \omega}$-brimmed models.  However, just on the
class of models, i.e., on the original ${\frak K}$, we know little.  This is
not surprising as, e.g. for elementary classes with
countable vocabulary, the solution of \L os
conjecture predates the main gap considerably.
\mn
$(c)$ \ub{Superstability}:

Having claimed that the superstability is a central dividing
lines, it is unsatisfactory to arrive at it here from categoricity assumptions
only.

That is, the detailed building of apparatus parallel to superstability
is here applied to cases in which we start assuming suitable
categoricity assumption, prove there are relevant good
$\lambda$-frames and continue.  (But if $\psi \in
\Bbb L_{\omega_1,\omega}$ or ${\frak K}$ is an abstract elementary
class which is
PC$_{\aleph_0}$ and $2^{\aleph_0} + \dot I(\aleph_1,{\frak K}) <
2^{\aleph_1}$ this is not so: by \sectioncite[\S3]{600} there is a good
$\aleph_0$-frame ${\frak s}$ whose $\aleph_1$-saturated models belongs
to Mod$_\psi$ but ${\frak s}$ is not necessarily uni-dimensional
(which is the ``internal" form of categoricity)).  
Probably the main weakness of
beautiful $\lambda$-frames as a candidate to being the true
superstable is the lack of non-structure results which are not
``local".  Presently, the results are about ``failure of
categoricity", see  \sectioncite[\S12]{705} where for beautiful ${\frak
s}$, which is not uni-dimensional, we prove non-categoricity in
$K^{\frak s}_\mu$ for every $\mu >$.  So natural candidate version of
solvability, see \cite{Sh:842}.
\mn
$(d)$ \ub{$\aleph_1$-compact structures}:

We may like to relax the definition of \aec \, to investigate
classes of structures satisfying some kind of countable
compactness, i.e., any reasonable countable set of demands has a
solution.  This will include 
``$\aleph_1$-saturated models" of an elementary class (even with countable
vocabulary) also complete metric spaces but those are closer to
elementary classes.

What we lose is closure under unions of $\omega$-chains.  For
elementary classes this corresponds to $\aleph_1$-saturated models
(more generally, LST$({\frak K})^+$-saturated) and we have stable instead of
superstable (the class of complete metric spaces is closer to
elementary classes).  We have considerable knowledge on the stable
case but much less
than on superstable ones.  In particular, even for elementary classes
with countable vocabulary the main gap for stable $\aleph_1$ 
models is not known.
\bn
$(e)$ \ub{Some unaesthetic points in Theorem \scite{E53-nl.2.2}}

One of them is that from \cite{Sh:576} we get (in \sectioncite[\S3]{600}) a good
$\lambda^+$-frame and not a good $\lambda$-frame.  Second, we use here
(in III) for simplicity in the non-structure results an extra
set theoretic assumption, though a very weak one.
\sn
\beginaside
Namely, the weak diamond ideal on $\lambda^+$ is not
$\lambda^{++}$-saturated.  The negation of this statement, if
consistent, has high
consistency strength.  In fact, my attempts to derive good $\lambda$-frames
from \cite{Sh:576} or dealing with weaker versions had delayed
\chaptercite{600} considerable.
\endaside
\bn
$(f)$ \ub{Lack of Counter-examples}:

By Hart-Shelah \cite{HaSh:323}, Shelah-Villaveces \cite{ShVi:648}
there are some examples for the categoricity spectrum being
non-trivial.  Still in many theorems on dividing lines
it is not proved that they are real, i.e., that there are examples.
\bn
$(g)$ \ub{Natural Examples}:

This bothers me even less than clause (f) but for many investigators 
the major drawback
is lack of ``natural examples", i.e., finding classes which are already
important where the theory developed on the structure side throw light
on the special case.  (E.g., for
simple theories, pseudo finite fields; for $\aleph_0$-stable theories,
differentially closed fields of characteristic zero; for countable
stable theories, differentially closed fields of characteristic $p>0$
(and even separably closed fields of charactertistic $p>0$)).
\nl
But see Zilber works, especially Ravello paper.
\bn
(B) \ub{The full half and half baked}: 

Some works throw some light on some of the points from (A), in
particular \chaptercite{734}, \cite{Sh:E46}, \cite{Sh:838}.  Concerning (a), in
\chaptercite{734} we assume an \aec \, ${\frak K}$ is categorical in large
enough $\mu$
and we investigate ${\frak K}_\lambda$ for $\lambda < \mu$ which are
carefully chosen, specifically we assume
\mr
\item "{$(*)_\lambda$}"  $(a) \quad$ cf$(\lambda) = \aleph_0$ which
means $\lambda = \Sigma\{\lambda_n:n < \omega\}$ for some $\lambda_n < \lambda$
\sn
\item "{${{}}$}"  $(b) \quad \lambda = \beth_\lambda$ which means that
for every $\kappa < \lambda$ not only $2^\kappa < \lambda$ but
$\beth_\kappa < \lambda$ \nl

$\quad$ where $\beth_\alpha$ is defined inductively
by iterating exponentiation, i.e., \nl

$\quad$ defining inductively $\beth_\alpha
= \aleph_0 + \Sigma\{2^{\beth_\beta}:\beta < \alpha\}$
\ermn
\ub{or} even
\mr
\item "{$(**)_\lambda$}"  $(a) + (b) + \lambda$ is the limit of
cardinals $\lambda'$ satisfying $(*)_\lambda$.
\ermn
Are such cardinals large?  Not in the set theoretic sense (i.e., 
provably in ZFC there are such cardinals), they are in some sense
analog to the tower function in finite combinatorics.  Ignoring ``few"
exceptional $\mu$, a result of \chaptercite{734} is the existence of a 
superlimit model in ${\frak K}_\lambda$; moreover the main theorem
\marginbf{!!}{\cprefix{734}.\scite{734-f.20}} of
\chaptercite{734} says that there is a good $\lambda$-frame ${\frak s}$
 with ${\frak K}_{\frak s} \subseteq {\frak K}$; the proof uses
infinitary logics.  Also if the categoricity spectrum
contains arbitrarily large cardinals then for some closed unbounded
class $\bold C$ of cardinals, $[\lambda \in \bold C \wedge \text{\rm
cf}(\lambda) = \aleph_0 \Rightarrow {\frak K}$ categorical in $\lambda$].
It seems reasonable that this can be combined with \chaptercite{705}, but there are
difficulties.

Having \marginbf{!!}{\cprefix{734}.\scite{734-f.20}} may still leave us wondering whether
we have more tangible argument that we have advance.  So we go back to
earlier investigations of such general contexts.  Now Makkai-Shelah
\cite{MaSh:285} deal with $T \subseteq \Bbb L_{\kappa,\omega}$
categorical in some $\mu$ big enough than $\kappa + |T|$ and develop
enough theory to prove that the categoricity spectrum in an
end-segment of the cardinals starting not too far, i.e. below
${}^\beth(2^{\kappa+|T|})^+$  \ub{but}, with two
extra assumptions.

First, $\kappa$ is a strongly compact cardinal.  This is natural as our
problem is that $\Bbb L_{\kappa,\omega}$ lack many of the good
properties of first order logic, and for strongly compact cardinals,
some form of compactness is regained (even for $T \subseteq \Bbb
L_{\kappa,\kappa}$).  Still this is undesirable. 

Second, we should assume that $\mu$ is a successor cardinal, this
exhibit that the theory we build is not good enough.  Now
Kolman-Shelah \cite{KlSh:362} + \cite{Sh:472} partially rectify the
first problem: $\kappa$ is required just to be a measurable cardinal
(instead of strongly compact), still measurable is not a small cardinal.
Moreover, there is an extra, quite heavy price - we deal with the
categoricity spectrum just below $\mu$ and say nothing on it above so
the categoricity spectrum is proved to be an interval instead of an
end-segment.  A parallel work \cite{Sh:394} replace measurability by
the assumption that our ${\frak K}$ is an abstract elementary class with
amalgamation; a major point there is trying to deal with the theory
problem of locality of types (and see Baldwin \cite{Bal0x}).  
Note that in both works we get
amalgamation of ${\frak K}$ below $\mu$.

We address both cases together, assuming only that our abstract
elementary class ${\frak K}$ has the amalgamation property \ub{below}
$\mu$.  We try to eliminate those two model theoretic drawbacks:
starting from a successor cardinal, and looking only below it, in
\marginbf{!!}{\cprefix{734}.\scite{734-am3.6}}, using \chaptercite{705}.  For this we prove
that suitable cases of failure of non-structure imply cases of 
$(< \mu,\kappa)$-locality\footnote{called ``tame" by many}
 for saturated models (which
means if $p \in {\Cal S}_{\frak K}(M),M \in {\frak K}_{< \mu}$ is saturated
then $\langle p \restriction N:N \le_{\frak K} M,\|N\| =
\kappa\rangle$ determine $p$).  We also show that every $M \in K_N$ is
quite saturated, using a generalization of the stability spectrum for
linear orders from \sectioncite[\S6]{734}.

Finally, we conclude (also for abstract elementary class) ${\frak K}$ 
with amalgamations assuming enough cases of $2^\lambda <
2^{\lambda^+}$ we can characterize the categoricity spectrum
(eliminating earlier restriction to successor cardinals).  This is
done showing \chaptercite{705} applies, so we need the existence of enough
$\lambda$, such that $\langle 2^{\lambda^{+n}}:n <
\omega\rangle$ is strictly increasing.

So we have eliminated the two thorny model theoretic problems and we
eliminated the use of large cardinals but we use this weak form of
GCH, we intend to deal with it in \cite{Sh:842}.

Considering clause (b) from (A), 
the main gap, it seems far ahead.  A more basic
short-coming is that in \sectioncite[\S12]{705} we get ``${\frak s}^{+
\omega}$ is $\lambda^{+ \omega}_{\frak s}$-beautiful" and ``for
beautiful $\mu$-frame ${\frak t}$ we can prove the main gap" but
this is just for, essentially, the class of 
$\lambda^{+ \omega}_{\frak s}$-saturated models.  
\bn
Concerning (A)(c), superstability, \cite{Sh:842} suggests ``${\frak K}$
is $(\lambda,\kappa)$-solvable" as the true generalization of
superstable (remembering superstability is schizophrenic in our
context); this is weaker than categoricity and we use this assumption
in \chaptercite{734}; it is O.K. to use it always but we delay this to
\cite{Sh:842}.  Essentially it means:
\mr
\item "{$\boxdot$}"   for some vocabulary 
$\tau_1 \supseteq \tau_{\frak K}$ of cardinality
$\kappa$ and $\psi \in \Bbb L_{\kappa^+,\omega}(\tau_1),\psi$ has a
model of cardinality $\ge \beth_{(2^\kappa)^+}$ and ($[M \models
\psi \wedge \|M\| = \lambda \Rightarrow M \restriction \tau$ is
superlimit in ${\frak K}$].
\ermn
A major justification for the parallelism with superstability
is that for elementary classes this is
equivalent to superstability.  

But in \cite{Sh:842}, \sectioncite[\S12]{705} needs to be reworked 
hopefully toward the needed continuation.

We can look at results from \cite{Sh:c} which were 
not regained in beautiful $\lambda$-frames.  Well, of
course, we are far from the main gap for the original ${\frak K}$ 
(\cite[XIII]{Sh:c}) and there are
results which are obviously more strongly connected to elementary
classes, particularly ultraproducts.  This leaves us with parts of
type theory: semi-regular types, weight, ${\bold P}$-simple
\footnote{The motivation is for suitable $\bold P$ (e.g. a single
regular type)  that on the one hand stp$(a,A) \pm
\bold P \Rightarrow \text{ stp}(a/E,A)$ is $\bold P$-simple for some
equivalence relation definable over $A$ 
and on the other hand if stp$(a_i,A)$ is $\bold P$-simple for
 $i < \alpha$ then $\Sigma\{w(a_i,A) \cup \{a_j:j<i\}):i < \alpha\}$
 does not depend on the order in which we list the $a_i$'s.  Note that
 $\bold P$ here is ${\Cal P}$ there.} types, ``hereditarily orthogonal
to $\bold P$" (the last two were defined and investigated in 
\cite[V,\S0 + Def4.4-Ex4.15]{Sh:a},
\cite[V,\S0,pg.226,Def4.4-Ex4.15,pg.277-284]{Sh:c}).  The  
more general case of (strictly) stable classes was
started in \cite[V,\S5]{Sh:c} and \cite{Sh:429} and much advanced in
Hernandes \cite{He92}.

Note that ``a type $q$ is $p$-simple (or ${\bold P}$-simple)" and ``$q$
is hereditarily orthogonal to $p$ (or ${\bold P}$)" are essentially
the \footnote{Note, ``foreign to $\bold P$" and ``hereditarily
orthogonal to $\bold P$ are equivalent.  Now ($\bold P = \{p\}$ for
ease)
\mr
\item "{$(a)$}"  $q(x)$ is $p(x)$-simple when for some set $A$, in
${\frak C}$ we have $q({\frak C}) \subseteq \text{ acl}(A \cup \bigcup
p_i({\frak C}))$
\sn
\item "{$(b)$}"  $q(x)$ is $p(x)$-internal when for some set $A$, in
${\frak C}$ we have $q({\frak C}) \subseteq \text{ dcl}(A \cup
p({\frak C}))$.
\ermn
Note
\mr
\item "{$(\alpha)$}"  internal implies simple
\sn
\item "{$(\beta)$}"  if we aim at computing weights it is better to
stress acl as it covers more
\sn 
\item "{$(\gamma)$}"  but the difference is minor and
\sn
\item "{$(\delta)$}"  in existence it is better to stress dcl, also it
is useful that $\{F \restriction (p({\frak C}) \cup q({\frak C}):F$ an
automorphism of ${\frak C}$ over $p({\frak C}) \cup \text{ Dom}(p)\}$
is trivial when $q(x)$ is $p$-internal but not so for $p$-simple
(though form a pro-finite group).} 
``internal" and ``foreign" in Hrushovski's profound works.

Some years ago \cite{Sh:839} started to deal with this to some extent.  
No problem to define weight, but for having ``simple" types we need
to be somewhat more liberal in the definition of \aec \, - allow
function symbols of infinite arity (= number of places) while
preserving the uniqueness of direct limit.  In the right form which
includes the case of $\aleph_1$-saturated models of a stable theory, we
generalize what was known (for elementary classes); see more in
\scite{E53-nl.5c.4} and before.

Lastly, considering (A)(e), to a large extent this is resolved as a
product of redoing and extending the non-structure theory of
\cite{Sh:576} in \cite{Sh:838}.  

In view of \sectioncite[\S5]{88r} it is natural to
weaken the stability demand in the definition of a good $\lambda$-frame
 to $M \in K_{\frak s} \Rightarrow 
|{\Cal S}^{\text{bs}}_{\frak s}(M)| \le \lambda^+_{\frak s}$ and this
is called a semi-good $\lambda$-frame. (The present way is to choose a
countable close enough set of types and redefine $(K,\le_{\frak K}$) so
we restrict the class of models.
Semi-good frames are introduced and investigated by Jarden-Shelah
\cite{JrSh:875}.  Concerning clause (A)(f), 
Baldwin-Shelah \cite{BlSh:862} expands our knowledge of examples
considerably.  Concerning clause (A)(g) see Zilber \cite{Zi0xa},
\cite{Zi0xb}.  

In \cite{Sh:F709} may try to axiomatize the end of
\sectioncite[\S5]{88r} and connect it to good $\aleph_0$-frames, 
\cite{Sh:E54} will say more 
on Chapter II.  In \cite{Sh:838} we also deal with the positive theory
of almost good frame and weak versions of $K^{3,\text{uq}}_{\frak
s}$.  Also \cite{Sh:F735} will consider redoing \chaptercite{705}
under weaker assumptions and getting more and \cite{Sh:F782} will continue
\chaptercite{734}, e.g. how the good $\lambda$-frame from
\sectioncite[\S4]{734} fit \chaptercite{705}.  Also \cite{Sh:F888} will try to
continue \cite{Sh:E56}, and \cite{Sh:F841} to continue \cite{Sh:838}.
\bn
(C) \ub{The white part of the map}:
\sn
So we would really like to know
\nl
\margintag{E53-cat.1.A}\ub{\stag{E53-cat.1.A} Problem}:  What can be the categoricity spectrum 
Cat-Spec$_{\frak K} = \{\lambda:{\frak K}$ is categorical$\}$ for an \aec? 

This seems too hard at present
and involves independence results.  Note also that easily
(by known results, see \cite{Ke70} or see 
(\cite[VII,\S5]{Sh:c}) for any $\alpha <
\omega_1$ for some \aec \, ${\frak K}$ (with LST$({\frak K}) =
\aleph_0$) we have: $\lambda \in \text{ Cat-Spec}_{\frak K}
\Leftrightarrow \lambda > \beth_\alpha$ (just let $\psi = \psi_1 \vee
\psi_2 \in \Bbb L_{\omega_1,\omega}(\tau),\psi_1$ has a model of
cardinality $\lambda$ iff $\lambda \le \beth_\alpha$ and $\psi_2$ says
that all predicates and function symbols are trivial).

Considering the history it seemed to me that the main question on our
agenda should be
\sn
\margintag{E53-Cat.1.B}\ub{\stag{E53-Cat.1.B} Conjecture}:  If ${\frak K}$ is an \aec \, then
either every large enough $\lambda$ belongs to Cat-Spec$_{\frak K}$
\ub{or} every large enough $\lambda$ does not belong to
Cat-Spec$_{\frak K}$ (provably in ZFC).  
\bn
After (or you may say if) this is resolved positively we should consider
\demo{\stag{E53-Cat.1.7} Conjecture}  1) If ${\frak K}$ is an a.e.c. with
LST$({\frak K}) = \chi$ then
\mr
\item "{$(a)$}"  Cat-Spec$_{\frak K}$ includes or is disjoint to
$[\beth_\omega(\chi),\infty)$ 
\nl
or even better
\item "{$(a)^+$}"  similarly for $[\lambda_\omega,\infty)$ where
$\lambda_0 = \chi,\lambda_{n+1} = \text{ min}\{\lambda:2^\lambda >
2^{\lambda_n}\},\lambda_\omega = \Sigma\{\lambda_n:n < \omega\}$
\ermn
probably more realistic are
\mr
\item "{$(b)$}"  similarly for $[\beth_{(2^\chi)^+},\infty)$, or at
least
\sn
\item "{$(c)$}" similarly $(\beth_{1,1}(\chi),\infty)$ or at least
$(\beth_{1,\omega^\omega}(\chi),\infty)$, see \sectioncite[\S0]{734}.
\ermn
This will be parallel in some sense to the celebrated investigations
of the countable models for (first order) countable $T$ categorical in
$\aleph_1$. 
\enddemo
\bn
Further questions are: (recall $\boxdot$ above)
\nl
\margintag{E53-nl.5c.1}\ub{\stag{E53-nl.5c.1} Question}:  What can be 
$\{(\lambda,\kappa):{\frak K}_\lambda$ is
$(\lambda,\kappa)$-solvable, $\lambda >> \kappa >> \text{ LST}({\frak K})\}$? 
\sn
Question \scite{E53-nl.5c.1} seems to us to be more profound than the 
categoricity spectrum as solvability is a form of superstability.  We
conjecture that the situation is as in \scite{E53-Cat.1.7}(c); note that
solvability seems close to categoricity and we have a start on it
(\chaptercite{734}, \cite{Sh:842}).
\bn
Still more easily defined (but a posteri too early for us) is:
\nl
\margintag{E53-splm.7}\ub{\stag{E53-splm.7} Question}:  1) What can be $\{\lambda:{\frak K}_\lambda$
has a superlimit model$\}$?  
\nl
2) Similarly for locally superlimit (see \marginbf{!!}{\cprefix{734}.\scite{734-0n.4}}).
\nl
3) For suitable $\Phi$ what can be $\{\lambda$: if $I$ is a linear order
of cardinality $\lambda$ then EM$_{\tau({\frak K})}(I,\Phi)$ is pseudo
superlimit$\}$? see \marginbf{!!}{\cprefix{734}.\scite{734-0n.5}}(3).

We conjecture it will be a variant of
\scite{E53-Cat.1.7} but will be harder and even:
\bn
\demo{\stag{E53-splm.14} Conjecture}  If $\lambda >
\beth_{1,1}(\text{LST}_{\frak K})$ (or $\lambda >
\beth_{1,\omega}(\text{LST}_{\frak K}$), \ub{then} ${\frak K}$ has a superlimit
model in $\lambda$ iff ${\frak K}$ is $(\lambda,
\text{LST}_{\frak K})$-solvable. 
\enddemo
\bn
We now return to $(D,\lambda)$-homogeneous models.
Of course, for special $D$'s we may be interested in some special
classes of models, but not necessarily the elementary sub-models
of ${\frak C}$.
Of course, parallely to the first
order case, the main gap for them is an important problem (e.g.
the class of existentially closed models of a universal first
order theory is a natural and important case).  But  the most
natural main case seems to me the ``$\frak C$ is 
$(D,\kappa)$-sequence homogeneous" context:
\bn
\margintag{E53-nl.5c.2}\ub{\stag{E53-nl.5c.2} Problem}:  Prove the main gap for the class of
$(D,\kappa)$-sequence-homogeneous $M \prec {\frak C}$;
considering what we know, we can assume $\kappa \ge \kappa(D)$, see
\cite{Sh:3} (and \S1(B)) and
concentrate on $\kappa \ge \aleph_1$  and we would like to prove that
\mr
\item "{$(a)$}"   either 
the number of such models of cardinality $\aleph_\alpha =
\aleph_\alpha^{<\kappa(D)}+\lambda(D)$ is small, i.e.,
$\le \beth_{\gamma(D)}(|\alpha|)$ for \footnote{of course,
$\beth_{\gamma(D)}(|\alpha|)$ may be $\ge 2^{\aleph_\alpha}$ in which
case this says little; this consistently occurs for every $\alpha \ge
\omega$.  But if G.C.H. holds, and if we ask on $\dot I \dot
E(\lambda,-)$ for the class we get clear cut results} some $\gamma(D)$ 
not depending on $\alpha$ \ub{or} the number is $2^{\aleph_\alpha}$
(where $\lambda(D)$ is the first ``stability cardinal" of $D$).
\sn
\item "{$(b)$}"  $\gamma(D)$ does not depend on $\kappa$.
\ermn
\bn
A parallel of ``the main gap for the class of
$\aleph_\varepsilon$-saturated models of a first order $T$" in
this context is dealt with in Hyttinen-Shelah \cite{HySh:676}, and
a parallel to the ``main gap for the class of model of a totally
transcendental first order $T$" in Grossberg-Lessman
\cite{GrLe0x}, and surely there is more to be said in those cases
\ub{but} in the problem above, even the case 
$\kappa = \aleph_1,{\frak C}$ saturated is not covered.
\bn
We hope eventually to find a stability theory for the 
``countably compact abstract elementary class" 
strong enough to prove as a special case 
the main gap for the
$\aleph_1$-saturated models of elementary classes (i.e., clause (d) of
(A)) as said above maybe \cite{Sh:839} help.

The reader may wonder: if not known for elementary classes why you
expect more from a general frame?  Of course, we do not know, but:
\bn
\margintag{E53-nl.5b.7}\ub{\stag{E53-nl.5b.7} Thesis}:  The better closure properties 
of the abstract frames should
help us, being able to, e.g., make induction on frames.
\bn
Hence
\nl
\margintag{E53-nl.5c.4}\ub{\stag{E53-nl.5c.4} Thesis}:  Some 
problems on elementary classes are better dealt with in some 
non-elementary contexts (close to abstract elementary class), 
as if we would like 
during the proof to consider some derived other classes, those
contexts give you more freedom.  In particular this may apply to
\mr
\item "{$(a)$}"   main gap for $|T|^+$-saturated models (the parallel of
$(D,\lambda)$-sequence-homogeneous above in \scite{E53-nl.5c.2} and (d) of
(A) and discussion on it in (B))
\sn
\item "{$(b)$}"    the main gap for the class of models of $T$ for an
\ub{uncountable} first order $T$.
\endroster
\bn
Note that \cite{Sh:300}, Chapter II has tried to materialize this, but
that program is not finished.
\sn
\margintag{E53-nl.5c.3}\ub{\stag{E53-nl.5c.3} Problem}:  Similar questions for 
the number of pairwise non-elementarily embeddable
$(D,\lambda)$-sequence homogeneous models.
\bn
In the case of the class of models (not the class of
$\aleph_1$-saturated models) for countable first order theories,
those two problems were solved together.

There are many other interesting questions in this context.  An important
one, of a different character is: \nl
\margintag{E53-nl.5c.5}\ub{\stag{E53-nl.5c.5} Problem}: 1) [Hanf number for sequence homogeneous]

Given a cardinal $\kappa$, what is the first $\lambda$ such that: if $T$ is
a complete first order theory, $D\subseteq D(T) = \{\sftp(\bar a,
\emptyset,M):M$ a model of $T,\bar{a} \in {}^{\omega >} M\}$ and there is a
$(D,\lambda)$-sequence-homogeneous model, \ub{then} for 
every $\mu>\lambda$ there is a $(D,\mu)$-sequence homogeneous model.
\nl
2) Similarly for $\{\kappa$: in ${\frak K}$ we have amalgamation for
models of cardinality $< \kappa$ (and $\kappa \ge 
\text{ LST}({\frak K}) > \aleph_0\}$).
\nl
3) Similarly for $(\Bbb D,\lambda)$-model homogeneous models (see
   \sectioncite[\S3]{300b}). 
 
Toward this we may define semi-beautiful classses as in
\sectioncite[\S12]{705} (or \cite{Sh:87a}, \cite{Sh:87b}) 
replacing the stable ${\Cal P}^-(n)$-systems by an abstract 
notion, omitting uniqueness and
the definability of types and retaining existence.  Semi-excellent classes
seem like an effective version
of having amalgamation, so it certainly implies it; such properties may serve
as what we actually have to prove to solve the problem \scite{E53-nl.5c.5}
above.  We may have to use more complicated frames: say classes
${\frak K}_n$ so that $M \in {\frak K}_n$ is actually 
a ${\Cal P}^-(n)$-system of models from ${\frak K}$.  (See
more in \cite{Sh:842}). 

Recall that a class ${\frak K}$ of structures with fixed vocabulary
$\tau$ is called \ub{universal} if it is closed under isomorphisms,
and $M \in {\frak K}$ \ub{if and only if}
every finitely generated submodel of $M$ belongs to ${\frak K}$.  So 
not every elementary class is a universal class, but many
universal classes are not first order (e.g., locally finite
groups). This investigation leads (see \cite{Sh:300}, Chapter II) to 
classes with an axiomatized notion of non-forking and much of \cite{Sh:c} was
generalized, sometimes changing the context (a case of Thesis
\scite{E53-nl.5c.4}), but, e.g., still:
\bn
\margintag{E53-nl.5c.6}\ub{\stag{E53-nl.5c.6} Problem}:  Prove the main gap for the universal
context.
\bn
\margintag{E53-nl.6c.7}\ub{\stag{E53-nl.6c.7} Question}:  Can we in \cite{Sh:576},
i.e. \cite{Sh:E46} weaken the
``categorical in $\lambda^+$" to ``has a superlimit model in
$\lambda^+$"?  

See on this hopefully \cite{Sh:F888}.
\sn
\margintag{E53-nl.6c.11}\ub{\stag{E53-nl.6c.11} Question}:  Do we use a parallel of
\sectioncite[\S12]{705} with existential side for serious effect?  (See
more in \cite{Sh:842}).
\newpage

\head {\S5 Basic knowledge} \endhead  \resetall \sectno=5
 \spuriousreset
\bn
(A) What knowledge needed and dependency of the chapters

The chapters were written separately, hence for better or for worse
there are some repetitions, hopefully helping the reader if he likes
to read only parts of this book.

\chaptercite{705} depends on \chaptercite{600} and \cite{Sh:E46} depends
somewhat, e.g. on \sectioncite[\S1]{600}, but in other cases there are
no real dependency.

In fact, reading \chaptercite{600}, \chaptercite{705} requires little
knowledge of model theory, they are quite self-contained, in 
particular you do not
need to know \chaptercite{88r}, Chapter II; this apply also to Chapter II
and to \cite{Sh:E46}.
Of course, if a claim proves that the axioms of good $\lambda$-frames
are satisfied by the class of
models of a sentence $\psi$ in a logic you have not heard about, it
will be a little loss for you to ignore the claim (this occurs in
\sectioncite[\S3]{600}). 
Still much of the material is motivated by parallelism to what we know in
elementary (= first order) contexts.  
Let me stress that neither do we see any merit
in not using large model theoretic background nor was its elimination
an a priori aim, but there is no reason to hide this fact from a potential
reader who may feel otherwise.

Also the set theoretic knowledge required in \chaptercite{600}, \chaptercite{705}
is small; still we use cardinals and 
ordinals of course, induction on ordinals, 
cofinality of an ordinal, so regular cardinals, see here below for what
you need.  A priori it seemed that somewhat more is needed in the proof of
the non-structure theorems, i.e., showing a class with a so-called
``non-structure property" has many, complicated models so cannot have
a structure theory.  But we circumvent this by quoting \cite{Sh:838},
or you can say delaying the proof.  That is, we
carry the construction enough to give a reasonable argument.  So the
reader can just agree to believe; similarly in Chapter II and in
\cite{Sh:E46}.

In \cite{Sh:838} itself, we rely somewhat on basics of
\sectioncite[\S1]{600}, and in the applications (\cite[\S4]{Sh:838}) we
somewhat depend on the relevant knowledge and for
\cite[\S5-\S8]{Sh:838} we assume the basics of \sectioncite[\S2]{600}.  Also
\cite[\S9,\S10,\S11]{Sh:838} are set theoretic, mainly use results on
the weak diamond which we quote.

The situation is different in \chaptercite{88r}.  Still you can read \S1,
\S2, \S3 of it ignoring some claims but in \S4,\S5 the infinitary logics $\Bbb
L_{\omega_1,\omega}(\bold Q)$ and its relatives and basic theorems on them
are important.  

For \chaptercite{734} you need basic knowledge of
infinitary logics and Ehenfeucht-Mostowski models, and in
\sectioncite[\S4]{734} (the
main theorem) we use the definition of good $\lambda$-frame from
\sectioncite[\S2]{600}. 
\bn
(B) Some basic definitions and notation

We first deal with model theory and then with set theory.
\bigskip

\definition{\stag{E53-6.1} Definition}  1) A vocabulary $\tau$ is a set of
function symbols (denoted by $G,H,F$) and relation symbols, (denoted by
$P,Q,R)$ (= predicates), to each such symbol a number of places (=
arity) is assigned (by $\tau$) denoted by
arity$_\tau(F)$, arity$_\tau(P)$, respectively.
\nl
2) $M$ is a $\tau$-model or a $\tau$-structure 
for a vocabulary $\tau$ means that $M$ consists of:
\mr
\item "{$(a)$}"  its universe, $|M|$, a non-empty set
\sn
\item "{$(b)$}"  $P^M$, the interpretation of a predicate $P \in \tau$
and $P^M$ is an arity$_\tau(P)$-place relation on $|M|$
\sn
\item "{$(c)$}"  $F^M$, the interpretation of a function symbol $F \in
\tau$ and $F^M$ is an arity$_\tau(F)$-place function from $|M|$ to
$|M|$ in the case of arity $0,F^M$ is an individual constant.
\ermn
3) We agree $\tau$ is determined by $M$ and denote it by $\tau_M$.  If
$\tau_1 \subseteq \tau_2,M_2$ a $\tau_2$-model, then $M_1 = M_2
\restriction \tau_1$, the reduct is naturally defined.
\nl
4) The cardinality of $M,\|M\|$, is the cardinality, number of elements
of the universe $|M|$ of $M$.  We may write $a \in M$ instead of $a
\in |M|$ and $\langle a_i:i < \alpha \rangle \in M$ instead $i <
\alpha \Rightarrow a_i \in M$, i.e., $\bar a \in {}^\alpha|M|$.
\nl
5) Let $M \subseteq N$ mean that

$$
\tau_M = \tau_N,|M| \subseteq |N|,P^M = P^N \restriction |M|,
F^M = F^N \restriction |M|
$$
\mn
for every predicate $P \in \tau_M$ and 
for every function symbol $F \in \tau_M$.
\nl
6) If $N$ is a $\tau$-model and $A$ is a non-empty subset of $|M|$
closed under $F^N$ for each function symbol $F \in \tau$, \ub{then}
$N \restriction A$ is the unique $M \subseteq N$ with universe $A$.
\enddefinition
\bigskip

\definition{\stag{E53-6.1.1} Definition}  1)  $K$ denotes a class of
$\tau$-models closed under isomorphisms, for some vocabulary $\tau =
\tau_K$.
\nl
2) ${\frak K}$ denotes a pair $(K,\le_{\frak K})$; $K$ as above (with
$\tau_{\frak K} := \tau_K$) and $\le_{\frak K}$ is a two-place
relation on $K$ closed under isomorphisms such that $M \le_{\frak K} N
 \Rightarrow M \subseteq N$.
\nl
3) $f$ is a $\le_{\frak K}$-embedding of $M$ into $N$ \ub{when} for some $N'
\le_{\frak K} N,f$ is an isomorphism from $M$ onto $N'$.
\nl
4) $K$ is categorical in $\lambda$ \ub{if} $K$ has one and
only one model up to isomorphism of cardinality $\lambda$.  If ${\frak
K} = (K,\le_{\frak K})$ we may say ``${\frak K}$ is categorical in $\lambda$".
\enddefinition
\bigskip

\definition{\stag{E53-6.1.2} Definition}  1) For a class $K$ (or ${\frak
K}$) of $\tau_K$-models
\mr
\item "{$(a)$}"  $K_\lambda = \{M \in K:\|M\| = \lambda\}$
\sn
\item "{$(b)$}"  ${\frak K}_\lambda = (K_\lambda,\le_{\frak K}
\restriction K_\lambda)$
\sn
\item "{$(c)$}"  $\dot I(\lambda,K) = \dot I(\lambda,{\frak K}) =
|\{M/\cong:M \in K_\lambda\}|$ so $K$ (or ${\frak K}$) is
 categorical in $\lambda$ iff $\dot I(\lambda,K)=1$
\sn
\item "{$(d)$}"  $\dot I \dot E(\lambda,{\frak K}) = \sup\{\mu$: there
is a sequence $\langle M_\alpha:\alpha < \mu  \rangle$ of members of
$K_\lambda$ such that $M_\alpha$ is not $\le_{\frak K}$-embeddable into
$M_\beta$ for any distinct $\alpha,\beta < \mu\}$.  But writing 
$\dot I \dot E(\lambda,{\frak K}) \ge \mu$ we mean the supremum is obtained if
not said otherwise.
\sn
\item "{$(e)$}"  $M \in {\frak K}$ is $(\le_{\frak
K},\lambda)$-universal if every $N \in {\frak K}_\lambda$ can be
$\le_{\frak K}$-embedded into it.  If $\lambda = \|M\|$ we may write
$\le_{\frak K}$-universal.  If ${\frak K}$ is clear from the context
we may write $\lambda$-universal or universal (for ${\frak K}$).
\ermn
We end the model-theory part by defining logics (this is 
not needed for \chaptercite{600}, \chaptercite{705}, \cite{Sh:E46} and Chapter
II except some parts of \chaptercite{300a}).
\enddefinition
\bigskip

\definition{\stag{E53-6.2} Definition}  A logic ${\Cal L}$ consists of:
\mr
\item "{$(a)$}"  function ${\Cal L}(-)$ (actually a definition) giving
for every vocabulary $\tau$ a set of so-called formulas $\varphi(\bar
x),\bar x$ a sequence of free variables with no repetitions
\sn
\item "{$(b)$}"  $\models_{\Cal L}$, satisfaction relation, i.e., for
every vocabulary $\tau$ and $\varphi(\bar x) \in {\Cal L}(\tau)$ and
$\tau$-model $M$ and $\bar a \in {}^{\ell g(\bar x)}M$ we have ``$M
\models_{\Cal L} \varphi[\bar a]"$ or in words 
``$M$ satisfies $\varphi[\bar a]"$; holds or fails.
\endroster
\enddefinition
\bn
As for set theory
\endaside
\definition{\stag{E53-6s.1} Definition}  1) A power = number of elements of
a set, is identified with the first ordinal of this power, that is a
cardinal.  Such
ordinals are called cardinals, $\aleph_\alpha$ is the $\alpha$-th
infinite ordinal.
\nl
2) Cardinals are denoted by $\lambda,\mu,\kappa,\chi,\theta,\partial$
(infinite if not said otherwise).
\enddefinition
\bigskip

\definition{\stag{E53-6s.2} Definition}  0) Ordinals are denoted by
$\alpha,\beta,\gamma,\delta,\varepsilon,\zeta,\xi,i,j$, but, if not
said otherwise $\delta$ denotes a limit ordinal.
\nl
1) An ordinal $\alpha$ is a limit ordinal if $\alpha > 0$ and
$(\forall \beta < \alpha)[\beta + 1 < \alpha]$.
\nl
2) For an ordinal $\alpha$, cf$(\alpha)$, the cofinality of $\alpha$,
is min$\{\text{\rm otp}(u):u \subseteq \alpha$ is unbounded$\}$;
it is a regular cardinal (see below), we can define the cofinality for
linear orders and again get a regular cardinal.
\nl
3) A cardinal $\lambda$ is regular if cf$(\lambda) = \lambda$,
otherwise it is called singular.
\nl
4) If $\lambda = \aleph_\alpha$ then $\lambda^+ = \aleph_{\alpha +1}$,
the successor of $\lambda$, so $\lambda^{++} = \aleph_{\alpha
+2},\lambda^{+ \varepsilon} = \aleph_{\alpha + \varepsilon}$.
\enddefinition
\bn
Recall:
\proclaim{\stag{E53-6s.3} Claim}  1) If $\lambda$ is a regular cardinal,
$|{\Cal U}_t| < \lambda$ for $t \in I$ and $|I| < \lambda$ \ub{then}
$\cup\{{\Cal U}_t:t \in I\}$ has cardinality $< \lambda$.
\nl
2) $\lambda^+$ is regular for any $\lambda \ge \aleph_0$ but
$\lambda^{+ \delta}$ is singular if $\delta$ is a limit ordinal $<
\lambda$ (or just $< \lambda^{+ \delta}$), and, obviously, $\aleph_0$
is regular but e.g. $\aleph_\omega$ is singular, in fact
$\aleph_\delta > \delta \Rightarrow \aleph_\delta$ is singular, but
the inverse is false. 
\endproclaim
\bn
Sometimes we use (not essential)
\definition{\stag{E53-6s.14} Definition/Claim}  1) ${\Cal H}(\lambda)$ is
the set of $x$ such that there is a set $Y$ of cardinality $< \lambda$
which is transitive (i.e. $(\forall y)(y \in Y \Rightarrow y \subseteq
Y)$ and $x$ belongs to $\lambda$.
\nl
2) Every $x$ belongs to ${\Cal H}(\lambda)$ for some $x$.
\nl
So for some purpose we can look at ${\Cal H}(\lambda)$ instead of the
universe of all sets.
\enddefinition
\newpage

\head {\S6 Index of symbols\footnotemark} \endhead  \resetall \sectno=6
\footnotetext{some will be used only in subsequent works; in
particular concerning forcing}
 \spuriousreset
\bn
$a \quad$ member of a model
\sn
$A \quad$ set of elements of model
\sn
${\frak A} \quad$ a ``complicated" model
\sn
$b \quad$ member of a model
\sn
$B \quad$ set of members of models
\sn
${\frak B} \quad$ a ``complicated" model
\sn
$c \quad$ member of model (also individual constant)
\sn
$\bold c \quad$ colouring, mainly \cite{Sh:838}
\sn
$C \quad$ set or elements of models \ub{or} a club
\sn
${\Cal C} \quad$ club of $[A]^{< \lambda}$,
\sn
${\frak C} \quad$ a complicated model, \ub{or} a monster
\sn
$d \quad$ member of model
\sn
$\bold d \quad$ expanded $I$-system, \sectioncite[\S12]{705}; ${\frak
u}$-free rectangle or triangle in \cite{Sh:838}
\sn
$D \quad$ diagram; set of $(< \omega)$-types in the first order sense
realized in a model, \nl

\qquad \,\chaptercite{88r}, \chaptercite{300b}
\sn
$\bold D \quad$ a function whose values are diagrams, \chaptercite{88r},
\chaptercite{300b} 
\sn
${\Bbb D} \quad$ diagram for model homogeneity, \chaptercite{88r}, so set of
isomorphism types of 
\nl

\quad models, also \chaptercite{300b}
\sn
${\frak D}$ a set of $\Bbb D$'s, \chaptercite{300b}
\sn
${\Cal D} \quad$ filter
\sn
${\Cal D}_\lambda \quad$ club filter on the regular cardinal $\lambda
> \aleph_0$
\sn
$e \quad$ element of a model \ub{or} a club
\sn
$\bold e \quad$ expanded $I$-system (used in continuations),
\sectioncite[\S12]{705};
\nl

\quad ${\frak u}$-free rectangle or triangle in \cite{Sh:838}
\sn
$E \quad$ a club 
\sn
$\Bbb E \quad$ filter
\sn
${\Cal E} \quad$ an equivalence relations, (e.g. 
${\Cal E}_M,{\Cal E}^{\text{at}}_M$ in \marginbf{!!}{\cprefix{600}.\scite{600-0.12}} for 
definition of type
\nl

\,\,\,and ${\Cal E}_{{\frak K},\chi},{\Cal E}^{\text{mat}}_{{\frak K},\chi}$
in \sectioncite[\S3]{300b})
\sn
$f \quad$ function (e.g., isormorphism, embedding usually)
\sn
$\bold f \quad$ function (\cite{Sh:838} in $(\bar M,\bar{\bold
J},\bold f) \in K^{3,\text{qt}}_{\frak u}$, also in \sectioncite[\S5]{600},
$(\bar M,\bar{\bold f}),(\bar M,\bar{\bold J},\bold f))$
\sn
$F \quad$ function symbol
\sn
$\Bbb F \quad$  amalgamation choice function
(\cite{Sh:838} also see \cite[\S3]{Sh:576}) 
\sn
$\bold F \quad$ function (complicated, mainly it witnesses a model
being limit, \sectioncite[\S3]{88r})
\sn
$g \quad$ function
\sn
${\frak g} \quad$ witness for almost every $(\bar M,\bar{\bold
J},\bold f)$ see \cite[1a.43-1a.51]{Sh:838}
\sn
$G \quad$ function symbol 
\sn
$\Game \quad$ game 
\sn
$h \quad$ function
\sn
${\frak h} \quad$ witnesses for almost every $(\bar M,\bar{\bold
J},\bold f) \in K^{\text{qt}}_{\frak u}$, see \cite[c.4A-c.4D]{Sh:838}
or

\qquad  \cite[1a.43-1a.51]{Sh:838}
\sn
$H \quad$ function symbol 
\sn
${\Cal H} \quad$ in ${\Cal H}(\lambda)$, rare here see \scite{E53-6s.14}
\sn
$i \quad$ ordinal/natural number
\sn
$I \quad$ linear order, partial order or index set
\sn
$\dot I \quad \dot I(\lambda,K)$, numbers on non-isomorphic models; 
\nl

$\dot I \dot E(\lambda,K)$ (see \chaptercite{88r}), also $\dot I(K)$, see
\cite{Sh:838} 
\sn
$\bold I \quad$ set of sequences or elements from a model, in particular:

\,\,\,$\bold I_{M,N} = \{c \in N:\ortp_{\frak s}(c,M,N) \in 
{\Cal S}^{\text{bs}}_{\frak s}(M)\}$, see \chaptercite{600}, \chaptercite{705}
\sn
$\check I[\lambda] \quad$ a specific normal ideal, see
\sectioncite[\S0]{88r}, marginal here
\sn
$\Bbb I \quad$ ideal 
\sn
${\Cal I}$ \,\, predense set in a forcing $\Bbb P$, very rare here 
\sn
$j \quad$ ordinal/natural number
\sn
$J \quad$ linear order, index set, \chaptercite{88r}
\sn
$\bold J \quad$ set of sequences or elements from a model
\sn
${\Bbb J} \quad$ ideal 
\sn
${\Cal J}$ \,\, predense set in a forcing $\Bbb P$, very rare here
\sn
$k \quad$ natural number
\sn
$K \quad$ class of model of a fix vocabulary $\tau_{\frak
K},K_\lambda$ is $\{M \in K:\|M\| = \lambda\}$
\sn
${\frak K} \quad$ is $(K,\le_{\frak K})$, usually \aec
\sn
$K^{3,x}_{\frak s}  \quad$ for $x = \{\text{bs,uq,pr,qr,vq,bu}\}$,
appropriate set of triples $(M,N,a)$ or $(M,N,\bold I)$,

\qquad  see \chaptercite{600}, \chaptercite{705}
\sn
$K^{3,\text{na}}_\lambda \quad$ for triples $(M,N,a)$, see
\cite{Sh:E46}
\sn
$K^{3,x}_{\frak u}  \quad$ set of triples
$(M,N,\bold J) \in \text{ FR}^\ell_{\frak u}$, see
\cite{Sh:838}
\sn
$\ell \quad$ natural number
\sn
$L \quad$ language (set of formulas, e.g., ${\Cal L}(\tau)$ but also
subsets of ${\Cal L}(\tau)$ which normally 
\nl 

\hskip10pt are closed under subformulas and first order operations), 
used in \chaptercite{88r}.
\sn
LST $\quad$ L\"owenheim-Skolem-Tarski numbers, mainly LST$({\frak K}) = 
\text{ LST}_{\frak K}$
\sn
${\Cal L} \quad$ logic, i.e., a function such that ${\Cal L}(\tau)$ is a
language for vocabulary $\tau$ (but also 
\nl

\hskip10pt a language mainly ${\Cal L}$ a
fragment of $\Bbb L_{\lambda^+,\omega}$, i.e., a subset closed
\nl

\hskip10pt under subformulas and the finitary operations)
\sn
$\prec_{\Cal L} \quad$ is used for $M \prec_{\Cal L} N$ iff $M \subseteq
N$ and for every $\varphi(\bar x) \in {\Cal L}(\tau_M)$ and $\bar a
\in {}^{\ell g(\bar x)}M$ 
\nl

\hskip10pt we have $M \models \varphi[\bar a]
\Leftrightarrow N \models \varphi[\bar a]$
\sn
$\Bbb L \quad$ first order logic and $\Bbb L_{\lambda,\kappa},\Bbb
L^\ell_{\lambda,\kappa}$, see \chaptercite{88r} so $\varphi(\bar x) \in
\Bbb L_{\lambda,\kappa}$ has \nl

$< \kappa$ free variables
\sn
$\bold L \quad$ the constructible universe
\sn
$m \quad$ natural number
\sn
$\bold m$ an $I$-system in \sectioncite[\S12]{705}
\sn
$M \quad$ model
\sn
$\bold M \quad$ complicated object, see \cite[\S3,\S4]{Sh:E46}
\sn
$n \quad$ natural number
\sn
$\bold n$ an $I$-system in \sectioncite[\S12]{705}, for continuation and in
\chaptercite{300f} 
\sn
$N \quad$ model
\sn
$\Bbb N \quad$ the natural numbers
\sn
$p \quad$ type
\sn
$\bold p \quad$ member of ${\Bbb P}$, a forcing condition, very rare here
\sn
$P \quad$ predicate
\sn
${\Cal P} \quad$ power set, family of sets,
\sn
$\bold P \quad$ family of types, \chaptercite{705}
\sn
$\Bbb P \quad$ forcing notion, very rare here 
\sn
$q \quad$ type
\sn
$\bold q \quad$ forcing condition, very rare here 
\sn
$Q \quad$ predicate
\sn
$\bold Q \quad$  a quantifier written $(\bold Q x)\varphi$, see
\chaptercite{88r}, if clear from the context means 
\nl

\hskip10pt $\bold Q^{\text{car}}_{\ge \aleph_1}$
\sn
$\bold Q^{\text{car}}_{\ge \kappa} \quad$ the quantifier there are $\ge
\kappa$ many
\sn
$\Bbb Q \quad$ the rationals
\sn
$r \quad$ type 
\sn
$\bold r \quad$  forcing condition, very rare here 
\sn
$R \quad$ predicate
\sn 
$\Bbb R \quad$ reals 
\sn
$s \quad$ member of $I,J$
\sn
${\frak s} \quad$ frame
\sn
$S \quad$ set of ordinals, stationary set many times
\sn
${\Cal S} \quad$ ${\Cal S}_{\frak K}(M)$ is a set of types in the sense of
orbits, ${\Cal S}^{\text{bs}}_{\frak s}(M)$ the basic types
\nl

\hskip10pt  (there are some alternatives to bs)
\sn
$\bold S \quad \bold S^\alpha_L(A,M)$ set of complete $(L,\alpha)$-types
over $M$, so a set of formulas,
\nl

\hskip10pt  used when we are dealing with a logic ${\Cal L}$, may use $\bold
S^\alpha_{\Cal L}(A,M)$ 
\sn
$\Bbb S \quad \Bbb S(M)$ is a set of pseudo types, are neither set of
formulas nor orbits, but formal
\nl

\hskip10pt  non-forking extension (for continuations, see \cite{Sh:842})
\sn
$t \quad$ member of $I,J$
\sn
$\sftp \quad$ type as set of formulas
\sn
$\ortp \quad$ type as an orbit, an equivalence class under mapping
\sn
$\bold t \quad$ type function
\sn
${\frak t} \quad$ frame
\sn
$T \quad$ first order theory, usually complete
\sn
${\Cal T} \quad$ a tree
\sn
$u \quad$ a set
\sn
${\frak u} \quad$ a nice construction framework, in \cite{Sh:838}
\sn
unif \quad in $\mu_{\text{unif}}(\lambda,2^{< \lambda})$, see
\marginbf{!!}{\cprefix{88r}.\scite{88r-0.wD}} or \cite[0z.6]{Sh:838}(6)
\sn
$U \quad$ a set
\sn
${\Cal U} \quad$ a set
\sn
$v \quad$ a set
\sn
$V \quad$ a set 
\sn
$\bold V \quad$ universe of set theory
\sn
$w \quad$ a set
\sn
$W \quad$ a set (usually of ordinals)
\sn
${\Cal W} \quad$ a class of triples $(N,\bar M,\bar{\bold J})$; see
\sectioncite[\S7]{705}
\sn
wd \quad in $\mu_{\text{wd}}(\lambda)$ see \sectioncite[\S0]{88r},
\cite[\S0]{Sh:838}
\sn
WDmId$_\lambda \quad$ the weak diamond ideal, see \marginbf{!!}{\cprefix{88r}.\scite{88r-0.wD}}
\sn
$x \quad$ variable (or element)
\sn
$\bold x \quad$ complicated object, in \cite{Sh:838} such that is a
sequence 

\qquad $\langle(\bar M^\alpha,\bar{\bold J}^\alpha,\bold
f^\alpha):\alpha < \alpha(*)\rangle$
\sn
$X \quad$ set
\sn
$y \quad$ variable
\sn
$\bold y \quad$ like $\bold x$
\sn
$Y \quad$ set
\sn
${\Cal Y} \quad$ a high order variable (see \sectioncite[\S3]{88r})
\sn
$z \quad$ variable
\sn
$Z \quad$ set
\sn
$\Bbb Z \quad$ the integers
\bn
\ub{Greek Letters}:
\bn
$\alpha \quad$ ordinal
\sn
$\beta \quad$ ordinal
\sn
$\gamma \quad$ ordinal
\sn
$\Gamma \quad$ various things; in \cite{Sh:E46} a set of models or types
\sn
$\delta \quad$ ordinal, limit if not clear otherwise
\sn
$\partial \quad$ cardinal
\sn
$\Delta \quad$ set of formulas (may be used for symmetric difference)
\sn
$\epsilon \quad$ ordinal
\sn
$\varepsilon \quad$ ordinal
\sn
$\zeta \quad$ ordinal
\sn
$\eta \quad$ sequence, usually of ordinals
\sn
$\theta \quad$ cardinal, infinite if not clear otherwise
\sn
$\vartheta \quad$ a formula, very rare
\sn
$\Theta \quad$ set of cardinals/class of cardinals
\sn
$\iota \quad$ ordinal (sometimes a natural number)
\sn
$\kappa \quad$ cardinal, infinite if not clear otherwise
\sn
$\lambda \quad$ cardinal, infinite if not clear otherwise
\sn
$\lambda({\frak K}) \quad$ is the L.S.T.-number of an \aec $(\ge
|\tau_{\frak K}|$ for simplicity), 
\nl

\qquad rare
\sn
$\Lambda \quad$ set of formulas, used in \chaptercite{734}, Chapter I
\sn
$\mu \quad$ cardinal, infinite if not said otherwise
\sn
$\nu \quad$ sequence, usually of ordinals
\sn
$\sigma \quad$ a term (in a vocabulary $\tau$)
\sn
$\Sigma \quad$ sum
\sn
$\pi \quad$ permutation
\sn
$\Pi \quad$ product
\sn
$\rho \quad$ sequence, usually of ordinals
\sn
$\varrho \quad$ sequence, usually of ordinals
\sn
$\tau \quad$ vocabulary (so ${\Cal L}(\tau),\Bbb L(\tau),\Bbb
L_{\lambda,\mu}(\tau)$ are languages)
\sn
$\xi \quad$ ordinal
\sn
$\Xi \quad$ a complicated object
\sn
$\Upsilon \quad$ ordinal and other objects
\sn
$\chi \quad$ cardinal, infinite if not said otherwise
\sn
$\varphi \quad$ formula
\sn
$\Phi \quad$ blueprint for EM-models
\sn
$\psi \quad$ formula
\sn
$\Psi \quad$ blueprint for EM-models
\sn
$\omega \quad$ the first infinite ordinal
\sn
$\Omega \quad$ a complicated object
\newpage


\nocite{ignore-this-bibtex-warning} 
\newpage
    
REFERENCES.  
\bibliographystyle{lit-plain}
\bibliography{lista,listb,listx,listf,liste}

\def\germ{\frak} \def\scr{\cal} \ifx\documentclass\undefinedcs
  \def\bf{\fam\bffam\tenbf}\def\rm{\fam0\tenrm}\fi 
  \def\defaultdefine#1#2{\expandafter\ifx\csname#1\endcsname\relax
  \expandafter\def\csname#1\endcsname{#2}\fi} \defaultdefine{Bbb}{\bf}
  \defaultdefine{frak}{\bf} \defaultdefine{=}{\B} 
  \defaultdefine{mathfrak}{\frak} \defaultdefine{mathbb}{\bf}
  \defaultdefine{mathcal}{\cal}
  \defaultdefine{beth}{BETH}\defaultdefine{cal}{\bf} \def\bbfI{{\Bbb I}}
  \def\mbox{\hbox} \def\text{\hbox} \def\om{\omega} \def\Cal#1{{\bf #1}}
  \def\pcf{pcf} \defaultdefine{cf}{cf} \defaultdefine{reals}{{\Bbb R}}
  \defaultdefine{real}{{\Bbb R}} \def\restriction{{|}} \def\club{CLUB}
  \def\w{\omega} \def\exist{\exists} \def\se{{\germ se}} \def\bb{{\bf b}}
  \def\equivalence{\equiv} \let\lt< \let\gt>
  \def\implies{\Rightarrow}\def\mathfrak{\bf}\def\germ{\frak} \def\scr{\cal}
  \ifx\documentclass\undefinedcs
  \def\bf{\fam\bffam\tenbf}\def\rm{\fam0\tenrm}\fi 
  \def\defaultdefine#1#2{\expandafter\ifx\csname#1\endcsname\relax
  \expandafter\def\csname#1\endcsname{#2}\fi} \defaultdefine{Bbb}{\bf}
  \defaultdefine{frak}{\bf} \defaultdefine{=}{\B} 
  \defaultdefine{mathfrak}{\frak} \defaultdefine{mathbb}{\bf}
  \defaultdefine{mathcal}{\cal}
  \defaultdefine{beth}{BETH}\defaultdefine{cal}{\bf} \def\bbfI{{\Bbb I}}
  \def\mbox{\hbox} \def\text{\hbox} \def\om{\omega} \def\Cal#1{{\bf #1}}
  \def\pcf{pcf} \defaultdefine{cf}{cf} \defaultdefine{reals}{{\Bbb R}}
  \defaultdefine{real}{{\Bbb R}} \def\restriction{{|}} \def\club{CLUB}
  \def\w{\omega} \def\exist{\exists} \def\se{{\germ se}} \def\bb{{\bf b}}
  \def\equivalence{\equiv} \let\lt< \let\gt>
\begin{thebibliography}{HShT 428}
\makeatletter \renewcommand{\@biblabel}[1]{[#1]} \makeatother
\def\eprintfootnotetext{References of the form {\tt math.XX/$\cdots$}
 refer to {\tt arXiv.org} }
\ifx\documentstyle\undefinedcontrolsequence
   \def\anyfootnote{\footnote{*}}
   \else\def\anyfootnote{\footnote}\fi
\def\eprintfn{\ifEprint\anyfootnote{\eprintfootnotetext}\fi\Eprintfalse }
\newif\ifEprint  \Eprinttrue

\bibitem[Bal88]{Bal88}John Baldwin.
\newblock {\em Fundamentals of Stability Theory}.
\newblock Perspectives in Mathematical Logic. Springer-Verlag, Berlin, 1988.

\bibitem[Bal0x]{Bal0x}John Baldwin.
\newblock {\em {Categoricity}}, volume to appear.
\newblock 200x.

\bibitem[BlSh 862]{BlSh:862}John Baldwin and Saharon Shelah.
\newblock {Examples of non-locality}.
\newblock {\em Journal of Symbolic Logic}, {\bf 73}:765--782, 2008.

\bibitem[Bl85]{Bl85}John~T. Baldwin.
\newblock {Definable second order quantifiers}.
\newblock In J.~Barwise and S.~Feferman, editors, {\em {Model Theoretic
  Logics}}, Perspectives in Mathematical Logic, chapter XII, pages 445--477.
  Springer-Verlag, New York Berlin Heidelberg Tokyo, 1985.

\bibitem[BETp06]{BETp06}John~T. Baldwin, Paul Eklof, and Jan Trlifaj.
\newblock {$N$ perp as an AEC}.
\newblock {\em Preprint}, 2006, revised 2007.

\bibitem[BKV0x]{BKV0x}John~T. Baldwin, David~W. Kueker, and Monica VanDieren.
\newblock {Upward Stability Transfer Theorem for Tame Abstract Elementary
  Classes}.
\newblock {\em Preprint}, 2004.

\bibitem[BLSh 464]{BLSh:464}John~T. Baldwin, Michael~C. Laskowski, and Saharon
  Shelah.
\newblock {Forcing Isomorphism}.
\newblock {\em {Journal of Symbolic Logic}}, {\bf 58}:1291--1301, 1993.
\newblock math.LO/9301208.

\bibitem[BlSh 156]{BlSh:156}John~T. Baldwin and Saharon Shelah.
\newblock {Second-order quantifiers and the complexity of theories}.
\newblock {\em {Notre Dame Journal of Formal Logic}}, {\bf 26}:229--303, 1985.
\newblock Proceedings of the 1980/1 Jerusalem Model Theory year.

\bibitem[BlSh 330]{BlSh:330}John~T. Baldwin and Saharon Shelah.
\newblock {The primal framework. I}.
\newblock {\em {Annals of Pure and Applied Logic}}, {\bf 46}:235--264, 1990.
\newblock math.LO/9201241.

\bibitem[BlSh 360]{BlSh:360}John~T. Baldwin and Saharon Shelah.
\newblock {The primal framework. II. Smoothness}.
\newblock {\em {Annals of Pure and Applied Logic}}, {\bf 55}:1--34, 1991.
\newblock Note: See also 360a below. math.LO/9201246.

\bibitem[BlSh 393]{BlSh:393}John~T. Baldwin and Saharon Shelah.
\newblock {Abstract classes with few models have `homogeneous-universal'
  models}.
\newblock {\em {Journal of Symbolic Logic}}, {\bf 60}:246--265, 1995.
\newblock math.LO/9502231.

\bibitem[BaFe85]{BaFe85}Jon Barwise and Solomon~Feferman (editors).
\newblock {\em {Model-theoretic logics}}.
\newblock Perspectives in Mathematical Logic. Springer Verlag, Heidelberg-New
  York, 1985.

\bibitem[BY0y]{BY0y}Itay Ben-Yaacov.
\newblock {Uncountable dense categoricity in cats}.
\newblock {\em J. Symbolic Logic}, {\bf 70}:829--860, 2005.

\bibitem[BeUs0x]{BeUs0x}Itay Ben-Yaacov and Alex Usvyatsov.
\newblock {Logic of metric spaces and Hausdorff CATs}.
\newblock {\em In preparation}.

\bibitem[BoNe94]{BoNe94}Alexandre Borovik and Ali Nesin.
\newblock {\em {Groups of finite Morley rank}}, volume~26 of {\em Oxford Logic
  Guide}.
\newblock The Clarendon Press, Oxford University Press, New York, 1994.

\bibitem[ChKe66]{ChKe66}Chen-Chung Chang and Jerome~H. Keisler.
\newblock {\em {Continuous Model Theory}}, volume~58 of {\em Annals of
  Mathematics Studies}.
\newblock Princeton University Press, Princeton, NJ, 1966.

\bibitem[ChKe62]{ChKe62}Chen chung Chang and Jerome~H. Keisler.
\newblock {Model theories with truth values in a uniform space}.
\newblock {\em Bulletin of the American Mathematical Society}, {\bf
  68}:107--109, 1962.

\bibitem[CoSh:919]{CoSh:919}Moran Cohen and Saharon Shelah.
\newblock {Stable theories and representation over sets}.
\newblock {\em preprint}.

\bibitem[Eh57]{Eh57}Andrzej Ehrenfeucht.
\newblock On theories categorical in power.
\newblock {\em {Fundamenta Mathematicae}}, {\bf 44}:241--248, 1957.

\bibitem[Fr75]{Fr75}Harvey Friedman.
\newblock {One hundred and two problems in mathematical logic}.
\newblock {\em Journal of Symbolic Logic}, {\bf 40}:113--129, 1975.

\bibitem[GbTl06]{GbTl06}R\"udiger G\"obel and Jan Trlifaj.
\newblock {\em {Approximations and endomorphism algebras of modules}},
  volume~41 of {\em de Gruyter Expositions in Mathematics}.
\newblock Walter de Gruyter, Berlin, 2006.

\bibitem[Gr91]{Gr91}Rami Grossberg.
\newblock { On chains of relatively saturated submodels of a model without the
  order property}.
\newblock {\em Journal of Symbolic Logic}, {\bf 56}:124--128, 1991.

\bibitem[GrHa89]{GrHa89}Rami Grossberg and Bradd Hart.
\newblock {The classification of excellent classes}.
\newblock {\em Journal of Symbolic Logic}, {\bf 54}:1359--1381, 1989.

\bibitem[GIL02]{GIL02}Rami Grossberg, Jose Iovino, and Olivier Lessmann.
\newblock {A primer of simple theories}.
\newblock {\em Archive for Mathematical Logic}, {\bf 41}:541--580, 2002.

\bibitem[GrLe0x]{GrLe0x}Rami Grossberg and Olivier Lessmann.
\newblock {The main gap for totally transcendental diagrams and abstract
  decomposition theorems}.
\newblock {\em Preprint}.

\bibitem[GrLe00a]{GrLe00a}Rami Grossberg and Olivier Lessmann.
\newblock {Dependence relation in pregeometries}.
\newblock {\em Algebra Universalis}, {\bf 44}:199--216, 2000.

\bibitem[GrLe02]{GrLe02}Rami Grossberg and Olivier Lessmann.
\newblock {Shelah's stability spectrum and homogeneity spectrum in finite
  diagrams}.
\newblock {\em Archive for Mathematical Logic}, {\bf 41}:1--31, 2002.

\bibitem[GrSh 259]{GrSh:259}Rami Grossberg and Saharon Shelah.
\newblock {On Hanf numbers of the infinitary order property}.
\newblock {\em {Mathematica Japonica}}, {\bf submitted}.
\newblock math.LO/9809196.

\bibitem[GrSh 238]{GrSh:238}Rami Grossberg and Saharon Shelah.
\newblock {A nonstructure theorem for an infinitary theory which has the
  unsuperstability property}.
\newblock {\em {Illinois Journal of Mathematics}}, {\bf 30}:364--390, 1986.
\newblock {Volume dedicated to the memory of W.W.~Boone; ed. Appel, K., Higman,
  G., Robinson, D. and Jockush, C.}

\bibitem[GrSh 222]{GrSh:222}Rami Grossberg and Saharon Shelah.
\newblock {On the number of nonisomorphic models of an infinitary theory which
  has the infinitary order property. I}.
\newblock {\em {The Journal of Symbolic Logic}}, {\bf 51}:302--322, 1986.

\bibitem[GrVa0xa]{GrVa0xa}Rami Grossberg and Monica VanDieren.
\newblock {Galois-stbility for Tame Abstract Elementary Classes}.
\newblock {\em submitted}.

\bibitem[GrVa0xb]{GrVa0xb}Rami Grossberg and Monica VanDieren.
\newblock {Upward Categoricity Transfer Theorem for Tame Abstract Elementary
  Classes}.
\newblock {\em submitted}.

\bibitem[HHL00]{HHL00}Bradd Hart, Ehud Hrushovski, and Michael~C. Laskowski.
\newblock The uncountable spectra of countable theories.
\newblock {\em Annals of Mathematics}, {\bf 152}:207--257, 2000.

\bibitem[HaSh 323]{HaSh:323}Bradd Hart and Saharon Shelah.
\newblock {Categoricity over $P$ for first order $T$ or categoricity for
  $\phi\in{\rm L}_ {\omega_ 1\omega}$ can stop at $\aleph_ k$ while holding for
  $\aleph_ 0,\cdots,\aleph_ {k-1}$}.
\newblock {\em {Israel Journal of Mathematics}}, {\bf 70}:219--235, 1990.
\newblock math.LO/9201240.

\bibitem[He74]{He74}C.~Ward Henson.
\newblock {The isomorphism property in nonstandard analysis and its use in the
  theory of Banach spaces}.
\newblock {\em Journal of Symbolic Logic}, {\bf 39}:717--731, 1974.

\bibitem[HeIo02]{HeIo02}C.~Ward Henson and Jose Iovino.
\newblock {Ultraproducts in analysis}.
\newblock In {\em Analysis and logic (Mons, 1997)}, volume 262 of {\em London
  Math. Soc. Lecture Note Ser.}, pages 1--110. Cambridge Univ. Press,
  Cambridge, 2002.

\bibitem[He92]{He92}A.~Hernandez.
\newblock {\em {On $\omega_1$--saturated models of stable theories}}.
\newblock PhD thesis, Univ. of Calif. Berkeley, 1992.
\newblock Advisor: Leo Harrington.

\bibitem[Hy98]{Hy98}Tapani Hyttinen.
\newblock {Generalizing Morley's theorem}.
\newblock {\em Mathematical Logic Quarterly}, {\bf 44}:176--184, 1998.

\bibitem[HySh 474]{HySh:474}Tapani Hyttinen and Saharon Shelah.
\newblock {Constructing strongly equivalent nonisomorphic models for
  unsuperstable theories, Part A}.
\newblock {\em {Journal of Symbolic Logic}}, {\bf 59}:984--996, 1994.
\newblock math.LO/0406587.

\bibitem[HySh 529]{HySh:529}Tapani Hyttinen and Saharon Shelah.
\newblock {Constructing strongly equivalent nonisomorphic models for
  unsuperstable theories. Part B}.
\newblock {\em {Journal of Symbolic Logic}}, {\bf 60}:1260--1272, 1995.
\newblock math.LO/9202205.

\bibitem[HySh 632]{HySh:632}Tapani Hyttinen and Saharon Shelah.
\newblock {On the Number of Elementary Submodels of an Unsuperstable
  Homogeneous Structure}.
\newblock {\em {Mathematical Logic Quarterly}}, {\bf 44}:{354--358}, 1998.
\newblock math.LO/9702228.

\bibitem[HySh 602]{HySh:602}Tapani Hyttinen and Saharon Shelah.
\newblock {Constructing strongly equivalent nonisomorphic models for
  unsuperstable theories, Part C}.
\newblock {\em {Journal of Symbolic Logic}}, {\bf 64}:634--642, 1999.
\newblock math.LO/9709229.

\bibitem[HySh 629]{HySh:629}Tapani Hyttinen and Saharon Shelah.
\newblock {Strong splitting in stable homogeneous models}.
\newblock {\em {Annals of Pure and Applied Logic}}, {\bf 103}:201--228, 2000.
\newblock math.LO/9911229.

\bibitem[HySh 676]{HySh:676}Tapani Hyttinen and Saharon Shelah.
\newblock {Main gap for locally saturated elementary submodels of a homogeneous
  structure}.
\newblock {\em Journal of Symbolic Logic}, {\bf 66}:1286--1302, 2001, no.3.
\newblock math.LO/9804157.

\bibitem[HShT 428]{HShT:428}Tapani Hyttinen, Saharon Shelah, and Heikki Tuuri.
\newblock {Remarks on Strong Nonstructure Theorems}.
\newblock {\em {Notre Dame Journal of Formal Logic}}, {\bf 34}:157--168, 1993.

\bibitem[HyTu91]{HyTu91}Tapani Hyttinen and Heikki Tuuri.
\newblock {Constructing strongly equivalent nonisomorphic models for unstable
  theories}.
\newblock {\em Annals Pure and Applied Logic}, {\bf 52}:203--248, 1991.

\bibitem[JrSh 875]{JrSh:875}Adi Jarden and Saharon Shelah.
\newblock {Good frames minus stability}.
\newblock {\em Preprint}.

\bibitem[Jn56]{Jn56}Bjarni J{\'o}nsson.
\newblock {Universal relational systems}.
\newblock {\em Mathematica Scandinavica}, {\bf 4}:193--208, 1956.

\bibitem[Jn60]{Jn60}Bjarni J{\'o}nsson.
\newblock { Homogeneous universal relational systems}.
\newblock {\em Mathematica Scandinavica}, {\bf 8}:137--142, 1960.

\bibitem[KM67]{KM67}H.~Jerome Keisler and Michael~D. Morley.
\newblock {On the number of homogeneous models of a given power}.
\newblock {\em Israel Journal of Mathematics}, {\bf 5}:73--78, 1967.

\bibitem[Ke70]{Ke70}Jerome~H. Keisler.
\newblock {Logic with the quantifier "there exist uncountably many"}.
\newblock {\em Annals of Mathematical Logic}, {\bf 1}:1--93, 1970.

\bibitem[Ke71]{Ke71}Jerome~H. Keisler.
\newblock {\em {Model theory for infinitary logic. Logic with countable
  conjunctions and finite quantifiers}}, volume~62 of {\em {Studies in Logic
  and the Foundations of Mathematics}}.
\newblock North--Holland Publishing Co., Amsterdam--London, 1971.

\bibitem[KiPi98]{KiPi98}Byunghan Kim and Anand Pillay.
\newblock {From stability to simplicity}.
\newblock {\em Bull. Symbolic Logic}, {\bf 4}:17--36, 1998.

\bibitem[KlSh 362]{KlSh:362}Oren Kolman and Saharon Shelah.
\newblock {Categoricity of Theories in $L_{\kappa,\omega}$, when $\kappa$ is a
  measurable cardinal. Part 1}.
\newblock {\em {Fundamenta Mathematicae}}, {\bf 151}:209--240, 1996.
\newblock math.LO/9602216.

\bibitem[LwSh 871]{LwSh:871}Michael~C. Laskowski and Saharon Shelah.
\newblock {Karp height of models of stable theories}.
\newblock 0711.3043.

\bibitem[LwSh 489]{LwSh:489}Michael~C. Laskowski and Saharon Shelah.
\newblock {On the existence of atomic models}.
\newblock {\em {Journal of Symbolic Logic}}, {\bf 58}:1189--1194, 1993.
\newblock math.LO/9301210.

\bibitem[LwSh 560]{LwSh:560}Michael~C. Laskowski and Saharon Shelah.
\newblock {The Karp complexity of unstable classes}.
\newblock {\em {Archive for Mathematical Logic}}, {\bf {40}}:69--88, 2001.
\newblock math.LO/0011167.

\bibitem[LwSh 687]{LwSh:687}Michael~C. Laskowski and Saharon Shelah.
\newblock {Karp complexity and classes with the independence property}.
\newblock {\em Annals of Pure and Applied Logic}, {\bf 120}:263--283, 2003.
\newblock math.LO/0303345.

\bibitem[Le0x]{Le0x}Olivier Lessmann.
\newblock {Abstract group configuration}.
\newblock {\em Preprint}.

\bibitem[Le0y]{Le0y}Olivier Lessmann.
\newblock {Pregeometries in finite diagrams}.
\newblock {\em Preprint}.

\bibitem[MaSh 285]{MaSh:285}Michael Makkai and Saharon Shelah.
\newblock {Categoricity of theories in $L_ {\kappa\omega},$ with $\kappa$ a
  compact cardinal}.
\newblock {\em {Annals of Pure and Applied Logic}}, {\bf 47}:41--97, 1990.

\bibitem[Mw85]{Mw85}Johann~A. Makowsky.
\newblock Compactnes, embeddings and definability.
\newblock In J.~Barwise and S.~Feferman, editors, {\em Model-Theoretic Logics},
  pages 645--716. Springer-Verlag, 1985.

\bibitem[MkSh 366]{MkSh:366}Alan~H. Mekler and Saharon Shelah.
\newblock {Almost free algebras }.
\newblock {\em {Israel Journal of Mathematics}}, {\bf 89}:237--259, 1995.
\newblock math.LO/9408213.

\bibitem[MoVa62]{MoVa62}M.~D. Morley and R.~L. Vaught.
\newblock {Homogeneous and universal models}.
\newblock {\em Mathematica Scandinavica}, {\bf 11}:37--57, 1962.

\bibitem[Mo65]{Mo65}Michael Morley.
\newblock {Categoricity in power}.
\newblock {\em Transaction of the American Mathematical Society}, {\bf
  114}:514--538, 1965.

\bibitem[Pi0x]{Pi0x}Anand Pillay.
\newblock {Forking in the category of existentially closed structures}.
\newblock In {\em {Connections between model theory and algebraic and analytic
  geometry}}, volume~6 of {\em Quad. Mat.}, pages 23--42. Dept. Math., Seconda
  Univ. Napoli, Caserta, 2000.

\bibitem[Sh:F888]{Sh:F888}Saharon Shelah.
\newblock {Categoricity in $\lambda$ and a superlimit in $\lambda^+$}.

\bibitem[Sh 300a]{Sh:300a}Saharon Shelah.
\newblock {\em {Chapter I}}.

\bibitem[Sh 300b]{Sh:300b}Saharon Shelah.
\newblock {\em {Chapter II}}.

\bibitem[Sh 300c]{Sh:300c}Saharon Shelah.
\newblock {\em {Chapter III}}.

\bibitem[Sh 300d]{Sh:300d}Saharon Shelah.
\newblock {\em {Chapter IV}}.

\bibitem[Sh 300e]{Sh:300e}Saharon Shelah.
\newblock {\em {Chapter V}}.

\bibitem[Sh 300f]{Sh:300f}Saharon Shelah.
\newblock {\em {Chapter VI}}.

\bibitem[Sh 300g]{Sh:300g}Saharon Shelah.
\newblock {\em {Chapter VII}}.

\bibitem[Sh 322]{Sh:322}Saharon Shelah.
\newblock {Classification over a predicate}.
\newblock {\em {preprint}}.

\bibitem[Sh 840]{Sh:840}Saharon Shelah.
\newblock {Model theory without choice: Categoricity}.
\newblock {\em Journal of Symbolic Logic}, {\bf submitted}.
\newblock math.LO/0504196.

\bibitem[Sh:e]{Sh:e}Saharon Shelah.
\newblock {\em {Non--structure theory}}, accepted.
\newblock {Oxford University Press}.

\bibitem[Sh 838]{Sh:838}Saharon Shelah.
\newblock {Non-structure in $\lambda^{++}$ using instances of WGCH}.
\newblock {\em Preprint}.
\newblock 0808.3020.

\bibitem[Sh:F782]{Sh:F782}Saharon Shelah.
\newblock {On categorical a.e.c. II}.

\bibitem[Sh 800]{Sh:800}Saharon Shelah.
\newblock {On complicated models}.
\newblock {\em Preprint}.

\bibitem[Sh:F841]{Sh:F841}Saharon Shelah.
\newblock {On $h$-almost good $\lambda$-frames: More on [SH:838]}.

\bibitem[Sh:F735]{Sh:F735}Saharon Shelah.
\newblock {Revisiting 705}.

\bibitem[Sh 842]{Sh:842}Saharon Shelah.
\newblock {Solvability and Categoricity spectrum of a.e.c. with amalgamation}.
\newblock {\em Preprint}.

\bibitem[Sh 839]{Sh:839}Saharon Shelah.
\newblock {Stable Frames and weight}.
\newblock {\em Preprint}.

\bibitem[Sh 868]{Sh:868}Saharon Shelah.
\newblock {When first order $T$ has limit models}.
\newblock {\em Notre Dame Journal of Formal Logic}, {\bf submitted}.
\newblock math.LO/0603651.

\bibitem[Sh 1]{Sh:1}Saharon Shelah.
\newblock {Stable theories}.
\newblock {\em {Israel Journal of Mathematics}}, {\bf 7}:187--202, 1969.

\bibitem[Sh 3]{Sh:3}Saharon Shelah.
\newblock {Finite diagrams stable in power}.
\newblock {\em {Annals of Mathematical Logic}}, {\bf 2}:69--118, 1970.

\bibitem[Sh 10]{Sh:10}Saharon Shelah.
\newblock {Stability, the f.c.p., and superstability; model theoretic
  properties of formulas in first order theory}.
\newblock {\em {Annals of Mathematical Logic}}, {\bf 3}:271--362, 1971.

\bibitem[Sh 48]{Sh:48}Saharon Shelah.
\newblock {Categoricity in $\aleph _{1}$ of sentences in $L_{\omega
  _{1},\omega}(Q)$}.
\newblock {\em {Israel Journal of Mathematics}}, {\bf 20}:127--148, 1975.

\bibitem[Sh 54]{Sh:54}Saharon Shelah.
\newblock {The lazy model-theoretician's guide to stability}.
\newblock {\em {Logique et Analyse}}, {\bf 18}:241--308, 1975.

\bibitem[Sh 56]{Sh:56}Saharon Shelah.
\newblock {Refuting Ehrenfeucht conjecture on rigid models}.
\newblock {\em {Israel Journal of Mathematics}}, {\bf 25}:273--286, 1976.
\newblock A special volume, Proceedings of the Symposium in memory of A.
  Robinson, Yale, 1975.

\bibitem[Sh:a]{Sh:a}Saharon Shelah.
\newblock {\em {Classification theory and the number of nonisomorphic models}},
  volume~92 of {\em {Studies in Logic and the Foundations of Mathematics}}.
\newblock {North-Holland Publishing Co., Amsterdam-New York, xvi+544 pp,
  \$62.25}, 1978.

\bibitem[Sh:93]{Sh:93}Saharon Shelah.
\newblock {Simple unstable theories}.
\newblock {\em {Annals of Mathematical Logic}}, {\bf 19}:177--203, 1980.

\bibitem[Sh 82]{Sh:82}Saharon Shelah.
\newblock {Models with second order properties. III. Omitting types for
  $L(Q)$}.
\newblock {\em {Archiv fur Mathematische Logik und Grundlagenforschung}}, {\bf
  21}:1--11, 1981.

\bibitem[Sh 87a]{Sh:87a}Saharon Shelah.
\newblock {Classification theory for nonelementary classes, I. The number of
  uncountable models of $\psi \in L_{\omega _{1},\omega }$. Part A}.
\newblock {\em {Israel Journal of Mathematics}}, {\bf 46}:212--240, 1983.

\bibitem[Sh 87b]{Sh:87b}Saharon Shelah.
\newblock {Classification theory for nonelementary classes, I. The number of
  uncountable models of $\psi \in L_{\omega _{1},\omega }$. Part B}.
\newblock {\em {Israel Journal of Mathematics}}, {\bf 46}:241--273, 1983.

\bibitem[Sh 200]{Sh:200}Saharon Shelah.
\newblock {Classification of first order theories which have a structure
  theorem}.
\newblock {\em {American Mathematical Society. Bulletin. New Series}}, {\bf
  12}:227--232, 1985.

\bibitem[Sh 205]{Sh:205}Saharon Shelah.
\newblock {Monadic logic and Lowenheim numbers}.
\newblock {\em {Annals of Pure and Applied Logic}}, {\bf 28}:203--216, 1985.

\bibitem[Sh 197]{Sh:197}Saharon Shelah.
\newblock {Monadic logic: Hanf numbers}.
\newblock In {\em {Around classification theory of models}}, volume 1182 of
  {\em {Lecture Notes in Mathematics}}, pages 203--223. {Springer, Berlin},
  1986.

\bibitem[Sh 155]{Sh:155}Saharon Shelah.
\newblock {The spectrum problem. III. Universal theories}.
\newblock {\em {Israel Journal of Mathematics}}, {\bf 55}:229--256, 1986.

\bibitem[Sh 88]{Sh:88}Saharon Shelah.
\newblock {Classification of nonelementary classes. II. Abstract elementary
  classes}.
\newblock In {\em {Classification theory (Chicago, IL, 1985)}}, volume 1292 of
  {\em {Lecture Notes in Mathematics}}, pages 419--497. {Springer, Berlin},
  1987.
\newblock {Proceedings of the USA--Israel Conference on Classification Theory,
  Chicago, December 1985; ed. Baldwin, J.T.}

\bibitem[Sh 225]{Sh:225}Saharon Shelah.
\newblock {On the number of strongly $\aleph_ \epsilon$-saturated models of
  power $\lambda$}.
\newblock {\em {Annals of Pure and Applied Logic}}, {\bf 36}:279--287, 1987.
\newblock See also [Sh:225a].

\bibitem[Sh 300]{Sh:300}Saharon Shelah.
\newblock {Universal classes}.
\newblock In {\em {Classification theory (Chicago, IL, 1985)}}, volume 1292 of
  {\em {Lecture Notes in Mathematics}}, pages 264--418. {Springer, Berlin},
  1987.
\newblock {Proceedings of the USA--Israel Conference on Classification Theory,
  Chicago, December 1985; ed. Baldwin, J.T.}

\bibitem[Sh 225a]{Sh:225a}Saharon Shelah.
\newblock {Number of strongly $\aleph_ \epsilon$ saturated models---an
  addition}.
\newblock {\em {Annals of Pure and Applied Logic}}, {\bf 40}:89--91, 1988.

\bibitem[Sh:c]{Sh:c}Saharon Shelah.
\newblock {\em {Classification theory and the number of nonisomorphic models}},
  volume~92 of {\em {Studies in Logic and the Foundations of Mathematics}}.
\newblock {North-Holland Publishing Co., Amsterdam, xxxiv+705 pp}, 1990.

\bibitem[Sh 284c]{Sh:284c}Saharon Shelah.
\newblock {More on monadic logic. Part C. Monadically interpreting in stable
  unsuperstable ${\scr T}$ and the monadic theory of ${}^ \omega\lambda$}.
\newblock {\em {Israel Journal of Mathematics}}, {\bf 70}:353--364, 1990.

\bibitem[Sh 429]{Sh:429}Saharon Shelah.
\newblock {Multi-dimensionality}.
\newblock {\em {Israel Journal of Mathematics}}, {\bf 74}:281--288, 1991.

\bibitem[Sh 394]{Sh:394}Saharon Shelah.
\newblock {Categoricity for abstract classes with amalgamation}.
\newblock {\em {Annals of Pure and Applied Logic}}, {\bf 98}:261--294, 1999.
\newblock math.LO/9809197.

\bibitem[Sh 576]{Sh:576}Saharon Shelah.
\newblock {Categoricity of an abstract elementary class in two successive
  cardinals}.
\newblock {\em {Israel Journal of Mathematics}}, {\bf 126}:29--128, 2001.
\newblock math.LO/9805146.

\bibitem[Sh 472]{Sh:472}Saharon Shelah.
\newblock {Categoricity of Theories in $L_{\kappa^* \omega}$, when $\kappa^*$
  is a measurable cardinal. Part II}.
\newblock {\em {Fundamenta Mathematicae}}, {\bf 170}:165--196, 2001.
\newblock math.LO/9604241.

\bibitem[Sh 603]{Sh:603}Saharon Shelah.
\newblock {Few non minimal types and non-structure}.
\newblock In {\em {Proceedings of the 11 International Congress of Logic,
  Methodology and Philosophy of Science, Krakow August'99; In the Scope of
  Logic, Methodology and Philosophy of Science}}, volume~1, pages 29--53.
  Kluwer Academic Publishers, 2002.
\newblock math.LO/9906023.

\bibitem[ShUs 837]{ShUs:837}Saharon Shelah and Alex Usvyatsov.
\newblock {Model theoretic stability and categoricity for complete metric
  spaces}.
\newblock {\em Israel Journal of Mathematics}, {\bf submitted}.
\newblock math.LO/0612350.

\bibitem[ShVi 648]{ShVi:648}Saharon Shelah and Andr\'es Villaveces.
\newblock {Categoricity may fail late}.
\newblock {\em {Journal of Symbolic Logic}}, {\bf submitted}.
\newblock math.LO/0404258.

\bibitem[ShVi 635]{ShVi:635}Saharon Shelah and Andr\'es Villaveces.
\newblock {Toward Categoricity for Classes with no Maximal Models}.
\newblock {\em {Annals of Pure and Applied Logic}}, {\bf 97}:1--25, 1999.
\newblock math.LO/9707227.

\bibitem[Sh:E46]{Sh:E46}{Shelah, Saharon}.
\newblock {Categoricity of an abstract elementary class in two successive
  cardinals, revisited}.

\bibitem[Sh:E54]{Sh:E54}{Shelah, Saharon}.
\newblock {Comments to Universal Classes}.

\bibitem[Sh:E56]{Sh:E56}{Shelah, Saharon}.
\newblock {Density is at most the spread of the square}.
\newblock 0708.1984.

\bibitem[Sh:F709]{Sh:F709}{Shelah, Saharon}.
\newblock Good$^*$ $ \lambda $-frames.

\bibitem[Str76]{Str76}Jacques Stern.
\newblock {Some applications of model theory in Banach space theory}.
\newblock {\em Annals of Mathematical Logic}, {\bf 9}:49--121, 1976.

\bibitem[Va02]{Va02}Monica~M. VanDieren.
\newblock {\em {Categoricity and Stability in Abstract Elementary Classes}}.
\newblock PhD thesis, Carnegie Melon University, Pittsburgh, PA, 2002.

\bibitem[Zi0xa]{Zi0xa}B.I. Zilber.
\newblock Dimensions and homogeneity in mathematical structures.
\newblock preprint, 2000.

\bibitem[Zi0xb]{Zi0xb}B.I. Zilber.
\newblock A categoricity theorem for quasiminimal excellent classes.
\newblock preprint, 2002.

\end{thebibliography}

\enddocument